\documentclass{article}
\usepackage{graphicx}
\usepackage{amsmath}
\usepackage{amsfonts}
\usepackage{amssymb}
\usepackage{amsthm}
\usepackage{mathtools}
\usepackage{longtable}
\usepackage{multirow}
\usepackage{color}
\usepackage{url}    
\usepackage{comment}
\usepackage{tikz}
\usetikzlibrary{shapes.geometric, arrows.meta, positioning, fit, backgrounds, calc, shadows, decorations.pathreplacing}
\usepackage{caption}
\usepackage{adjustbox}
\usepackage{natbib} 
\usepackage{booktabs}
\usepackage{graphicx}
\usepackage{float}
\usepackage{placeins}
\usepackage{tabularx}
\usepackage{enumitem}

\newcommand{\pddl}[1]{\texttt{\small #1}}
\usepackage{geometry}
\DeclareMathOperator*{\argmax}{argmax}

\newcommand{\reg}{\mathrm{reg}}
\newcommand{\rl}{\mathrm{rl}}

\newcommand{\rev}[1]{{\color{blue}#1}}

\tikzset{
    robot/.pic={
        \fill[gray!30] (-0.3,-0.3) rectangle (0.3,0.3);
        \draw[thick] (-0.3,-0.3) rectangle (0.3,0.3);
        \draw[thick, fill=white] (0,0.4) circle (0.2); 
        \fill[gray] (-0.4,-0.3) circle (0.12); 
        \fill[gray] (0.4,-0.3) circle (0.12); 
        \draw[thick] (-0.15,0.5) -- (-0.25,0.7); 
    },
    treegraph/.pic={
        \fill[black] (0,0.25) circle (0.06);
        \draw[thick] (0,0.25) -- (-0.25,-0.1);
        \draw[thick] (0,0.25) -- (0.25,-0.1);
        \fill[black] (-0.25,-0.1) circle (0.06);
        \fill[black] (0.25,-0.1) circle (0.06);
        \draw[thick] (-0.25,-0.1) -- (-0.35,-0.35);
        \draw[thick] (-0.25,-0.1) -- (-0.15,-0.35);
    },
    hierarchy/.pic={
        \fill[black] (0,0.3) circle (0.08);
        \draw[thick, ->] (-0.05,0.2) -- (-0.25,-0.05);
        \draw[thick, ->] (0,0.2) -- (0,-0.05);
        \draw[thick, ->] (0.05,0.2) -- (0.25,-0.05);
        \fill[gray] (-0.3,-0.15) circle (0.08);
        \fill[gray] (0,-0.15) circle (0.08);
        \fill[gray] (0.3,-0.15) circle (0.08);
    },
    miniplot/.pic={
        \draw[->, thin] (-0.3,-0.3) -- (0.3,-0.3);
        \draw[->, thin] (-0.3,-0.3) -- (-0.3,0.3);
        \draw[thick, blue] (-0.3,-0.15) .. controls (-0.1,0.2) .. (0.3,0.05);
    },
    clockface/.pic={
        \draw[thick, fill=white] (0,0) circle (0.3);
        \draw[thick] (0,0) -- (0,0.2);
        \draw[thick] (0,0) -- (0.15,-0.1);
        \foreach \angle in {0,90,180,270} \draw (\angle:0.25) -- (\angle:0.3);
    },
    smiley/.pic={
        \fill[yellow!80!orange, draw=black, thick] (0,0) circle (0.3);
        \fill[black] (-0.1, 0.05) circle (0.04);
        \fill[black] (0.1, 0.05) circle (0.04);
        \draw[thick] (-0.15, -0.1) .. controls (-0.05, -0.2) and (0.05, -0.2) .. (0.15, -0.1);
    },
    crash/.pic={
        \draw[thick, fill=red!80!yellow, line join=round]
            (0,0.4) -- (0.08,0.15) -- (0.35,0.25) -- (0.15,0.08) -- (0.4,0) --
            (0.15,-0.08) -- (0.25,-0.35) -- (0.08,-0.15) -- (0,-0.4) --
            (-0.08,-0.15) -- (-0.35,-0.25) -- (-0.15,-0.08) -- (-0.4,0) --
            (-0.15,0.08) -- (-0.25,0.35) -- (-0.08,0.15) -- cycle;
    }
} 

\title{Mission-Aligned Learning-Informed Control of Autonomous Systems: Formulation and Foundations}
\author{Vyacheslav Kungurtsev, Alessandro Di Frenna,  \\ Gustav {\v S}{\' i}r,
Monicah Cherop Naibei, \\ Haozhe Tian, Homayoun Hamedmoghadam \\ Akhil Anand, Sebastien Gros,   \footnote{VK (corresponding author: \protect\url{vyacheslav.kungurtsev@fel.cvut.cz}) and GS are with Czech Technical University in Prague. HT and HH are with Imperial College London. AA and SG are with the Norwegian University of Science and Technology. VK and GS acknowledge funding support to the Prosocial AI STEM initiative through the Institute for Advanced Conscious Studies in Los Angeles, CA, USA, as well as the Czech National Science Foundation under Project 24-11664S. HT acknowledges the funding by UK Research and Innovation (UKRI) under the UK government’s Horizon Europe funding guarantee (grant number 101084642), and support from IBM through an IBM PhD Fellowship. Authors acknowledge financial support from Imperial College London through an Imperial College Research Fellowship grant awarded to HH.}}
\date{March 2026}

\begin{document}

\maketitle

\begin{abstract}
Research, innovation and practical capital investment have been increasing rapidly toward the realization of autonomous physical agents. This includes industrial and service robots, unmanned aerial vehicles, embedded control devices, and a number of other realizations of cybernetic/mechatronic implementations of intelligent autonomous devices. In this paper, we consider a stylized version of robotic care, which would normally involve a two-level Reinforcement Learning procedure that trains a policy for both lower level physical movement decisions as well as higher level conceptual tasks and their sub-components. In order to deliver greater safety and reliability in the system, we present the general formulation of this as a two-level \rev{expert defined} optimization scheme which incorporates control \rev{with physics modeling of the system and environment} at the lower level, and classical planning \rev{with complex task logic} at the higher level, integrated with a capacity for learning. This synergistic integration of multiple methodologies---control, classical planning, and RL---presents an opportunity for greater insight for algorithm development, leading to more efficient and reliable performance. Here, the notion of reliability pertains to physical safety and interpretability into an otherwise black box operation of autonomous agents, concerning users and regulators. This work presents \rev{the principled integration of the two expert-defined domains. The methodology is validated in a simulation of robotic care aid, with complex medication and meal preparation tasks at the upper level and locomotion at the lower level. The validation confirms the framework's adaptability to human preferences and symbolic to continuous feedback}.  
\end{abstract}

\section{Introduction}

Consider a setting in which a robot is a permanent aide to a movement and health-impaired human. This robot is supposed to provide for the individual's needs, while generally attempting to benefit their health and well-being. It should also ensure that the individual is safe, and that any unexpected behavior should trigger immediately notifying a human specialist with the most critical information. There are various possible relevant objects in the environment, and the robot either moves next to the human to assist movement or retrieves the object for the human. There are a set of regular expected activities, eating, grooming, exercising, scribing dictation, etc. that the robot performs in accordance with: the desires of the human, the human's well-being as according to the responsible health authorities, the physical safety of the human, and any regulations. The robot can communicate and should suggest any activity on the schedule to the human, who can refuse or stop the activity at any time. Warnings can be issued, e.g., ``if we move away now your living room will have this pile and may look messy,'' but ultimately the human decides its behavior. 

In service of this task, the robot is able to lean to roll forward and back as well as turn in full rotation, and rotate its arm along three axes. Objects are picked up with grasping, or relevant task-oriented objects are magnetized. Objects in the environment define collisions to be avoided, which requires including the speed and position in the state space. The robot has a propulsion system, which is its forward control, and a braking mechanism. The manufacturing of this propulsion system involves some new features which are not system-identified to full certainty, and there is a small but noticeable nonlinear gyroscopic effect on its motion. This effect will have to be learned in the course of its operation. 

This problem requires the incorporation of learning, which can be performed with standard RL procedures, along with tokenized interpretable actions and tasks for trustworthiness, as well as physical mechatronic control of the robot's actuators. This combination is famously handled in the present day with Vision-Language-Action models, e.g.~\cite{ma2024survey}. However, long term autonomous robotics for sensitive domains with potential critical risks, such as physical care, is still a long distance, in research and development output, from trustworthy operation at scale~\cite{christoforou2020upcoming,pandy2025ai}. The fundamental problems associated with Deep RL for real world physical application is seminally summarized in \cite{dulac2019challenges}. The problems observed in practice include a) poor out-of-distribution generalization, b) lack of formal guarantees of accuracy and performance, essential for trust in end users and regulators, and c) strict adherence to reasonable protocols of safety. These issues are fundamental: recall that RL is, ultimately, a statistical model, one that does not structurally incorporate the vast knowledge base of the natural sciences. As such, the mathematical representation of the model lacks the mechanistic intuition appropriate for understanding cause and effect~\cite{kungurtsev2025cause} and thus potential consequences of, e.g., hazardous actions.

Inspired by the recent flurry of activity of Reinforcement Learning-based Model Predictive Control (RLMPC)~\cite{reiter2025synthesis}, the present paper develops a methodology for incorporating domain knowledge in order to enhance the efficiency and reliability of autonomous deployment. Indeed, it has been observed in several applications that combining Reinforcement Learning with Model Predictive Control can achieve superior performance than both MPC and RL alone, see e.g.~\cite{reiter2025synthesis}. 

We present a formalism that defines a bilevel optimization problem incorporating Reinforcement Learning (RL)-assisted Model Predictive Control in the bottom layer, and RL-assisted classical planning at the higher level. 
The overarching architecture of this proposed framework is depicted in Figure~\ref{fig:methodology_overview}.

\begin{figure}[t]
\centering
\begin{adjustbox}{width=\textwidth}
\begin{tikzpicture}[
  node distance=0.6cm and 0.8cm,
  richblock/.style={
      rectangle, draw, thick, fill=white, drop shadow,
      text width=3.2cm, text centered, rounded corners, 
      minimum height=2.5cm, font=\footnotesize
  },
  arrow/.style={thick, -{Stealth[scale=1.2]}},
  dashedarrow/.style={thick, dashed, -{Stealth[scale=1.2]}}
]
  \node (scheduler) [richblock, fill=blue!5] {
      \textbf{Scheduler} \\ \vspace{0.1cm}
      \tikz\pic[scale=0.8]{clockface}; \\ \tiny Slow Scale \\ Max $\mathbf{F}(\{F_T\})$
  };
  \node (planner) [richblock, fill=orange!5, right=of scheduler] {
      \textbf{Planner} \\ \vspace{0.1cm}
      \tikz\pic[scale=0.8]{treegraph}; \\ \tiny Discrete Logic \\ $a_t \in \mathcal{A}, s_t \in \mathcal{S}$
  };
  \node (mpc) [richblock, fill=green!5, right=1.2cm of planner] {
      \textbf{RL-MPC} \\ \vspace{0.1cm}
      \tikz\pic[scale=0.8]{miniplot}; \\ \tiny Continuous \\ $v(\tau), z(\tau)$
  };
  \node (env) [richblock, fill=gray!10, right=of mpc] {
      \textbf{Agent inside Environment} \\ \vspace{0.1cm}
      \tikz\pic[scale=0.8]{robot}; \\ \tiny Physical World \\ $\hat{z}_{t+1}, \hat{R}_{t+1}$
  };

  \draw [arrow] (scheduler) -- node[midway, above, font=\tiny, align=center] {Task\\$T$} (planner);
  \draw [arrow] (planner) -- node[midway, above, font=\tiny, align=center] {Action\\$a_t, p(a_t)$} (mpc);
  \draw [arrow] (mpc) -- node[midway, above, font=\tiny] {$v(\tau)$} (env);
  
  \draw [dashedarrow] (env.south) -- ++(0,-0.5) -| node[pos=0.25, above, font=\tiny] {Mapped State $\hat{s}_{t+1}$} (planner.south);
  \draw [dashedarrow] (env.south) -- ++(0,-1.0) -| node[pos=0.25, above, font=\tiny] {Task Eval $F_T$} (scheduler.south);
\end{tikzpicture}
\end{adjustbox}
\caption{High-level methodology overview. A Scheduler (left) drives a Planner and MPC, ultimately controlling the physical Agent (right), with feedback loops returning state and reward information for respective updates.}
\label{fig:methodology_overview}
\end{figure}

Human-aligned and safe autonomous robotics is a technology of significant interest to many researchers, companies, and governments, with population aging creating many potential labor shortages for assistance care and health for the elderly, the potential for automated manufacturing, and other reasons. As it is also a representative example of the conceptual components of such an operation, it will be a running example throughout this line of work, although the same principles apply to similar domains such as automated manufacturing, self-driving vehicles and other embodied AI. In this work we focus on a simplified stylistic model of long run autonomous robotic human care and focus on developing the foundations of the mathematical formulation for the appropriate optimization and data driven sequential decision making framework. We present a formal schema to define the autonomous decision procedures of an agent using multiple relevant disciplines associated with mathematical engineering.

\paragraph{AI and Autonomous Systems}

Designing and managing systems towards optimality or looser satisfaction criteria over time is a generic cross-disciplinary problem that has attracted attention from numerous communities. In the case wherein the model of the states and actions can be well represented with either a finite set or fairly low dimensional and bounded continuous domain while the system dynamics and the utility of different circumstances is unknown, Reinforcement Learning (RL) has risen to prominence as the preeminent methodology~\rev{\cite{sutton2018reinforcement,bertsekas1996neuro}}. RL, like Classification (e.g., image recognition), and Generative Learning (e.g., Large Language Models), has been a central component in the accelerating, groundbreaking developments of AI that have characterized the technological landscape for the past decade. Neural networks leverage the raw computational power of highly parallel GPU hardware permitting training of deep models with millions to billions of trainable parameters. This forms highly accurate statistical models of objects, processes and phenomena when there are vast quantities of data present or available. 

In recent years, a variety of approaches for opening the deep and opaque models developed by AI training have appeared aiming to improve their interpretability and reliability. In the case of image recognition and LLMs, these concerns are largely a matter of ensuring an agreeable user experience. However, RL involves making decisions. In the context of games of pure entertainment, this relative lack of criticality of risks persists. However, when RL agents make decisions in contexts with real-world consequences, a notable degree of justified concern appears. An opaque black-box model simply cannot provide formal guarantees of safety or performance reliability. Rather, it is simply a model that efficiently provides a most likely answer, in some formal probabilistic sense. Low probability but catastrophic events could happen and the risk of such rare tail events is difficult to quantify -- naturally given low sample counts for rare events, it is challenging to use purely empirical methods to predict and forecasts such risks and simulate contingency plans. \rev{Moreover, RL pipelines are associated with table lookup based finite state and action spaces or neural network based parametrization. In many circumstances, the autonomous decisions are meant to exhibit more complex discrete structure, as defined by, for instance, rules in propositional or first order logic. This motivates hybrid decision-making stacks that complement learned policies with explicit high-level task representations from classical planning~\cite{ghallab2004automated,lavalle2006planning}}.


\paragraph{Contributions}

The general premise of this work is as follows: it is currently difficult if not impossible to ensure that RL can systematically avoid or minimize the risk of damaging living beings or critical agents and infrastructure via action or instability. Probabilistic learning techniques, importance sampling in unsafe regions, human in the loop, and other strategies that remain within the RL paradigm~\cite{garcia2015comprehensive} may present clever mechanisms to dramatically lower the probability of said risks, but they are fundamentally insufficient. Comprehensive guarantees of safe and stable operation must arise from a paradigm that is separate from probabilistic inductive learning.

In particular, crisp hard science-based models of system dynamics and foreseeable risks (e.g. of collision) are the standard, and a comprehensively reliable approach to endow RL systems with true guarantees of safe and stable operation. Said guarantees would significantly expand the set of domains and circumstances to which RL can be applied, and at the same time provide a transparent and interpretable means by which regulators can certify the safety of an RL system/agent. This will be done by embedding Model Predictive Control within the lower layer. Moreover, at the highest level of abstraction as to the agent's operation, the decisions must be interpretable and associated with the mission and tasks associated with the AI. These should make logical and semantic sense to a user, and as far as the general agreement and sense of sequences of higher level operations. This will be done by embedding classical planning into the top layer.

In this paper we present a set of technical procedures from the methodological fields of process control and planning that incorporate structural assurances into deep neural networks in a systematic way. Along with addressing the aforementioned concerns of interpretability and reliability, these sets of computational techniques significantly improve the quality and precision of the RL model and enable magnitudes of greater statistical power for any available datasets used to train the optimal decision model.

In summary, the Contributions of this manuscript include:
\begin{itemize}
 \item A bilevel integrated data-driven optimization scheme for operationalizing high level interpretable token-based decisions together with low level physical control behavior. \rev{Specifically the rich discrete ontologies and expert-defined rules associated with verifiation, safety, human-alignment and interpretability are enabled with a upper level Classical Planning layer. Meanwhile the physics-based modeling and world representation define a nonlinear-programming based Model Predictive Control operation at the lower level.}
 \item A novel methodology for enriching Classical Scheduling and Planning with actions defining the operation of control problems and real-time streaming data adjusting the components Planning ontology. \rev{By formally assigning each action to the triggering of a formulation and real-time algorithmic computation of Model Predictive Control, this enables decision-integration across high level discrete and level continuous action spaces.}
 \item A novel methodology of complete end-to-end Planning-Control differentiation through the use of fuzzification operators, \rev{Fuzzy maps are used to provide information from the state estimation and identification performed alongside real-time control to informative discrete quantities for the Planning layer to update its high level conception in accordance with changes and discrepancies in models of the system and environment.}
 \item A systems engineering formalism that facilitates crisp guarantees of trustworthy interpretable actions and safe control behavior simultaneously with learning to facilitate adaptation to an uncertain and changing environment. \rev{While Deep Reinforcement Learning is eschewed as an end-to-end operator, its techniques can still be used as far as obtaining real-time information in the form of rewards that update the cost functions of Planner in the same methodology as in Q-Learning.}
 \item Comprehensive validation through a standing example of Robotic Care throughout the manuscript together with a detailed numerical illustration. \rev{The simulation incorporates discrete planning for complex tasks of medication delivery and food preparation and delivery together with a robot's continuous locomotion through the environment, in a manner that is adaptive to the preferences of the human cared for. }
\end{itemize}

\rev{The rest of this paper proceeds as follows. Section~\ref{sec:rl} describes the predominant baseline for complex autonomous systems of hierarchical RL. Section~\ref{sec:MPC} reviews Model Predictive Control - which facilitates expert physics based modeling into control, and recent developments incorporating learning into its methodology. Section~\ref{sec:planning} reviews Classical Planning and Scheduling. Section~\ref{sec:fuzzy} reviews fuzzy maps, which enable integration between discrete and continuous quantities. Section~\ref{sec:intergrating} presents the framework of our contribution by which symbolic expert systems associated with planning are operationally intertwined with physics-based MPC. Section~\ref{sec:experiments} describes the illustrative simulation we developed for this work modeling robotic care. Section~\ref{sec:results} describes the results of our simulation. This is followed by a general Discussion Section~\ref{sec:discussion} that presents our considerations for desiderata in extensions and developments for the implementation of our framework in the real world as well as current limitations in doing so. We conclude in Section~\ref{sec:conclusion}.}

\section{Hierarchical Reinforcement Learning}\label{sec:rl}
In this Section we provide the background for and the full formulation, as well as a description thereof, of a Bilevel Reinforcement Learning problem as it could be applied to the practical consideration above regarding robotic human care. To this end, we present a summary of the main principles of RL, its extension to continuous action spaces, and bilevel and generally hierarchical RL problems and solution algorithms.

\subsection{Reinforcement Learning - Background}

Formally, a standard RL problem (two prominent texts being~\cite{sutton2018reinforcement,bertsekas1996neuro}) is defined to consist of a set of states $\mathcal{X}$ and a set of actions $\mathcal{A}$. The objective is to learn a policy $\pi:\mathcal{X}\to \mathcal{M}(\mathcal{A})$, where $\mathcal{M}(\mathcal{A})$ denotes the set of multinomial distributions over the action set $\mathcal{A}$. The state evolves in time through an unknown transition probability tensor $\mathbf{P}$ with components $P(S_{t+1} = s'\vert S_t = s, A_t= a)$. The reward $R(s,s', a)$ defines some utility the agent receives upon that particular transition to $s'$. In RL, techniques to learn the Value function $V^{\pi}(s)$ and Action-Value function $Q^{\pi}(s, a)$ are iterated together with establishing the policy $\pi$ in either performing the most rewarding action $a$ or randomly explore and learn more about the system. 

RL searches for the optimal decision-making policy by interacting with a system modeled as a Markov Decision Process (MDP).
Due to its data-driven nature, RL is particularly powerful in complex or stochastic systems, for which explicit models are difficult to construct~\cite{recht2019tour}. RL aims to learn a policy $\pi: \mathcal{S}\times\mathcal{A}\rightarrow[0,1]$ that minimizes the cumulative reward, or return, defined as:
\begin{align}
\begin{split}
    \eta(\pi) =  \mathbb{E}_{s(0), a(0), \dots}\left[\sum_{t=0}^\infty\gamma^t R(s(t), a(t))\right],
\end{split}
\label{eq:eta}
\end{align}
where $s(0)\sim \rho(s(0))$, $a(t)\sim\pi(a(t)| s(t))$, and $s(t+1)\sim P(s(t+1)|s(t),a(t))$. At each state, the expected cumulative reward given a policy $\pi$ can be estimated via the value function $V^\pi(s(t))$ and the state-action value function $Q^\pi(s(t), a(t))$, defined below:
\begin{align}
\begin{split}
    &V^\pi(s(t))=\mathbb{E}_{a(t), s(t+1),\dots}\left[\sum_{i=0}^\infty\gamma^{t+i}R(s(t), a(t))\right], \\
    &Q^\pi(s(t), a(t))=\mathbb{E}_{s(t+1), \dots}\left[\sum_{i=0}^\infty\gamma^{t+i}R(s(t), a(t))\right],
\end{split}
\label{eq:vf}
\end{align}
where $a(t)\sim\pi(a(t)| a(t))$, and $s(t+1)\sim f(s(t+1)|s(t),a(t))$. For Eq.~\ref{eq:eta}\&\ref{eq:vf}, it is seen that $\eta(\pi)=\mathbb{E}_{s(0)\sim \rho(s(0))}\left[V^\pi(s(0))\right]$. With the system modeled as an MDP, state transition depends only on the current state and action, thus, $V^\pi(s(t))$ and $Q^\pi(s(t), a(t))$ can be recursively estimated from observed interactions $(s(t), a(t), s(t+1), R(s(t), a(t)))$ using the Bellman equations:
\begin{align}
\begin{split}
    &V^\pi(s(t)) = \mathbb{E}_{a(t), s(t+1)}\left[R(s(t), a(t)) + \gamma V^\pi(s(t+1))\right], \\
    &Q^\pi(s(t), a(t)) = \mathbb{E}_{s(t+1)}\left[R(s(t), a(t)) + \gamma V^\pi(s(t+1))\right], \\
\end{split}
\label{eq:value-bellman}
\end{align}
where  $s(t)\sim\pi(s(t)| a(t))$ and $s(t+1)\sim P(s(t+1)|s(t),a(t))$. The Bellman equations define a recursive update rule for the expected return of state-action pairs, enabling iterative refinement based on observed interactions.

\subsection{Reinforcement Learning with Function Approximations}

The update rule in Eq.~\ref{eq:value-bellman} treats states $\mathcal{S}$ and state-action pairs $\mathcal{S}\times\mathcal{A}$ as table entries that need to be estimated. Consequently, methods that directly estimate $ V^\pi(s(t)) $ and $ Q^\pi(s(t), a(t))$ are referred to as tabular RL. However, in many real-world applications, the state space $ \mathcal{S} $ and action space $ \mathcal{A} $ are extremely large or continuous, and it becomes infeasible to apply tabular RL due to the curse of dimensionality.

To handle large or continuous state and action spaces, modern RL methods employ function approximation for value functions. Specifically, the policy $\pi$, state-value function $V^\pi$, and the action-value function $Q^\pi$ are approximated using parameterized functions, $\pi_\theta$, $V_\phi$, and $Q_\psi$, with trainable parameters $\theta$, $\phi$ and $\psi$, respectively. This parameterization facilitates generalization to unobserved states and state-action pairs, allowing learning across $\mathcal{S}$ and $\mathcal{A}$ without explicit enumeration. The function approximation also enables a series of RL algorithms termed the policy gradient methods, which directly optimize the policy $\pi$ to maximize the objective in Eq. \ref{eq:eta} \cite{sutton1999policy}.

Function approximation (in RL), has a number of challenges, such as overfitting and reduced stability during training \cite{packer2018assessing, mnih2013playing}. In the tabular setting, where state and action spaces are discrete, each policy update in RL is guaranteed to improve performance \cite{watkins1992q, kakade2002approximately}. However, when function approximation is used, such guarantees become difficult to establish. To address stability issues in training RL with deep neural networks as function approximators, various approaches have been proposed. For example, the Replay Buffer \cite{mnih2013playing} helps stabilize training by sampling batches of historical interactions. Trust-region policy optimization (TRPO) \cite{schulman2015trust} employs a second-order method to constrain policy updates and prevent performance degradation. Clipped double Q-learning \cite{fujimoto2018addressing} mitigates overestimation bias in $Q_\psi$ learning by leveraging target Q-networks. Still, the performance of RL has been outstanding for many settings - by using deep neural networks as function approximators and the recently proposed stabilizing techniques, RL can now address many problems that were previously intractable \cite{levine2016end, zhu2017target, degrave2022magnetic}.

\subsection{Bilevel Reinforcement Learning Formulation}

Hierarchical, or Multilevel (in generalizing the classic Bilevel case) optimization concerns the solution of problems whereby the criteria exhibit a dependent layer structure. That is, the reward function for the top layer depends on the value of another variable which is in turn the solution of a different optimization, one which in turn depends on the operations of the earlier decision variables. See, e.g.~\cite{migdalas2013multilevel} for a collection of representative works. When all of the levels define convex optimization problems, reliable single loop algorithms can be assured through careful step level-weighing strategies~\cite{shafiei2024trilevel}. 

In the case wherein the cost/reward functions are nonconvex, convergence guarantees are less prevalent and firm. Standard Stochastic Approximation multilevel extensions are described in~\cite{dereich2019general}. Hierarchical Stochastic Approximation has also been considered for (nonsmooth) risk measures~\cite{crepey2024adaptive}. Typically, these algorithms require a careful tuning of stepsize annealing, with a hierarchy of decay rates for stepsizes associated with each level. However, see a recent single-timescale analysis in~\cite{ghadimi2020single} and~\cite{shen2022single}.

However, in consideration of Reinforcement Learning, there is an additional bilevel structure at each layer. In particular, the policy is what is of interest, but the policy depends on the Value function associated with each state, which itself must be learned. Since this dependency is directional, two timescales is the standard approach, e.g.~\cite{zeng2024two}. For an application in multi-stage stochastic programming, which can be considered an approximation to DP, with Stochastic Approximation, see~\cite{lan2021dynamic}. The two timescale approach lends itself to theoretical analysis of the learning procedure using Stochastic Differential Equations, due to the work in e.g.~\cite{yin2012continuous,zhang2017optimal,yin2009asymptotic,nguyen2024second}. Extensions to the hierarchical case, that is, considering both the hierarchy of DP problems together with the necessity of two time scale methods for RL, will be necessary for implementing this work.

For Reinforcement Learning specifically~\cite{kreidieh2019inter} considers a boss-worker hierarchical RL with derivation of a cooperative policy gradient for the boss. The work~\cite{zhu2022hierarchical} defines navigation to avoid obstacles at a lower level, and goal-driven movement at the upper level, for a robotics RL engine. The captioning of the image with a language layer and a RL visual layer is described in~\cite{xu2019multi}

To define the two tiered problem of managing a robotic aid to a human, we use a discrete action space for the set of chosen tasks, such as to clean, feed, etc. for the upper problem and a continuous action space for the optimal control defining navigation and operation at the lower level. We will use deterministic policies, meaning that there is one unique optimal action chosen upon the presence of any given state. Of course, algorithmically, we may still randomize the action chosen for exploration purposes.  

The Upper Level RL Problem aims to learn a deterministic policy $\pi: \mathcal{X}\to \mathcal{A}$ that minimizes the cumulative stage cost defined as:
\begin{align*}
    \eta(\pi) = \mathbb{E}_{x(0,x(1),\cdots)}
    \left[\sum_{t=0}^\infty\gamma^t F\left(x(t), a(t),\left\{z(\tau),v(\tau\}\right\}_{t-1<\tau\le t}(x(t-1), a(t-1))\right)\right],
\end{align*}
where 
\[
x(0)\sim \rho(x(0)), a(t)=\pi(s(t))
\]
and 
\[
x(t+1)\sim f(x(t+1)|x(t),a(t),\{z(\tau),v(\tau)\}_{t<\tau\le t+1})
\]
and $z(\tau)$ is a time varying continuous state space encoding the location, velocity and other canonical coordinates. Note that this $f$ encodes a probability distribution over a set of states $x(t+1)\in\mathcal{X}$. The inner problem control $v(\tau)$ is the time varying control action driving a deterministic process $z(\tau+h)=g(z(\tau),v(\tau))$ where $h\ll 1$ is the fast time scale. Note that this process depends on the state and action of the upper level chosen at the previous time step. 

The continuous control action $v(\tau)$ solves a lower level problem:
\begin{align}
\begin{split}
    \zeta(\pi;\eta,t) = \mathbb{E}_{x(t+1)}
    \left[\sum_{\tau=t,t+h,\cdots}^{t+1}\beta^{\tau} G\left(z(\tau), v(\tau),x(t+1),a(t+1),x(t),a(t))\right)\right],
\end{split}
\label{eq:etabil2}
\end{align}
where $G$ is a functional that includes either $F$ and other local physical criteria, or tracking a solution that optimizes the top level functional. 
This interaction between the slow-timescale discrete planner and the fast-timescale continuous controller is illustrated in Figure \ref{fig:bilevel_architecture}.

\begin{figure}[t]
\centering
\resizebox{\textwidth}{!}{%
\begin{tikzpicture}[
  box/.style={rectangle, draw, thick, fill=white, drop shadow, minimum width=4.5cm, minimum height=4.2cm, align=center, font=\footnotesize},
  arr/.style={-{Stealth[scale=1.2]}, thick}
]

  \node[box, fill=orange!10] (upper) at (0,0) {
      \textbf{Upper Level RL} \\ \textit{(Discrete, Slow $t$)} \\[0.1cm]
      \tikz{\pic[scale=0.8]{hierarchy};} \\[0.2cm]
      \textbf{State $x(t) \in \mathcal{X}$:} Discrete Context \\
      \textbf{Action $a(t) \in \mathcal{A}$:} Feed, Clean, ... \\[0.2cm]
      Policy $\pi: \mathcal{X} \to \mathcal{A}$ \\
      $\min \mathbb{E} \left[\sum \gamma^t F(\dots)\right]$
  };
  
  \node[box, fill=cyan!10] (lower) at (6.5,0) {
      \textbf{Lower Level Control} \\ \textit{(Continuous, Fast $\tau$)} \\[0.1cm]
      \tikz{\pic[scale=0.8]{miniplot};} \\[0.2cm]
      \textbf{State $z(\tau)$:} Location, Velocity \\
      \textbf{Control $v(\tau)$:} Motor Inputs \\[0.2cm]
      Policy $\zeta: z(\tau) \to v(\tau)$ \\
      $\min \mathbb{E} \left[\sum \beta^\tau G(\dots)\right]$
  };
  
  \node[box, fill=gray!5, minimum width=5.0cm] (dyn) at (13.2,0) {
      \textbf{Environment \& Rewards} \\
      \textcolor{gray!80!black}{Evaluate Costs $F$ and $G$} \\[0.2cm]
      \rule{0pt}{2.2cm} 
  };
  
  \pic[scale=0.7] (rob) at (11.5, -0.5) {robot};
  
  \pic[scale=0.8] (happy) at (13.5, 0.2) {smiley};
  \node[right, font=\footnotesize, text=green!60!black, align=left] at (13.8, 0.2) {\textbf{+ Reward} \\ \scriptsize Accomplished};
  
  \pic[scale=0.8] (sad) at (13.5, -1.0) {crash};
  \node[right, font=\footnotesize, text=red, align=left] at (13.8, -1.0) {\textbf{- Penalty} \\ \scriptsize ~~~Collision};
  
  \draw[->, thick, dashed, draw=green!60!black] (11.9, -0.2) .. controls (12.3, 0.2) and (12.7, 0.2) .. (13.1, 0.2);
  \draw[->, thick, dashed, draw=red] (12, -0.7) .. controls (12.3, -1.0) and (12.7, -1.0) .. (13.1, -1.0);

  \draw[arr] (upper) -- node[midway, above, font=\footnotesize] {Action $a(t)$} (lower);
  \draw[arr] (lower) -- node[midway, above, font=\footnotesize] {Control $v(\tau)$} (dyn);
  
  \draw[arr, dashed] (dyn.north) -- ++(0,0.6) -| node[pos=0.25, above, font=\footnotesize] {Discrete State $x(t+1)$} (upper.north);
  \draw[arr, dashed] (dyn.south) -- ++(0,-0.6) -| node[pos=0.25, below, font=\footnotesize] {Continuous State $z(\tau+h)$} (lower.south);

\end{tikzpicture}%
}
\caption{\rev{Hierarchical Bilevel RL Architecture: Visualizing the decoupling of timescales. The upper level policy assigns discrete abstract actions (e.g., feeding a human) to the lower level. The lower level evaluates continuous state feedback at a high frequency to generate physical control outputs. The resulting environmental trajectory determines the reward signals evaluated by the cost functions $F$ and $G$.}}
\label{fig:bilevel_architecture}
\end{figure}

In this case the upper level value function $V^\pi(x(t))$ becomes:
\begin{align}
\begin{split}
    V^\pi(x(t))=\mathbb{E}_{x(0),x(1),\cdots}\left[\sum_{i=0}^\infty\gamma^{t+i}F(x(t), a(t),\left\{z(\tau),v(\tau\}\right\}_{t-1<\tau\le t}(x(t-1), a(t-1)))\right],\\ a(t)=\pi(x(t))
\end{split}
\label{eq:vfupp}
\end{align}
However, this would assume $\left\{z(\tau),v(\tau\}\right\}_{t-1<\tau\le t}(x(t-1), a(t-1))$ already optimizes the Value function, that is satisfies the Bellman equations
\begin{align}
\begin{split}
   V^{\zeta}(z(\tau)) = \mathbb{E}_{ x(t+1)}\left[\beta^{\tau} G\left(z(\tau), v(\tau),x(t),a(t),x(t+1),a(t+1))\right) + \beta V^{\zeta}(z(\tau+h))\right],
\end{split}
\label{eq:belldown}
\end{align}
for all $t$. Observe the inherent coupling, the upper and lower level solutions depend on each other, introducing a couple Dynamic Programming Fixed Point Problem. However, one can consider that in real-time operation, one can predict the state for future control problems, and so always have precomputed solutions one is consistently updating. 


\section{Safe RL with Model Predictive Control}\label{sec:MPC}

A significant challenge in applying RL to real-world problems is safety. Although approaches based on Constrained Markov Decision Processes (CMDPs) \cite{altman2021constrained, achiam2017constrained} allow minimization of the expected cumulative stage cost while following a set of safety constraints, they often require the exploration of unsafe trajectories during training.
To avoid safety violations during training in particular (and to enhance learning efficiency in general), an emerging solution is the integration of domain knowledge or system constraints into the RL learning process \cite{garcia2015comprehensive}. This approach has been shown to outperform data-driven RL methods with soft penalties discouraging constraint violation, e.g.~\cite{ray2019benchmarking}.

The relevant domain for ensuring safe operation  according to the laws of physics is \emph{Model Predictive Control} (MPC). Indeed, the primary historical distinction between the two fields, model predictive control and reinforcement learning, developed under their separate domain native considerations, is the use of explicit background knowledge informed from known scientific principles, on the one hand, and (often very) large sample experiments querying a black box comprehensive roll-out simulation of the entire procedure. In practice, one is often in between the ideal scenario for MPC and RL, that is, a partially complete physical understanding of some of the mechanics is available. As such, part of the system dynamics are known, and part have to be ascertained with data.

\subsection{Model Predictive Control - Background}

Formally, MPC consists of repeated solution of an optimal control problem (OCP) defined over some finite horizon.

With MPC, one solves the discretized nonlinear optimization problem defining an approximation of the near future system behavior, applies the first control, lets the system evolve, measures the state again, and solves a new OCP with the horizon shifted one step forward. Note two important distinctions from RL: 1) the lack of uncertain quantities or random variables in the model, and 2) the presence of functions defining constraints. These constraints can enforce safety, e.g., by indicating the geometry of the environment as far the presence of objects an agent can collide with. 

MPC provides a policy by solving an OCP at each discrete time instant based on the current system state $s$, on a finite receding horizon. The MPC problem can be cast as,
\begin{subequations}
\label{eq:MPC0} 
\begin{align}
\min_{x,u}&\quad V_T \left(x_T\right) + \sum_{t=0}^{T-1}\, F\left(x_t,u_t\right)\,, \label{eq:cost0}\\
\mathrm{s.t.} &\quad x_{t+1} = f\left(x_t,u_t\right),\quad x_0 = s\,, \label{eq:dyn0}\\
&\quad h\left(x_t,u_t\right) \leq 0,\quad u_t \in \mathcal{A}\,, \quad  \label{eq:const0} 
\end{align}
\end{subequations}
where $T$ is the prediction horizon, $F$ is the stage cost, $V_T$ is the terminal cost, $f$ is the dynamics and $h$ is the inequality constraint. Problem \ref{eq:MPC0} produces a complete profile of control inputs $u^\star = \{u_0^\star,\ldots, u_{N-1}^\star\}$ and corresponding state predictions $x^\star= \{x_0^\star,\ldots, x_{N}^\star\}$. However, only the first element $u_0^\star$ of the input sequence $u^\star$ is applied to the system \cite{MPCbook}. At the next physical sampling time, a new state $s$ is received, and problem \ref{eq:MPC0} is solved again, producing a new $u^\star$ and a new $u_0^\star$. MPC (\ref{eq:MPC0}) hence yields a policy
\begin{align}
\label{eq:MPC0:Policy} 
\pi_\mathrm{MPC}\left(s\right) = u_0^\star\,,
\end{align} 
with $u_0^\star$ solution of \ref{eq:MPC0} for $s$ given. For $\gamma\approx 1$, policy (\ref{eq:MPC0:Policy}) can provide a good approximation of the optimal policy $\pi^\star$ for an adequate choice of prediction horizon $T$, terminal cost $V_N$ and if the MPC model $f$ approximates the true dynamics sufficiently well.

\subsection{Model Predictive Control and RL as Dynamic Programming}

Ultimately, one can consider both model predictive control and reinforcement learning to be problems that simplify the general Dynamic Programming, or Stochastic Control, problem. This resemblance and the ultimate possibility of unifying analysis of algorithmic structure in these historically distinct research communities was pointed out in~\cite{powell2014clearing,powell2019unified}.

\begin{equation}\label{eq:DP}
\begin{aligned}
\inf\limits_{u(t) \in \mathcal{C}^u \subseteq \mathbb{R}^{n_u}, \forall t} \ & \mathbb{E}_{\xi} \bigg[ \sum_{t=1}^T F(x(t, \xi(1:t)), u(t, \xi(1:t)), \xi(t)) \bigg] \\
\text{subject to } & x(1, \xi(1)) = f(u(1, \xi(1)), x_0, \xi(1)), \\
& x(t, \xi(1:t)) = f(u(t, \xi(1:t)), x(t-1, \xi(1:t-1)), \xi(t)), \quad t = 2, \dots, T, \\
& h(u(t, \xi(1:t)), x(t, \xi(1:t))) \leq 0, \quad t = 1, \dots, T, 
\end{aligned}
\end{equation}

Stochastic MPC can be considered a look-forward policy approximated with SAA. RL can be considered developing an empirical model wherein the constraints and the dynamics are completely encoded as embedded in the deep representation of the Value/Action-Value function.  Thus, it can be seen that both schemes are different means of approximating a dynamic programming problem. In MPC, a continuous-time control-state system is discretized, and an iterated sequence of lookahead problems is solved. In RL, trajectory sampling as information gathering is a central component of the algorithmic process, however only problems with finite state and action spaces become tractable. This choice of approximation/modeling limitation comes from the domain, and for our purposes, in order to ensure physical constraints on the system, we will use discrete time but continuous state-action space representations . However, rather than being some latent variable or structure in the value function, the influence of  the actions and states are now subject to nonlinear equation systems. Thus domain physics knowledge constrains the continuous space sufficiently so as to become computationally tractable, at least with approximate techniques.

\subsection{Model Predictive Control Incorporation within RL}
In recent years the integration of Reinforcement Learning (RL) with Model Predictive Control (MPC), known as RL-based MPC (RLMPC), has gained significant attention for addressing the limitations of conventional MPC methods in dynamic and uncertain environments. This approach leverages MPC as a function approximator within RL \cite{reiter2025synthesis} and has been proven to yield optimal solutions to Markov Decision Processes (MDPs) even when the predictive model is inaccurate \cite{gros2019data}. The RLMPC framework proposed in \cite{gros2019data} combines the strengths of both paradigms: MPC ensures explainable policies with guarantees, while RL improves closed-loop performance of MPC by adapting it directly from data to compensate for model inaccuracies. Furthermore, RLMPC can encode explicit safety constraints to provide a safe-RL framework \cite{reiter2025synthesis}. These advancements position RLMPC as a powerful framework for handling complex, constrained, and adaptive decision-making problems across various applications.

RLMPC approaches have been increasingly applied to real-world systems, especially in energy systems and robotics. In energy systems, RLMPC has been applied to optimize peak power management and energy flexibility in smart grids, including multi-agent residential systems with local renewable energy production and storage \cite{reiter2025synthesis}. In robotics, RLMPC has been used for autonomous surface vessels (ASVs) to enhance their performance in tasks such as collision avoidance, docking, and trajectory tracking \cite{reiter2025synthesis}. 

The performance of an MPC scheme relies heavily on the accuracy of its predictive model, yet many real-world systems are inherently difficult to model precisely \cite{ anand2024data}. This model dependency is a key limitation of MPC, as selecting the model that maximized the closed-loop performance is non-trivial, and a more accurate model does not always guarantee better performance. However, it has been shown that MPC can still achieve optimality under an inaccurate model by appropriately adjusting its parameters \cite{gros2019data}. The integration of Reinforcement Learning (RL) with MPC aims to overcome this limitation by using RL to adaptively improve the closed-loop performance. Recent studies have established that when focusing on closed-loop performance, an MPC scheme can be adjusted holistically, allowing for the design of optimal policies without relying on a completely accurate model of the system. As another strategy, we can also see the use of explicit information learning methods from Partially Observed MDPs applied to learning an uncertain model in~\cite{matias2022simultaneous}. 

Now we will introduce the necessary background on RLMPC. Consider an MPC-based policy, 
\begin{align}
\pi_{{\theta}}( s) =  u^\star_0\,, \label{eq:MPC:Pi:fromU} 
\end{align}
where $ u^\star_0$ is part of the solution of:
\begin{subequations}
\label{eq:MPC:Pi}
\begin{align} 
 x^\star, u^\star=\mathrm{arg} \min_{ x, u}&\quad V_{{\theta}} \left( x_N\right) + \sum_{k=0}^{N-1}\, F_{{\theta}}\left( x_k, u_k\right), \label{eq:MPC:Pi:Cost}\\
\mathrm{s.t.} &\quad  x_{k+1} =  f_{{\theta}}\left( x_k, u_k\right),\quad   x_0 =  s\,,  \label{eq:MPC:Pi:Dyn} \\
&\quad  h_{{\theta}}\left( x_k, u_k\right) \leq 0,\quad  u_k \in \mathcal{A}\,. \label{eq:MPC:Pi:Const}
\end{align}
\end{subequations}
This MPC formulation is identical to \eqref{eq:MPC0}, but the cost, constraints, and dynamics underlying the MPC scheme are now all parameterized in ${\theta}$, to the exception of the input constraint $ u_k \in \mathcal{A}$. Besides the MPC-based policy \eqref{eq:MPC:Pi:fromU}, MPC can be construed as a (possibly local) model of the action-value function $Q^\star$. We can define a parameterized action-value function $Q_{{\theta}}$ based on the MPC scheme \eqref{eq:MPC:Pi}, as a model of the optimal action-value function $Q^\star$ as follows:
 
\begin{subequations}
\label{eq:MPC:Q}
\begin{align}
Q_{{\theta}}( s, a) = \min_{ x, u}&\quad \eqref{eq:MPC:Pi:Cost},\\
\mathrm{s.t.} &\quad \eqref{eq:MPC:Pi:Dyn} - \eqref{eq:MPC:Pi:Const},\quad  u_0 =  a \label{eq:MPC:Q:II}\,,
\end{align}
\end{subequations}

where a constraint $ u_0 =  a$ included in \eqref{eq:MPC:Q:II} is the only difference to \eqref{eq:MPC:Pi}. MPC \eqref{eq:MPC:Q} is a valid model of $Q^\star$ in the sense that it satisfies the following Bellman relationships on the true system
\begin{subequations}
\label{eq:Bellman0} 
\begin{align}
V^\star\left( s\right) &=\min_{ a}\, Q^\star\left( s, a\right)\,, \label{eq:Bellman0:V} \\
Q^\star\left( s, a\right) &= F\left( s, a\right) + \gamma \mathbb E\left[V^\star\left( s_+\right)\,|\, s, a\, \right]\,. 
\end{align}
\end{subequations}
and the optimal policy
\begin{align}
\label{eq:OptPolicy0} 
 \pi^\star\left( s\right) =\mathrm{arg} \min_{ a}\, Q^\star\left( s, a\right)\,. 
\end{align}
i.e.:
\begin{subequations}
\begin{align}
V_{{\theta}}( s) &=\min_{ a}\,Q_{{\theta}}( s, a)\,,\label{eq:MPC:V}\\ 
\pi_{{\theta}}\left( s\right) &= \mathrm{arg}\min_{ a} Q_{{\theta}}( s, a), \label{eq:MPC:PiViaMin}
\end{align}
\end{subequations}
where $ V_{{\theta}}( s)$ is the optimal cost resulting from solving MPC \eqref{eq:MPC:Pi}. If the MPC parameters $\theta$ are such that $Q_{{\theta}} = Q^\star$, then the MPC scheme \eqref{eq:MPC:Pi} delivers the optimal policy $ \pi^\star$ through \eqref{eq:MPC:Pi:fromU}, i.e. $\pi_{{\theta}}= \pi^\star$. An important question, then, is how effectively an MPC scheme can approximate $Q^\star$ at least in a neighborhood of $ a =  \pi^\star\left( s\right)$.  The central result in RL-based MPC provided in Theorem 1 in \cite{gros2019data} answers this question and establishes it under some mild conditions, \eqref{eq:MPC:Q} provides an exact model of $Q^\star$ even if the predictive model \eqref{eq:MPC:Pi:Dyn} is inaccurate. The Theorem 1 in \cite{gros2019data}  states the following: 

\textit{Considering the parameterized stage cost, terminal cost, and constraints in \eqref{eq:MPC:Pi} are universal function approximations (i.e., can approximate a given function accurately) with adjustable parameters $\theta$. Then there exist parameters $\theta^\star$ s.t. the following identities hold $\forall \gamma$:
\begin{enumerate}
  \item $V_{\theta^\star}( s)=V^\star( s),\,\forall  s\in\Omega$\label{eq:VV}
  \item ${\pi}_{\theta^\star}( s)=\pi^\star( s),\,\forall  s\in\Omega$\label{eq:piV}
  \item $ Q_{\theta^\star}( s, a)=Q^\star( s, a),\,\forall  s\in\Omega$, for the inputs $ a\in \mathcal{A}$ such that $\lvert V^\star( f_{\theta^\star}( s, a))\rvert<\infty$\label{eq:QQ}
\end{enumerate}
if the set 
\begin{align}\label{eq:assum:Theo}
\Omega&:=\left\{ s\in \mathcal{S}\,\,\Big|\,\,\left| V^\star( x^\star_t)\right|<\infty, \ \forall\, t \leq  T\right\}\,,
\end{align}
is non-empty, where $ x^\star_t$ is the optimal state solution of \eqref{eq:MPC:Pi}.}

This in turn entails that MPC \eqref{eq:MPC:Pi} can provide the optimal solution to the true MDP even if the MPC model is inaccurate. This can be achieved by selecting the appropriate stage cost, terminal cost, and constraints in the MPC scheme. The applicability of this results is very broad and extends, e.g., to robust MPC, stochastic MPC, mixed-integer MPC, and Economic MPC (EMPC), and is valid both for the discounted and undiscounted settings \cite{reiter2025synthesis}. The role of RL in that context is to provide these learning tools .  RL-based MPC  integrate RL tools that aim at minimizing $J( \pi_{\theta})$, or at achieving $Q_{{\theta}} \approx Q^\star$ from the data. There are alternatives to RL for adjusting the MPC parameters, such as e.g.
Bayesian optimization. However, RL is often regarded as the most effective technique, especially if the number of parameters to adjust is not very small. 

Both Q-learning and policy gradient approaches can be used for RL-based MPC \cite{reiter2025synthesis}. In case of Q learning aim to satisfy (\ref{eq:QQ}), where  $Q_{\theta}$ is delivered by \eqref{eq:MPC:Q} and $ \min_{ a^\prime}Q_{\theta}( s, a^\prime)= V_{{\theta}}( s) $ is delivered by (\ref{eq:MPC:PiViaMin}), i.e. by the optimal cost of MPC \eqref{eq:MPC:Pi} solved at state $ s$. Alternatively policy gradient approach aim to adjust the parameters $\theta$ of the MPC as a policy \eqref{eq:MPC:Pi:fromU} to minimize cumulative cost $J( \pi_{\theta})$ directly to satisfy \eqref{eq:piV}. This is the standard approach for systems with continuous action spaces for which we can distinguish between two methods. Policy gradient methods estimate the policy gradient $\nabla_{\theta}J( \pi_{\theta})$ from data to update the parameters $\theta$ in a gradient descent fashion. Alternatively, \textit{direct policy search methods} build a surrogate model of $J( \pi_{\theta})$ from data and use it to propose new policy parameters $\theta$.  However, direct policy search tends to scale poorly with the size of the policy parameters. 

One common approach in policy gradient methods is the \textit{deterministic policy gradient} by~\cite{silver2014deterministic} with the following expression for the gradient of the cumulative cost:
\begin{align}
\nabla_{\theta}J( \pi_{\theta}) = \mathbb E\left[\nabla_{\theta} \pi_{\theta}\left( s\right) \nabla_{ a} Q_{ \pi_{\theta}}\left( s,\pi_{\theta}( s)\right)  \right]\,.
\end{align}
Here $\pi_{\theta}$ is delivered by MPC~\eqref{eq:MPC:Pi}. The \textit{critic} $Q_{ \pi_{\theta}}$ is typically built separately using a generic function approximator and policy evaluation techniques. It has been shown that the DPG method can be simplified in the case of RLMPC by approximating the required value function through a first-order expansion of the value function estimate obtained from the MPC policy itself \cite{anand2023painless}.

\subsubsection{Sensitivities}
Many RL methods, including Q-learning and policy gradient methods, require the sensitivities of the function approximators $\pi_{\theta}$, $Q_{\theta}$, and $V_{\theta}$. When approximated by an MPC scheme, these are typically continuous but only piecewise smooth. 
However, when the MPC achieves Linear Independence Constraint Qualification (LICQ) and Second Order Sufficient Condition (SOSC), non-smooth points correspond to weakly active constraints in the MPC. These points form a set of zero measures for a well-formulated MPC scheme. Because RL methods always use the sensitivities inside expected value operators, their contribution to the learning disappears as long as the state transition density of the true MDP is bounded. The non-smoothness of MPC schemes may superficially appear as an issue for RL, but fortunately, in all practical cases, it is not. This issue can be ignored when the MPC response is continuous.

The arguments above do not necessarily hold anymore if the MPC does not satisfy SOSC and, as a result, produces discontinuous policies, e.g., a \textit{bang-bang} response. For example, designing the MPC as a Linear Program might lead to such a situation, as well as processes with switching operations or more decisions than state equations. We can then still use Q-learning since $Q_{\theta}$ and $V_{\theta}$ typically remain continuous and piecewise smooth. However, the discontinuities in the policy lead to problems in policy gradient methods. In addition, extending optimization methods that perform MPC for problems with degenerate constraints, those that do not satisfy the LICQ, like~\cite{kungurtsev2017predictor}, to incorporate RL, remains an important open problem. 

\subsubsection{Real-time MPC}
An MPC problem at a give initial state $s$ is a Nonlinear Program of the form
\begin{subequations}
\label{eq:MPCNLP}
\begin{align}
\min_{w}&\quad \Phi(w) \\
\mathrm{s.t.}&\quad g(w,s)=0\label{eq:MPCNLP:eq} \\
&\quad h(w) \leq 0
\end{align}
\end{subequations}
where $w$ gathers the discrete-time state and input trajectories. The equality constraints \eqref{eq:MPCNLP:eq} typically ensure that the state trajectories abide by the physical model \eqref{eq:dyn0} used in formulating the MPC problem. 

Defining the Lagrange function associated to \eqref{eq:MPCNLP}:
\begin{align}
\mathcal L(w,\lambda,\mu) = \Phi(w) + \lambda^\top g(w,s) + \mu^\top h(w)
\end{align} 
for the dual variables $\lambda,\mu$, a solution to NLP \eqref{eq:MPCNLP} satisfies the KKT conditions:
\begin{subequations}
\label{eq:NLPKKT}
\begin{align}
\nabla_w\mathcal L(w,\lambda,\mu) = 0\\
g(w,s)=0 \\
h(w) \leq 0,\quad \mu \geq 0 \\
\mu^\top h(w) = 0 
\end{align}
\end{subequations}
for dual variables $\lambda,\mu$. If NLP \eqref{eq:MPCNLP} holds a sufficiently differentiable cost $\Phi$ and constraints $f,h$, it can be solved via Successive Quadratic Programming (SQP), whereby problem \eqref{eq:MPCNLP} is successively approximated by the Quadratic Program (QP) 
\begin{subequations}
\label{eq:MPCQP}
\begin{align}
\min_{\Delta w}&\quad \frac{1}{2}\Delta w^\top H \Delta w + \nabla_w\Phi(w)^\top \Delta w \\
\mathrm{s.t.}&\quad g(w,s) + \nabla_w g(w,s)^\top \Delta w=0\label{eq:MPCQP:eq} \\
&\quad h(w) + \nabla_w h(w)^\top \Delta w\leq 0
\end{align}
\end{subequations}
where matrix $H$ is the Hessian with respect to $w$ of the Lagrange function $\mathcal L$ at $w,\lambda,\mu$. The solution of QP \eqref{eq:MPCQP} provides the primal-dual update:
\begin{align}\label{eq:update}
w \leftarrow w + \alpha \Delta w,\quad \lambda \leftarrow  \lambda + \alpha \Delta \lambda,\quad \lambda \leftarrow  \mu + \alpha \Delta \mu
\end{align}
where $\Delta\lambda$ and $\Delta\mu$ are given by the change in the dual variables associated to \eqref{eq:MPCQP} from the previous QP solution, and $\alpha\in (0,1]$ is a step size selected via a line-search strategy. SQP iterates solving QP \eqref{eq:MPCQP} and updating the primal-dual variables \eqref{eq:update} until convergence to a solution $w^\star(s), \lambda^\star(s), \mu^\star(s)$ to \eqref{eq:NLPKKT} for a given $s$. The SQP procedure start with an initial guess $w_\mathrm{init}$ at which the first QP \eqref{eq:MPCQP} is formed. If $w^\mathrm{init}$ is sufficiently close to $w^\star$, then full steps $\alpha=1$ can be used in \eqref{eq:update}, and the SQP iteration converges quadratically to $w^\star(s)$.

A key observation for real-time MPC is that the initial condition $s_k$ at a given discrete time step $k$ are typically close to their prediction formed at time $k-1$. More specifically, considering the predicted optimal state $x_1^\star(s_{k-1})\subset w^\star(s_{k-1})$ at time $k-1$ is typically a good approximation of $s_{k}$. This observation leads to the \textit{shifting strategy} to provide an initial guess for NLP \eqref{eq:MPCNLP} at $s_k$ from the solution obtained at time $k-1$ for $s_{k-1}$. The shifting strategy consists in using the optimal state-input trajectories $x^\star_{0,\ldots,N}(s_{k-1})$ $u^\star_{0,\ldots,N-1}(s_{k-1})$ contained in $w^\star(s_k)$ to form the state-input initial guesses:
\begin{align}
\label{eq:Shift}
    x^\mathrm{init}_i(s_k) &=  x^\star_{i+1}(s_{k-1}),\quad u^\mathrm{init}_i(s_k) &=  u^\star_{i+1}(s_{k-1}),\quad i=0,\ldots,N-1
\end{align}
The new state and input $x^\mathrm{init}_N(s_k)$, $u^\mathrm{init}_{N-1}(s_k)$ cannot be built from \eqref{eq:Shift}, and are often simply guessed using 
\begin{subequations}
\label{eq:ShiftEnd}
\begin{align}
u^\mathrm{init}_{N-1}(s_k) &= u^\star_{N-2}(s_{k-1}),\\ 
x^\mathrm{init}_{N}(s_k) &= f(x^\star_{N-1}(s_{k-1}),u^\star_{N-2}(s_{k-1}))
\end{align}
\end{subequations}
More advance strategies are possible but do not always make a significant difference in the performance of the SQP iteration. Initial guesses for the dual variables $\lambda,\mu$ can also be built using a similar process.

A key to ensuring that $x_1^\star(s_{k-1})\approx s_{k}$ is to make sure that the sampling time of the MPC scheme is fast enough to ensure that the prediction error cannot accumulate during the physical time interval between the discrete times $k-1$ and $k$. When this strategy is sufficiently effective, the SQP procedure using the initial guess \eqref{eq:Shift} often fully converges in a very few iterations. An extreme version of SQP for real-time MPC is then to perform a single SQP iteration at each time sample $k$ of MPC, and rely on the strong quadratic contraction resulting from the update \eqref{eq:update} from the initial guess \eqref{eq:Shift}-\eqref{eq:ShiftEnd} to keep the (not fully converged) primal-dual solutions close to their optimal values $w^\star(s), \lambda^\star(s), \mu^\star(s)$. This principle is at the core of the \textit{Real-Time Iteration} (RTI). In addition to circumventing multiple SQP iterations, RTI seeks to minimize the time spent between obtaining a state measurement $s_k$ and delivering the corresponding approximation to the optimal solution $w^\star(s_k)$. Indeed, delivering $w^\star(s_k)$ requires performing two computationally expensive operations: forming QP \eqref{eq:MPCQP}, i.e. computing $H,\nabla_w\Phi(w),g(w,s_k), \nabla_w g(w,s),h(w),\nabla_w h(w)$, and solving it. A key observation behind the RTI is that it is possible to form QP \eqref{eq:MPCQP} \textit{before} receiving $s_k$. Hence the RTI performs the update of the MPC solution in two phases. The \textit{preparation phase} consists in performing the shifting \eqref{eq:Shift}-\eqref{eq:ShiftEnd} and preparing QP \eqref{eq:MPCQP}, while the \textit{feedback phase} consists in solving QP \eqref{eq:MPCQP} once $s_k$ is received. 

\subsubsection{Exploration in RLMPC}
Exploration is a core requirement for all RL methods, i.e., regularly applying actions $\mathbf a\neq \mathbf\pi_{\mathbf\theta}(\mathbf s)$ to the real system to gather the information necessary to improve the policy. Given the constraints in \eqref{eq:MPC:Pi:Const}, a natural question is how to generate exploration that does not jeopardize the MPC feasibility.  A straightforward way to address this is by manipulating the MPC cost such that the resulting action differs from the original policy provided by \eqref{eq:MPC:Pi}. One possible modification is to add a gradient over the initial action $\mathbf u_0$, which yields the MPC scheme
\begin{subequations}\label{eq:MPC:Explo}
\}
\begin{align}
\phi_{\mathbf\theta}\left(\mathbf s,\mathbf d\right) = \min_{\mathbf x,\mathbf u}&\, \mathbf d^\top\mathbf u_0 + T_{\mathbf{\theta}} \left(\mathbf x_N\right) + \sum_{k=0}^{N-1}\, L_{\mathbf{\theta}}\left(\mathbf x_k,\mathbf u_k\right) \label{eq:RMPC:Pi:Cost}\\
\mathrm{s.t.} &\quad  \eqref{eq:MPC:Pi:Dyn},\,\eqref{eq:MPC:Pi:Const}\,,
\end{align}
\end{subequations}
where $\mathbf d\in \mathbb R^m$ is a vector of the size of the action, possibly selected randomly. Because only the MPC cost is modified, the solution of \eqref{eq:MPC:Explo} is feasible for  \eqref{eq:MPC:Pi}. One can then use MPC~\eqref{eq:MPC:Explo} to produce feasible exploration for policy gradient methods and to generate action-value functions for Q-learning. In that context, \eqref{eq:MPC:Explo} produces a feasible action with exploration $\mathbf a = \mathbf u_0^\star$ where $\mathbf u_0^\star$ is the solution of \eqref{eq:MPC:Explo} and depends on $\mathbf d$. The optimal cost $\phi_{\mathbf\theta}$ of \eqref{eq:MPC:Explo} delivers the action-value function
\begin{align}
Q_{\mathbf\theta}\left(\mathbf s,\mathbf a\right)  = \phi_{\mathbf\theta}\left(\mathbf s,\mathbf d\right) - \mathbf d^\top\mathbf u_0^\star\,.
\end{align}
 see \cite{gros2019data} for more details. However, the feasibility of \eqref{eq:MPC:Explo} does not guarantee that applying the resulting policy to the real system will not violate the constraints \eqref{eq:MPC:Pi:Const}. For example, consider the effect of stochastic dynamics or model errors. This issue of safe exploration is discussed in the next section.

An additional technique is presented in~\cite{matias2022simultaneous}. Here, one solves the modified problem:

\begin{subequations}
\}
\begin{align} 
\phi_{\mathbf\theta}\left(\mathbf s,\mathbf d\right) = \min_{\mathbf x,\mathbf u}&\,  T_{\mathbf{\theta}} \left(\mathbf x_N\right) + \sum_{k=0}^{N-1}\, \left\{L_{\mathbf{\theta}}\left(\mathbf x_k,\mathbf u_k\right) +JS(\mathbf x_k;Q)\right\}\label{eq:RMPC:Pi:Cost}\\
\mathrm{s.t.} &\quad  \eqref{eq:MPC:Pi:Dyn},\,\eqref{eq:MPC:Pi:Const}\,,
\end{align}
\end{subequations}
where $JS(\mathbf x_k;Q)$ is the Jensen-Shannon divergence, an information metric, associated with the addition information provided to $Q$ upon potential visit to state $\mathbf x_k$. This includes the information content as by the reduction in uncertainty of the RL functions explicitly.

\subsubsection{Data Drive Predictive Control}
In scenarios where explicit state-space models are unavailable, recent advancements have extended the framework of RLMPC to operate without predefined notions of state or explicit dynamics models \cite{reiter2025synthesis}. This broadens the applicability of RLMPC to systems with partially or completely unknown dynamics. Similarly, other studies have introduced end-to-end learning approaches for Nonlinear MPC (NMPC) policies, enabling the derivation of control strategies directly from data without relying on first-principles models \cite{reiter2025synthesis}. These methods are particularly relevant in settings where deriving accurate system models is impractical or infeasible.

Computational complexity is another major concern when combining RL with MPC, especially when gradient-based methods are used. These methods often require frequent evaluations of computationally expensive MPC schemes, limiting their scalability to large datasets. To address this, approaches have been proposed to learn parameterized MPC schemes directly from offline data, bypassing the need to solve MPC problems repeatedly over the collected dataset \cite{reiter2025synthesis}. Broadly speaking, however, this is an important open area of research.

\subsubsection{Safe Reinforcement Learning with RLMPC}
The research community is increasingly considering the application of RL to safety-critical systems. A straightforward way to define safety is through a set of critical constraints which should not be violated, typically expressed as functions of the system state. A classic approach to safe RL is to learn the policy in silico and employ a \textit{pessimistic} model of the real system, meaning that it overestimates the probability of violating critical constraints. During learning, one can then use high penalties for violations of the critical constraints, e.g. through the use of barrier functions . RL will then naturally adjust the policy to avoid these penalties. 

In MPC, using a pessimistic model of the system is common within the Robust MPC (RMPC) methods. Starting from a pessimistic model, RMPC builds a safe policy by ensuring that even the worst-case predictions satisfy the critical constraints at all future times. To that end, the deterministic model \eqref{eq:MPC:Pi:Dyn} is replaced by a model that describes the evolution of sets enclosing all possible future trajectories, leading to an expression of the form
\begin{align}
\label{eq:RMPC:DispersionModel}
\mathbb X_{k+1} = \mathbf f_{\mathbf\theta}(\mathbb X_{k}, \mathbf\pi_\mathrm{c}(\mathbb X_{k},\mathbf x_k,\mathbf u_k) ) \oplus \mathbb W_{\mathbf\theta}\,.
\end{align}
One then designs the policy $\mathbf\pi_\mathrm{c}$ to \textit{manage} the growth of the sets $\mathbb X_{k}$, usually by operating on their deviation from a reference trajectory $\mathbf x_k$. The set $\mathbb W_{\mathbf\theta}$ represents process noise to capture the prediction uncertainties. It is possible to learn \eqref{eq:RMPC:DispersionModel} from data through \textit{set-membership identification}. We can then build the following MPC scheme to enforce satisfaction of constraints~\eqref{eq:MPC:Pi:Const} explicitly:
\begin{subequations}
\label{eq:RMPC:Pi}
\begin{align} 
\min_{\mathbf x,\mathbf u,\mathbb X}&\quad T_{\mathbf{\theta}} \left(\mathbf x_N\right) + \sum_{k=0}^{N-1}\, L_{\mathbf{\theta}}\left(\mathbf x_k,\mathbf u_k\right) \label{eq:RMPC:Pi:Cost:1}\\
\mathrm{s.t.} &\quad \mathbf x_{k+1} = \mathbf f_{\mathbf{\theta}}\left(\mathbf x_k,\mathbf u_k\right),\quad \eqref{eq:RMPC:DispersionModel}  \label{eq:RMPC:Pi:Dyn} \\
&\quad \mathbf h_{\mathbf{\theta}}\left(\mathbb X_{k},\mathbf\pi_\mathrm{c}(\mathbb X_{k},\mathbf x_k,\mathbf u_k)\right) \leq 0,\label{eq:RMPC:Pi:Const}  \\
&\quad \mathbf\pi_\mathrm{c}(\mathbb X_{k},\mathbf x_k,\mathbf u_k) \in \mathcal A \\
&\quad \mathbb X_{N} \in \mathbb T_{\mathbf\theta},\quad \mathbf x_0 = \mathbf s,\quad \mathbb X_0 = \mathbf s \,.\label{eq:RMPC:Pi:Const:TerminalInit}
\end{align}
\end{subequations}
MPC \eqref{eq:RMPC:Pi} generates a safe policy $\mathbf\pi_{\mathbf{\theta}}= \mathbf u^\star_0$ by construction for $\mathbb T_{\mathbf\theta}$ adequately chosen and if \eqref{eq:RMPC:DispersionModel} accounts for the worst case situations observed in the data, see \cite{MPCbook}. Application of RL to adjust RMPC schemes has been developed recently \cite{reiter2025synthesis}. In the context of RL for RMPC-based policies, it is useful to stress that there is a \textit{hard} separation between learning for safety and learning for closed-loop performance. Indeed, safety is learned via set-membership identification on \eqref{eq:RMPC:DispersionModel} and then enforced in the RMPC scheme by construction. Closed-loop performance is optimized using RL in parallel.

For the sake of clarity, we ought to underline that the adjustment of the constraints \eqref{eq:RMPC:Pi:Const} in the RMPC scheme and of set $\mathbb W_{\mathbf\theta}$ ought to be done with care in order to preserve safety. In particular, the adjustment of set $\mathbb W_{\mathbf\theta}$ must ensure that model \eqref{eq:RMPC:DispersionModel} accounts for all past data points in the set-membership sense. Arguably, safety-critical constraints in \eqref{eq:RMPC:Pi:Const} ought not to be modified.

\textbf{Safe Exploration.} 
Safe exploration is difficult to produce without a model of the system in the form \eqref{eq:RMPC:DispersionModel}, which can predict the worst-case evolution of the system, and assess the impact of the exploration on the system safety. However, even with a model \eqref{eq:RMPC:DispersionModel} of the system, it can be expensive to verify the safety of an input differing from the safe policy, and even more expensive to build the set of safe inputs.

It is straightforward to apply the feasible exploration approach discussed earlier to the RMPC formulation \eqref{eq:RMPC:Pi}. Given that we only modify the RMPC cost, the resulting solution is feasible for \eqref{eq:RMPC:Pi}, therefore enabling safe exploration if \eqref{eq:RMPC:Pi} yields a safe policy.

\textbf{Safe Policies and Safe Learning.}
There are two ways to implement safe learning online. One is in a batch fashion, meaning that parameter updates require collecting a minimum amount of transition data. The alternative is to perform parameter updates at every time step. In safety-critical applications, safety must always be preserved in either case.

However, while taking actions from a unique, safe, and stable policy ensures the stability and safety of the system, taking actions from a sequence of safe and stable policies may not. That is because a sequence of policies does not necessarily inherit the properties of the individual policies. We discuss these conditions in detail in \cite{gros2022learning}.


\subsection{Soft Constraints}
In many real-world safe control applications, the ``estimated'' environment models are available (or can be built) to derive sub-optimal control policy priors \cite{liepe2014framework, hippisley2017development}. Control policy priors of reasonable accuracy are essential to many applications for critical systems where control failure is consequential, especially systems with a single-life setting where failure even during training is unacceptable. Consider the example of devising a control policy that prescribes the optimal drug dosages for regulating a patient's health status. This is a critical system with a single-life setting where no harm to the patient is tolerated during policy exploration. From available records of other patients, an estimated patient model can be built to predict the response to different drug dosages and ensure adherence to the safety bounds (set based on clinical knowledge). However, a new patient's response can deviate from the estimated model, which poses a significant challenge in control adaptability and patient treatment performance.

In applications analogous to the above example, i.e., control in critical systems, MPC derived from the estimated models can be used as a decision oracle that informs the RL learning process and regulates RL policy updates \cite{cheng2019control}. MPC's regulatory capabilities are leveraged in the recently proposed RL-AR (Reinforcement Learning with Adaptive Regularization) \cite{tianreinforcement},  an algorithm for safe training and deployment of RL that features convergence to the optimal RL policy. RL-AR comprises two parallel agents and an adaptive focus module: (i) The \textit{safety regularizer} agent follows a deterministic policy $\pi_\reg: \mathcal{S}\rightarrow\mathcal{A}$ proposed by an MPC; (ii) The \textit{off-policy RL} agent is an adaptive agent with $\pi_\rl: \mathcal{S}\rightarrow\mathcal{M}(\mathcal{A})$ that can learn from an acting policy that is different from $\pi_\rl$; (iii) The \textit{focus module} learns a state-dependent weight $\beta: \mathcal{S} \rightarrow [0,1]$ for combining the deterministic $a_\reg(s)=\pi_\reg(s)$ and the stochastic $a_\rl(s)\sim \pi_\rl(s)$ actions proposed by the two agents. The safety regularizer component has a static built-in estimated environment model $\tilde{f}:\mathcal{S}\times\mathcal{A} \rightarrow \mathcal{S}$ (that is different from the actual environment), while the off-policy RL agent and focus module are dynamically updated using observed interactions in the actual environment.

The RL-AR workflow is as follows: (i) $\pi_\reg(s)$ generates $a_\reg(s)$, which hard-codes safety constraints in the optimization problem over a period forecasted by $\tilde{f}$. The forecasting ability anticipates and prevents running into unsafe states for the critical system; (ii) $\pi_\rl(s)$ generates $a_\rl(s)$ to allow stochastic exploration and adaptation to the actual environment; (iii) $\beta(s)$ applies the following policy combination:
\begin{align}
\begin{split}
     a_\beta(s) = \beta(s) a_\reg(s) + (1-\beta(s)) a_\rl(s).
\end{split}
\label{eqmixp}
\end{align}
The state-dependent focus weight $\beta(s)$ is initialized to $\beta(s)\geq 1-\epsilon$, $\forall s\in\mathcal{S}$, hence prioritizing the safe $\pi_\reg$ before $\pi_\rl$ learns a viable policy. As more interactions are observed for a state $s$ and the expected return of $\pi_\rl(s)$ improves, $\beta(s)$ gradually shifts the focus from the initially suboptimal $\pi_\reg(s)$ to $\pi_\rl(s)$. To update $\beta(s)$, the following objective is used:
\begin{align}
\begin{split}
        \beta'(s) = \argmax_{\beta \in [0, 1]} \mathbb{E}\left[Q^{\pi_\beta}(s, \beta a_\reg(s)+(1-\beta) a_\rl(s))\right].
\end{split}
\label{eqattp}
\end{align}

In \cite{tianreinforcement}, it is analytically shown that if the reward function $R$ and transition probability $P$ of the MDP are Lipschitz continuous, then when a state $s$ has not been sufficiently explored and its corresponding $\beta(s)$ updates are limited, the sub-optimality of the RL policy $\pi_\rl$ has a limited impact on the expected return of the combined policy $\pi_\beta$, which is the actual acting policy. Meanwhile, as learning progresses and $\beta(s)$ approaches 0, RL-AR allows unbiased convergence to the optimal policy. RL-AR was shown effective in many safe RL applications ranging from medical to chemical systems, with a significant tolerance of estimated model deviation from the actual environment (see Section 4 in \cite{tianreinforcement}).

\subsection{Learning State Dynamics}

In the current integration of MPC with RL, the predominant framework in the control community is using discretized differential equations to describe the state evolution and thus hard physics defining boundaries of safety and other such considerations. Safety is guaranteed with the presence and enforcement, in a dedicated optimization procedure, of deterministic functional constraints of a single uniquely defined discretized dynamics, with the learning performed to inform the reward to objective layer. 

However, consider the possibility that the physical dynamics have not been completely identified and quantified. Yet we stay a distance from the completely agnostic framework of NN representations, but instead consider that the equations have multiple subcomponent terms, some known exactly and some partially, and thus must be learned from data.

The broad framework can be referred to as \emph{Grey Box} modeling. The purpose therein is to use any and all a priori physics and first principles based knowledge to formulate the dynamics of the process of interest. Meanwhile, the remaining unknown system dynamics is empirically modeled with data. 

Indeed, all of the component functions in the MPC, including the dynamics $f$ and the constraints $c$, can be subject to learning, training and inference. The use of physics informed biases in the parametrized neural model for $f$ for instance could be used to enhance learning the precise differential system of equations. 

As a recent example of NN-informed Grey-box identification, consider~\cite{dong2022optimization,dong2024descent}. Broadly speaking, a powerful and effective approach is to consider a semilinear system,
\[
\mathcal{L}(u) +\mathcal{N}(u)=f,\,\nabla\cdot A\nabla u+\Phi(u;\theta) = f
\]
where $\mathcal{L}$ is a linear operator known to precision, and $\mathcal{N}$ is nonlinear and must be learned, in general, to produce a neural network model $\Phi(\cdot;\theta)$. For computing the solution of the MPC, the full estimated model is used to define the dynamics.

Robotics presents a natural application of hybrid RL-MPC online decision tools. Robotics has been a motivating application of RL as due to its enormous potential utility across a variety of domains. See for example the early survey~\cite{kormushev2013reinforcement}. Given their increasing role in automated manufacturer, home care, and other domains with vulnerable people or material, safety assurances become a priority for trusted stakeholder uptake. For a recent summary of work enforcing safety in robotics applications, see, e.g.~\cite{brunke2022safe}. 

As a field of control, robotics presents a number of challenges as far as integrating the mechanical operative features in a manner that is accurate in modeling while also efficient with data. Entire fields of Mechatronics and Cybernetics include deep investigations as to the efficient and reliable interplay of physical mechanics, electronics, control and computing required for robot operation. One notable feature that can assist the sample complexity and degrees of freedom searched by solution methods is the presence of digitized movement choices, that is rotation along a finite set of distinct axes in freespace. This  suggests the computational potential of scientific deep learning. Included is the use of geometric biases, e.g. the rotational movements of a robot often has a cyclic group structure which can be explicitly enforced in the modeling of the problem. The use of geometric learning and numerical methods can well be integrated in a control formulation, e.g. \cite{duruisseaux2023lie}.

\section{\rev{Task-Level Planning and Scheduling}}\label{sec:planning}

\subsection{Scheduling - Background}

Scheduling is a fundamental problem in AI that involves allocating resources, time, and tasks to achieve specific objectives under given constraints~\cite{dick2022scheduling}. Formally, a scheduling problem consists of a set of tasks \( T = \{t_1, t_2, \dots, t_n\} \), where each task \( t_i \) has associated properties such as its duration \( d(t_i) \), start time \( s(t_i) \), and finish time \( f(t_i) = s(t_i) + d(t_i) \). Additionally, there is a set of resources (processors) \( R = \{r_1, r_2, \dots, r_m\} \), where each resource \( r_j \) may have capacity constraints \( c(r_j) \), limiting the number of tasks that can use \( r_j \) simultaneously.

Constraints also play a critical role in scheduling problems. Precedence constraints ensure that certain tasks cannot start until others have finished, formalized as \( s(t_j) \geq f(t_i) \), where \( t_i \) must complete before \( t_j \) begins. Resource constraints ensure that, at any time \( t \), the number of tasks using a resource \( r_j \) does not exceed its capacity \( c(r_j) \). Temporal constraints, such as lower bounds on start times, impose additional limits, for instance, \( s(t_i) \geq \tau_i \) for a given \( \tau_i \).

The objective of a scheduling problem is to optimize a specific criterion (objective function). Common objectives include minimizing the makespan, defined as \( \min \max_{t_i \in T} f(t_i) \), which represents the total time required to complete all tasks~\cite{dick2022scheduling}. Other objectives include minimizing resource usage by reducing total cost or utilization and maximizing throughput, which aims to complete the maximum number of tasks within a given time limit. 

Scheduling problems can be generally classified based on their structural properties. In single-machine scheduling, tasks are executed on a single resource, often with the goal of minimizing makespan or ``lateness''. Parallel-machine scheduling extends this to multiple identical or unrelated resources~\cite{baker2018principles}. Job shop scheduling introduces a more complex structure, where tasks are grouped into jobs, each requiring a specific sequence of operations on different machines~\cite{baker2018principles}. Flow shop scheduling simplifies this by requiring that all jobs follow the same sequence of machines. More generally, constraint-based scheduling allows for arbitrary combinations of time, precedence, and resource constraints~\cite{baptiste2001constraint}.

Solving scheduling problems in AI is generally based on techniques from constraint satisfaction programming (CSP), linear programming, heuristic search, and metaheuristic optimization (e.g. genetic algorithms), which are employed to tackle the underlying combinatorial complexity~\cite{dick2022scheduling}.
Applications of scheduling in AI then span various domains. In manufacturing, scheduling generally optimizes machine and resource allocation to maximize production efficiency. In industrial robotics, multi-robot systems rely on scheduling to allocate tasks while respecting time and resource limits.
Staffing problems, such as in healthcare systems, are often framed as scheduling to allocate medical staff, equipment, and facilities, e.g. to maximize patient throughput and minimize wait times~\cite{ernst2004staff}.

\subsection{Classical Planning - Background}

Planning is another established framework for sequential decision making with a long tradition in AI~\cite{ghallab2004automated}.
Planning shares core similarities with MPC and RL in their common objective of finding optimal action sequences. 
Like MPC, planning operates under the assumption of a \textit{known} environment model, enabling proactive decision-making.
Similar to most RL, planning focuses on discrete state and action spaces, facilitating the application of \textit{search} algorithms.
However, planning distinguishes itself through its {goal-centric} approach and the use of \textit{relational} representations of states and actions, a key difference from the numerical representations typically used in MPC and the enumerated representations common in (tabular) RL.
This relational approach allows planning to handle more complex, abstract relationships and generalize over different objects and environments more effectively.

\subsubsection{Relational Logic}
\label{sec:logic}

First-order (predicate) logic is the cornerstone of all formal knowledge representations and reasoning in AI, providing the formal language for expressing facts, rules, and relationships between concepts~\cite{huth2004logic}.

\textit{\textbf{Syntactically}}, in first-order logic (FOL), we define a \textit{vocabulary} (signature) $\Sigma$ consisting of: 
a set of \textit{constant symbols} $\mathcal{C}$ (e.g., ``Room$_1$'',``Patient$_2$''); 
a set of \textit{variable symbols} $\mathcal{V}$ (e.g., position ``$x$'', distance``$d$''); 
a set of \textit{function symbols} $\mathcal{F}$ (e.g., ``$\text{neighborOf}(x)$''), each with an associated \textit{arity} $n \in \mathbb{N}$ denoted by $f/n$ where $f \in \mathcal{F}$; 
and a set of \textit{predicate symbols} $\mathcal{P}$ (e.g., ``$\text{Asleep}(x)$'' or ``$\text{Nursing}(x, y)$''), also with associated arities $P/n$ where $P \in \mathcal{P}$.
The set of \textit{terms} $\mathcal{T}(\Sigma)$ over $\Sigma$ is then defined recursively: $\mathcal{C} \subset \mathcal{T}(\Sigma)$, $\mathcal{V} \subset \mathcal{T}(\Sigma)$; and if $f/n \in \mathcal{F}$ and $t_1, \dots, t_n \in \mathcal{T}(\Sigma)$, then $f(t_1, \dots, t_n) \in \mathcal{T}(\Sigma)$ (e.g., ``$\text{nurseOf(Patient$_1$)}$''). 
The set of \textit{formulae} $\mathcal{F}(\Sigma)$ over $\Sigma$ is then also defined recursively as follows. 
If $P/n \in \mathcal{P}$ and $t_1, \dots, t_n \in \mathcal{T}(\Sigma)$, then $P(t_1, \dots, t_n) \in \mathcal{F}(\Sigma)$ is an atomic formula; if $\phi, \psi \in \mathcal{F}(\Sigma)$, then $\neg \phi, (\phi \land \psi), (\phi \lor \psi), (\phi \rightarrow \psi), (\phi \leftrightarrow \psi) \in \mathcal{F}(\Sigma)$, and if $\phi \in \mathcal{F}(\Sigma)$ and $x \in \mathcal{V}$, then $\forall x \phi$ and $\exists x \phi \in \mathcal{F}(\Sigma)$ (e.g., ``$\forall x (\text{Man}(x) \rightarrow \text{Mortal}(x))$'').

\textit{\textbf{Semantics}} in logic is then defined by the means of \textit{interpretations}~\cite{huth2004logic}. An {interpretation} $\mathcal{I} = (D, I)$ consists of a non-empty \textit{domain} $D$ and an \textit{interpretation function} $I$: 
for each $c \in \mathcal{C}$, $I(c) \in D$; for each $f/n \in \mathcal{F}$, $I(f): D^n \rightarrow D$; and for each $P/n \in \mathcal{P}$, $I(P): D^n \rightarrow \{\text{True}, \text{False}\}$. 
A \textit{variable assignment} (substitution) $v: \mathcal{V} \rightarrow D$ then maps variables to elements in $D$. 
The interpretation of a term $t$ under $\mathcal{I}$ and $v$ is denoted $t^{\mathcal{I},v}$. Truth of a formula $\phi$ is denoted $\mathcal{I}, v \models \phi$. $\mathcal{I}$ \textit{satisfies} $\phi$, denoted $\mathcal{I} \models \phi$, if $\mathcal{I}, v \models \phi$ for all $v$. 

The \textit{relational} logic is then an important function-free ($\mathcal{F} = \emptyset$) subset of FOL that underlies the representation formalism of classical planning.

\subsubsection{Planning}

The problem of planning is normally formalized in the language of \textit{predicate logic}~\cite{huth2004logic}. The world is represented by a set of {predicates} $\mathcal{P}$, where each $p \in \mathcal{P}$ has arity $\text{arity}(p)$. A state $s \in \mathcal{S}$ is then a truth assignment $s: \mathcal{P} \to \{\text{True}, \text{False}\}$. For a ground predicate $p(c_1, \dots, c_{\text{arity}(p)})$ over the respective constants $c_i$, $s(p(c_1, \dots, c_{\text{arity}(p)}))$ denotes its truth value.

Actions (controls) are operators $a \in \mathcal{A}$ with precondition $\text{pre}(a)$ and effect $\text{eff}(a) = (\text{eff}^{+}(a), \text{eff}^{-}(a))$, representing added and deleted predicates. Action $a$ is \textit{applicable} in $s$ if $s \models \text{pre}(a)$. Applying $a$ to $s$ results in a new state $s' = {\psi}(a, s)$:

\begin{equation}\label{eq:planning_state_transition}
s'(p) = \begin{cases}
\text{True} & \text{if } p \in \text{eff}^{+}(a) \\
\text{False} & \text{if } p \in \text{eff}^{-}(a) \\
s(p) & \text{otherwise}
\end{cases}
\end{equation}

A planning problem is defined by an initial state $s_0$ and a \textit{goal condition} $\phi_G$, represented generally by a logical formula.
A plan $\mathbf{a} = (a_0, \dots, a_{n-1}) \in \mathcal{A}^n$ achieves the goal if applying it to $s_0$ results in a state $s_n$ satisfying $\phi_G$:

\begin{equation}\label{eq:goalplan}
s_n \models \phi_G, \quad \text{where } s_{t+1} = \psi(a_t, s_t), \quad t = 0, \dots, n-1.
\end{equation}

Planning algorithms can be generally categorized as either satisficing or (less commonly) optimal, depending on whether they guarantee to achieve the goal or also finding the shortest plan.
Somewhat similar to MPC (Eq.~\ref{eq:MPC0}), optimal planning can be formulated as aiming to minimize a cost:

\begin{equation}
\mathbf{a}^* = \arg\min_{\mathbf{a}} C(\mathbf{a}), \quad C(\mathbf{a}) = \sum_{t=0}^{n-1} C(s_t, a_t) + C_T(s_n),
\end{equation}
subject to the state transitions defined by Eq. \eqref{eq:planning_state_transition} and the goal satisfaction~\eqref{eq:goalplan}. 
The state transition function $\psi$ (Eq. \eqref{eq:planning_state_transition}) can be thought of analogously to MPC (Eq. \ref{eq:MPC0}) and RL, but operates on relational structures instead of numerical vectors. Similarly, the reward function $R$ or cost $C$ is analogous to the cost function in optimal planning, but planning explicitly focuses on achieving the goal (terminal) state, which is most commonly the only reward.
Connecting these fields further, model-based RL~\cite{moerland2023model} often utilizes planning by learning a model and then using it for search, bridging the gap between data-driven learning and model-based control.


Conceptually, classical planning could be viewed as a special case of MPC where (i) the transition dynamics (\(\mathcal{T}(s, a)\)) are discrete and deterministic, (ii) the cost function prioritizes goal achievement (typically $0$ if \(x \models \mathcal{G}\) otherwise \(\infty\)), and (iii) the optimization horizon spans until the goal is reached.


\subsection{Planning extensions (Multicriteria and Online Planning)}

As the classical planning lacks the expressiveness needed for real-world dynamic systems, there are also extensions such as \textit{numeric}, \textit{temporal}, and \textit{stochastic} planning that bring it a bit closer to the dynamic, real-time nature of MPC. \textit{Numeric planning} incorporates continuous state variables and numeric fluents, partially enabling modeling dynamic systems akin to MPC’s difference equations. \textit{Temporal planning} introduces explicit reasoning about time, including action durations and deadlines. \textit{Stochastic planning} then addresses uncertainty by modeling probabilistic (or worst-case) transitions.

Also, classical PDDL (Planning Domain Definition Language)~\cite{lavalle2006planning} planning typically involves finding a sequence of actions \( \pi = \langle a_1, a_2, \dots, a_T \rangle \) that transitions the system from an initial state \( s_0 \) to a goal state \( s_n \models \phi_G \). Extensions following the numeric and optimal planning then led to introducing \textit{plan metrics} \( F(\pi) \), functions that evaluate the quality of a plan based on numeric fluents \( \{f_1, f_2, \dots, f_k\} \) in the domain (e.g., fuel, distance, etc.). These metrics allow differentiating between plans that achieve the same goal \( G \) but with varying numeric performance. For example, a plan \( \pi \) is optimal if \( \pi^* = \arg\min_\pi F(\pi) \), where \( F(\pi) \) might minimize action cost, time, or resource usage. The metric \( F \) can combine the multiple criteria through functions over fluents, such as \( F(\pi) = \sum_{i=1}^k w_i F_i(\pi) \), yielding context-dependent optimal plans. This approach, more reminiscent of the function optimization in MPC and RL, is supported in extended PDDL versions, however, remains rarely supported by actual planners.

Further, classical PDDL planning operates under the assumption of a static, fully observable environment, generating a complete plan \( \pi = \langle a_1, a_2, \dots, a_T \rangle \) offline to transition from an initial state \( s_0 \) to a goal state \( s_g \), where \( s_g \models G \) for a fixed goal condition \( G \). However, real-world environments often exhibit and unpredictable behavior, which can be addressed by interleaving planning and execution to achieve a form of ``online planning.'' This can be modeled through iterative \textit{replanning}, where the system updates the problem state \( s_t \) to reflect the current environment and recomputes the optimal plan \( \pi^* = \arg\min_\pi F(\pi \mid s_t) \) from the current \( s_t \). 
A possible hybrid approach is then to combine PDDL with inherently dynamic frameworks like the MPC or RL, where PDDL handles high-level planning, and the dynamic framework manages the uncertainties by optimizing policies \( \pi^* \) over the underlying state transitions \( P(s_{t+1} \mid s_t, a_t) \) in real-time.

\subsection{Planning (and RL) with (Graph and) Relational Neural Networks}

Classical planning has been extended in recent years to attempt to leverage the computational and representational power of neural modeling and training. As in MPC, the techniques of RL can be used to learn both the functions defining the classical planning problem, the transition and cost, through data, or even be used to compute an optimal policy. Before continuing to describe some of the main methodological tools available, we also indicate another contemporary AI tool component: the use of structured data representations. By using Graphs, and more generally structured relations, we can endow an inductive bias into a wide variety of use cases that facilitates computation as well as interpretability. 

Recent advancements in structured deep learning models, particularly Graph Neural Networks (GNNs)~\cite{zang2023graphmp}, have proven effective for classical planning tasks, which can often be naturally encoded in the form of various graphs. Given a graph \( G = (V, E) \), where \( V \) represents entities, such as domain objects, and \( E \) represents relationships between them, a GNN learns the embeddings of nodes \( \mathbf{v}_i^{(K)} \) by iteratively passing messages between nodes. The message-passing update rule is generally given by:

\begin{equation}
\mathbf{v}_i^{(k+1)} = \phi \left( \mathbf{v}_i^{(k)}, \bigoplus_{j \in \mathcal{N}(i)} \psi(\mathbf{v}_j^{(k)}, \mathbf{e}_{ij}) \right),
\end{equation}
where \( \mathcal{N}(i) \) denotes the neighbors of node \( i \), \( \mathbf{e}_{ij} \) is the, possible feature-associated, edge between nodes \( i \) and \( j \), \( \psi \) is the message function, and \( \phi \) updates the node embedding. After \( K \) iterations, the learned node embeddings \( \mathbf{v}_i^{(K)} \) are aggregated and used to predict heuristics \( h(s) \), which are then used in search algorithms. The most prominent in planning, i.e. searching for a sequence of actions that leads from an initial state \( s_0 \) to a goal state \( s_g \), is the classic A* search that maintains a priority queue of states \( s \), ordered by the evaluation function \( f(s) \), which is defined as \( f(s) = g(s) + h(s)\),
where \( g(s) \) is the known cost to reach state \( s \) from the initial state \( s_0 \), and \( h(s) \) is the heuristic estimate of the cost to reach the goal state \( s_g \). At each step, A* selects the state with the lowest \( f(s) \), expands its neighbors, and updates their \( f \)-values accordingly. 
With enough data samples, the GNN can often learn a heuristic \( h(s) \) based on the graph structure of the planning problem, improving the A* search efficiency by providing more accurate cost-to-go estimates, leading to fewer expanded nodes and faster planning. Structured neural networks can also be applied to learn the integration of scheduling and planning~\cite{ma2020online}

\paragraph{Relational NNs}
However, standard GNNs can not integrate domain's \textit{background knowledge} (BK)  directly into the learned model heuristics. 
To address this, one may use more expressive ``neuro-symbolic'' frameworks based on the language of relational logic (Sec.~\ref{sec:logic}). One such prominent example framework are Lifted Relational Neural Networks (LRNNs)~\cite{sourek2018lifted}, which extend GNNs by incorporating \textit{relational rules} that guide the learning process. Similarly to RL, a planning policy \( \pi(s) \) can be a (learned) model that selects the (optimal) action \( a \) given a state \( s \), i.e.,

\begin{equation}
\pi(s) = \arg\max_{a \in A} C(a \mid s),
\end{equation}
where \( A \) is the set of \textit{all} possible actions, and \( C(a \mid s) \) is some (learned) scoring function. 
With LRNNs, however, domain knowledge can also be encoded, e.g. in the form of a prior nondeterministic policy \( \sigma(s) \) that returns merely a subset of feasible actions \( \sigma(s) \subseteq A \). Similarly to the T-S fuzzy control formalism (Sec.~\ref{sec:fuzzy}), this policy can be defined by rules of the form:

\begin{equation}
\text{condition}(s) \to \text{action}(a),
\end{equation}
where \( \text{condition}(s) \) is a \textit{relational} logic formula over the state \( s \), and \( \text{action}(a) \) specifies the feasible actions according to the domain knowledge capturing some desired domain constraints, such as safety. This nondeterministic policy \( \sigma(s) \) then effectively restricts the learning hypothesis space of the model \( \pi(s) \). For instance, in the popular Blocksworld domain~\cite{lavalle2006planning}, \( \sigma(s) \) may specify that misplaced blocks are unstacked onto the table before restacking them in the correct order, which ensures a \textit{satisficing} but potentially suboptimal policy.
The policy \( \pi(s) \) can then be further refined by optimizing the scoring function \( F(a \mid s) \) that ranks actions in \( \sigma(s) \) based on their likelihood of improving the given plan quality metrics, such as minimizing cost. This allows the model to focus on optimizing the solution quality while adhering to the given BK constraints, instead of solving the planning problem from scratch in potentially unsafe ways.

In general, integrating logical (relational) background knowledge can be used for domain-guided planning, where learned policies \( \pi(s) \) are shaped by both training data and expert-provided domain knowledge. This offers advantages over existing, e.g. GNN-based, learning-to-plan methods by ensuring that the learned policies are both satisficing and aligned with real-world constraints.

\paragraph{Reinforcement Learning with Logic and Relational Policies}

Related ideas have also been explored in the reinforcement learning setting with other neuro-symbolic frameworks, such as $\delta$ILP~\cite{evans2018learning}, which aims to relax the problem of Inductive Logic Programming (ILP), or rule learning, into a differentiable setting via (constrained) rule templates where each rule is associated with a learnable weight, similar to the LRNNs, and trained with gradient descent.
Building on that, frameworks such as~\cite{jiang2019neural} integrate $\delta$ILP with RL to enhance interpretability and generalization through the logical formalism. The RL problem is modeled as an MDP with logic components \((\mathcal{M}, p_S, p_A)\), where \(\mathcal{M} = (\mathcal{S}, \mathcal{A}, T, R)\) defines the state space, action space, transitions, and rewards. A state encoder \(p_S : \mathcal{S} \to 2^\mathcal{G}\) maps states to sets of logical atoms, and an action decoder \(p_A : [0,1]^{|\mathcal{D}|} \to [0,1]^{|\mathcal{A}|}\) translates the logic valuations into action probabilities.
Using $\delta$ILP, the policies are represented as valuation vectors \(e \in [0, 1]^{|\mathcal{G}|}\), where each element indicates the confidence of a logical atom being true. Inference then occurs via a recurrent neural model \(f_\theta\), which iteratively updates the valuations through time. Logical rules are encoded as weighted first-order clauses of the form 
\[
\alpha \leftarrow \alpha_1, \ldots, \alpha_n,
\]
where \(\alpha\) is the head atom, and \(\alpha_1, \ldots, \alpha_n\) are body atoms, again similarly to LRNN or the T-S fuzzy systems. Trainable clause weights, constrained by a softmax function, then ensure differentiability, making the policy
\[
\pi(s) = p_A(f_\theta(p_S(s))),
\]
combining logic reasoning with standard RL optimization using policy gradients to maximize the expected returns.

\subsection{Learning for Planning}
The previous section described how to integrate the crisp structure planning provides as far as representations of policies into an otherwise data-driven scheme for performing Reinforcement Learning. 

Here, consider now the converse, that one has a planning problem, and because of unknown functions, learning must be performed. We can see the resemblance to Dynamic Programming again as with MPC. Rather, planning shared the discrete action and state space of classical DP, unlike MPC. However, although an objective function can be present in both, the presence and moreover priority of the goal state is particular to planning. 

Consider again the components of the planning problem functions. Recall the equations associated with the planning problem:
\[
s_n \models \phi_G, \quad \text{where } s_{t+1} = \psi(a_t, s_t), \quad t = 0, \dots, n-1.
\]

\[
\mathbf{a}^* = \arg\min_{\mathbf{a}} C(\mathbf{a}), \quad C(\mathbf{a}) = \sum_{t=0}^{n-1} C(s_t, a_t) + C_T(s_n,n),
\]

Now, consider that, unlike in classical planning, the functions $\psi$ and $F$ are not known, but must be learned. This can be done by observing the sequences of real states $\{\hat{s}_1,\hat{s}_2,\cdots,\hat{s}_n\}$ as well as obtaining rewards $R_t$ during the process. For learning, parametrization of the problem functions can be performed,
\[
\psi(a_t,s_t):=\psi(\theta^{\psi},a_t,s_t),\,
F(a_t,s_t):=F(\theta^{F},a_t,s_t),\,F_T(s_n,n):=F_T(\theta^{T},s_n,n)
\]
and these functions can be learned directly, which would correspond to a model based approach of reinforcement learning. 

Alternatively, a purely model-free approach can solve for the solution of the planning problem without learning the components of the function. In the standard way, an action value function $V(s_t,a_t)$ can be directly learned, or simply a policy gradient approach used, to compute and define a policy 
\[
\pi(s) = p_A(f_\theta(p_S(s))),
\]
as described above. In this case, however, the original structure of the planning problem is hardly maintained, i.e., this approach is barely distinguishable from standard RL. 

Later in this work, we do consider uncertain realizations of the dynamics and rewards in the context of care planning. However, we will use a lower level state representation, describing detailed physical dynamics, and run by an MPC-based solver, to assist in, ultimately more accurate and robust, data-driven descriptions of the planning operation and its solution.

\section{Fuzzy Sets and Fuzzy Control - Background}\label{sec:fuzzy}

Fuzzy sets and fuzzy control provide a mathematical framework for dealing with imprecision and vagueness. Traditional set theory, based on binary membership where an element either belongs to a set or not, is often insufficient for capturing the vagueness inherent in many real-world systems. Fuzzy set theory, introduced by~\cite{zadeh1965fuzzy}, extends the classical sets by allowing elements to have varying degrees of \textit{membership}, thereby modeling the uncertainty and imprecision.

\subsubsection{Fuzzy Sets}

A fuzzy set $\tilde{A}$ in some universal set of discourse $X$ is defined by its \textit{membership} function $\mu_{\tilde{A}} : X \to [0, 1]$, where $\mu_{\tilde{A}}(x)$ indicates the degree to which an element $x \in X$ belongs to $\tilde{A}$. Formally, a fuzzy set can be represented as:

\begin{equation}
    \tilde{A} = \{ (x, \mu_{\tilde{A}}(x)) \mid x \in X \},
\end{equation}
where $\mu_{\tilde{A}}(x) = 1$ indicates full membership, $\mu_{\tilde{A}}(x) = 0$ indicates non-membership, and $0 < \mu_{\tilde{A}}(x) < 1$ represents a partial membership.
Common examples of membership functions include triangular, trapezoidal, Gaussian, and sigmoid functions~\cite{bai2006fundamentals}.

Basic set operations on these fuzzy sets are then commonly defined as 

\begin{itemize}
    \item {Union:} $\mu_{\tilde{A} \cup \tilde{B}}(x) = \max(\mu_{\tilde{A}}(x), \mu_{\tilde{B}}(x))$
    \item {Intersection:} $\mu_{\tilde{A} \cap \tilde{B}}(x) = \min(\mu_{\tilde{A}}(x), \mu_{\tilde{B}}(x))$
    \item {Complement:} $\mu_{\tilde{A}^c}(x) = 1 - \mu_{\tilde{A}}(x)$
\end{itemize}

In practice, fuzzy sets are then used to model ``linguistic variables'', where terms such as ``low, medium, or high'' are represented by fuzzy sets with appropriately defined membership functions. For example, in a temperature control system, the variable ``temperature'' could be modeled using fuzzy sets such as $\text{Low}$, $\text{Medium}$, and $\text{High}$.

\subsubsection{Fuzzy Logic}

Fuzzy logic then extends classical Boolean logic by incorporating \textit{degrees of truth}. In fuzzy logic, the truth value of a statement is again not restricted to binary values ($0$ or $1$) but can take any value in the interval $[0, 1]$.
To that aim, fuzzy logic operators extend classical logical operators to operate on the fuzzy sets~\cite{bai2006fundamentals}. 

These fuzzy logic operators are then typically defined using t-norms $T: [0, 1]^2 \rightarrow [0, 1]$ (for conjunction, satisfying commutativity, associativity, monotonicity, and having 1 as identity) and t-conorms $S: [0, 1]^2 \rightarrow [0, 1]$ (for disjunction, with analogous properties). 
Common choices include: minimum t-norm ($T(a, b) = \min(a, b)$), product t-norm ($T(a, b) = a \cdot b$), Łukasiewicz t-norm ($T(a, b) = \max(0, a + b - 1)$), maximum t-conorm ($S(a, b) = \max(a, b)$), probabilistic sum t-conorm ($S(a, b) = a + b - a \cdot b$), and Łukasiewicz t-conorm ($S(a, b) = \min(1, a + b)$). 
Negation is then commonly defined as $\neg a = 1 - a$, although there are different choices as well.

\subsubsection{Fuzzy control}


Fuzzy control utilizes fuzzy set theory and fuzzy logic to design \textit{control} systems~\cite{passino1998fuzzy}. A fuzzy controller typically consists of four main components. First, \textit{Fuzzification} involves converting crisp input values into fuzzy sets using the membership functions. Second, \textit{Inference} applies fuzzy rules to evaluate outputs based on fuzzy inputs. Third, \textit{Aggregation} combines the results from all applicable rules to form a single fuzzy output set and, fourth, \textit{Defuzzification} converts the fuzzy output set into a crisp output value.  

The fuzzy rules are typically expressed in the form of ``IF-THEN'' statements, where the antecedent (IF part) describes the input conditions using fuzzy sets, and the consequent (THEN part) determines the output action.  

\paragraph{Mamdani system}
The Mamdani fuzzy inference system~\cite{mamdani1975experiment} is one of the earliest and most widely used fuzzy control methodologies. It is based on the idea of representing both antecedents and consequents as fuzzy sets. A typical rule in the Mamdani model is expressed as:  
\begin{equation}
    \text{IF } x_1 \text{ is } \tilde{A}_1 \text{ AND } x_2 \text{ is } \tilde{A}_2 \text{ THEN } y \text{ is } \tilde{B},
\end{equation}
where $x_1, x_2$ are input variables, $\tilde{A}_1, \tilde{A}_2$ are fuzzy sets for the antecedents, and $\tilde{B}$ is the fuzzy set describing the consequent.  

In Mamdani systems, the inference process proceeds as follows. For each rule, the degree of activation is computed by evaluating the membership functions of the antecedents. Commonly, the \textit{min} operator is used to compute the degree of fulfillment for the conjunction (t-norm) of antecedents (joined by AND), i.e.,  
\begin{equation}
    \mu_{\tilde{R}_i}(x_1, x_2) = \min(\mu_{\tilde{A}_{i1}}(x_1), \mu_{\tilde{A}_{i2}}(x_2)).
\end{equation}
The output fuzzy set $\tilde{B}_i$ of each rule is then scaled by the degree of activation $\mu_{\tilde{R}_i}$, e.g. using again the \textit{min} operator:  
\begin{equation}
    \mu_{\tilde{B}_i}(y) = \min(\mu_{\tilde{B}}(y), \mu_{\tilde{R}_i}).
\end{equation}
The scaled output fuzzy sets $\tilde{B}_i$ from all rules are then combined, typically using the \textit{max} operator:  
\begin{equation}
    \mu_{\tilde{B}_{\text{agg}}}(y) = \max_i \mu_{\tilde{B}_i}(y).
\end{equation}
Finally, the aggregated fuzzy output set $\tilde{B}_{\text{agg}}$ is converted into a crisp value through defuzzification. A common method is the centroid method, where the crisp output $y^*$ is computed as:  
\begin{equation}
    y^* = \frac{\int y \cdot \mu_{\tilde{B}_{\text{agg}}}(y) \, dy}{\int \mu_{\tilde{B}_{\text{agg}}}(y) \, dy}.
\end{equation}  

The Mamdani approach is particularly intuitive and well-suited for human-interpretable control systems, where the rules are derived from expert knowledge.  

\paragraph{Takagi-Sugeno system}
The Takagi-Sugeno (T-S) fuzzy model~\cite{takagi1985fuzzy} is the other prominent formalism that takes a different approach by employing crisp mathematical functions in the consequents of the rules. A typical T-S rule is expressed as:  
\begin{equation}
    \text{IF } x_1 \text{ is } \tilde{A}_1 \text{ AND } x_2 \text{ is } \tilde{A}_2 \text{ THEN } y = g(x_1, x_2),
\end{equation}
where $g(x_1, x_2)$ is typically a linear function:  
\begin{equation}
    g(x_1, x_2) = a_1 x_1 + a_2 x_2 + b,
\end{equation}
with $a_1, a_2, b \in \mathbb{R}$ as parameters.  

In the T-S model, the inference process begins by computing the degree of activation $\mu_{\tilde{R}_i}(x_1, x_2)$ for each rule, as in the Mamdani system. The weighted output for each rule is then calculated based on its activation degree:  
\begin{equation}
    y_i = \mu_{\tilde{R}_i}(x_1, x_2) \cdot g(x_1, x_2).
\end{equation}
The outputs of all rules are combined using a weighted average:  
\begin{equation}
    y^* = \frac{\sum_{i=1}^n \mu_{\tilde{R}_i}(x_1, x_2) \cdot g_i(x_1, x_2)}{\sum_{i=1}^n \mu_{\tilde{R}_i}(x_1, x_2)}.
\end{equation}  

The T-S model is computationally efficient and well-suited for control problems where a mathematical model of the system dynamics is available.  



\section{Integrating the Symbolic and Control Layers}\label{sec:intergrating}

Consider now the robotic aide described in the Introduction of the manuscript. Here we describe an integration of the low level control operation with the higher level conceptual discrete planning and scheduling decision making.

Figure~\ref{fig:schema} presents the high level picture of the operation, with Figure~\ref{fig:integration_timescales} detailing the operations across different time-scale granularities. Over a long, slow time scale, a Scheduler decides upon a sequence of tasks. When the time comes for each task, a planning problem corresponding to that task is formulated and solved. In the course of executing the solution of the planning problem, each action taken generates a continuous control problem that is solved at a fast time scale that yields the true reward and next state, which may or may not be the one predicted by the planner. Models of rewards and transitions are learned, with information flowing across layers in a continuous, sequential manner.

\begin{figure}[t]
\centering
\begin{adjustbox}{width=\textwidth}
\begin{tikzpicture}[node distance=1cm and 1cm,
    modblock/.style={rectangle, draw, thick, fill=white, drop shadow, text width=2.8cm, align=center, minimum height=1.5cm, font=\small},
    bus/.style={rectangle, draw, dashed, thick, fill=yellow!10, minimum height=1cm, minimum width=14cm, rounded corners}
    ]
    \node [modblock, fill=blue!5] (sched) at (0,0) {\textbf{Scheduler} \\ \tikz\pic[scale=0.5]{clockface};};
    \node [modblock, fill=orange!5] (plan) at (4,0) {\textbf{Planner} \\ \tikz\pic[scale=0.5]{treegraph};};
    \node [modblock, fill=green!5] (ctrl) at (8,0) {\textbf{MPC} \\ \tikz\pic[scale=0.5]{miniplot};};
    \node [modblock, fill=gray!20] (sys) at (12,0) {\textbf{Robot} \\ \tikz\pic[scale=0.5]{robot};};
    
    \node [bus] (learn) at (6,-2.5) {
        \textbf{Learning Updates} \hspace{1cm} 
        Parameters: $\theta^\psi$ (Transition), $\theta^C$ (Cost), $\theta^V$ (Value)
    };
    
    \draw [->, thick] (sched) -- node[above, font=\footnotesize] {$T$} (plan);
    \draw [->, thick] (plan) -- node[above, font=\footnotesize] {$p(a_t)$} (ctrl);
    \draw [->, thick] (ctrl) -- node[above, font=\footnotesize] {$\hat{v}$} (sys);
    
    \draw [->, dashed, blue, thick] (sys.south) -- node[right, font=\footnotesize] {Data} (sys.south |- learn.north);
    
    \draw [->, dashed, blue] (learn.north -| sched.south) -- (sched.south);
    \draw [->, dashed, blue] (learn.north -| plan.south) -- (plan.south);
    \draw [->, dashed, blue] (learn.north -| ctrl.south) -- (ctrl.south);
    
    \draw [->, dashed, red] (sys.north) -- ++(0,0.5) -| node[pos=0.25, above=0.2] {Feedback Loop} (ctrl.north);
    \draw [->, dashed, blue] (sys.north) -- ++(0,0.7) -| node[pos=0.25, above] {} (sched.north);
\end{tikzpicture}
\end{adjustbox}
\caption{\rev{High-level Integration Schema: A learning module (bottom) collects data from the robot to update parameters of the decision-making modules (top).}}
\label{fig:schema}
\end{figure} 
\begin{figure}[t]
\centering
\resizebox{\textwidth}{!}{%
\begin{tikzpicture}[
  module/.style={rectangle, draw, thick, fill=white, drop shadow, rounded corners, align=center, font=\small, minimum height=1.2cm},
  txt/.style={font=\small},
  arrow/.style={->, thick, -{Stealth[scale=1.2]}},
  dashedarrow/.style={->, thick, dashed, -{Stealth[scale=1.2]}}
]

  \draw[-{Stealth[scale=1.2]}, line width=1.5pt] (0, 6) -- (14.5, 6) node[above left, font=\bfseries] {Time};

  \draw[dashed, thick, gray] (7.5, 6.5) -- (7.5, -2.4) node[below, font=\bfseries, text=black] {Current Time $t$};

  \node[module, fill=blue!5] (sched) at (9, 5) {\textbf{Scheduler} \\ \tikz\pic[scale=0.4]{clockface};};
  \node[module, fill=gray!10, minimum height=0.6cm] (taskB) at (13, 4.5) {Task B};
  
  \draw[dashedarrow] (sched.east) -- (taskB.west);

  \node[module, fill=orange!2, minimum height=0.6cm] (taskA) at (4, 4.5) {Task A};
  \node[module, fill=orange!5] (planA) at (4, 3.0) {\textbf{Planning Problem A} \\ \footnotesize $\phi_G, \psi, C, C_T$ \\ \tikz\pic[scale=0.4]{treegraph};};
  
  \draw[dashedarrow] (sched.west) -- (taskA.east);
  \draw[arrow] (taskA.south) -- (planA.north);

  \node[txt] (s_past) at (1, 1.5) {$\hat{s}_{t-1}, \hat{a}_{t-1}$};
  \node[txt] (s_curr) at (4, 1.5) {$\hat{s}_t, \hat{a}_t$};
  \node[txt] (s_next) at (9.5, 1.5) {$\bar{s}_{t+1}, \bar{a}_{t+1}$};
  \node[txt] (s_goal) at (13.5, 1.5) {$\bar{s}_{T} \in \phi_G$};

  \draw[arrow] (0, 1.5) -- (s_past);
  \draw[arrow] (s_past) -- (s_curr);
  \draw[dashedarrow] (s_curr) -- (s_next);
  \draw[dashedarrow] (s_next) -- (s_goal);

  \node[module, fill=green!5] (mpc) at (4, -0.5) {\textbf{MPC} \\ \footnotesize $W, G, f, h, \tilde{A}$ \\ \tikz\pic[scale=0.4]{miniplot};};
  \draw[arrow] (s_curr) -- node[left, font=\footnotesize] {Trigger} (mpc);

  \node[txt] (z_curr) at (3, -2.) {$\hat{z}_{\tau(0)}, \hat{v}_{\tau(0)}$};
  \node[txt] (z_mid)  at (5.5, -2.) {$\hat{z}_{\tau(1)}, \hat{v}_{\tau(1)}$};
  \node[txt] (z_dots) at (7.5, -2.) {$\dots$};
  \node[txt] (z_next) at (9.5, -2.) {$\hat{z}_{\tau(H)}$}; 

  \draw[arrow] (1.5, -2.) -- (z_curr);
  \draw[arrow] (z_curr) -- (z_mid);
  \draw[arrow] (z_mid) -- (z_dots);
  \draw[arrow] (z_dots) -- (z_next);

  \draw[arrow] (z_next) to[out=100, in=-100] node[right, font=\footnotesize, align=left] {(De-)Fuzzify ($\mu$) \\ Aggregate ($T^s$)} (s_next);

  \draw[dashedarrow, red] (9.5, 1.8) .. controls (9.5, 2.5) and (6.5, 2.5) .. node[above, font=\footnotesize, align=center, text=red] {State Discrepancy \\ $\to$ Re-plan} (planA.east);

\end{tikzpicture}%
}
\caption{A temporal model of integrating the Planning and Control Layers. The slow-timescale Scheduler assigns tasks (e.g., Task A) which spawn Planning Problem(s). The discrete Planner issues actions $\hat{a}_t$ that trigger the fast-timescale MPC. Continuous physical states $\hat{z}$ map back to discrete symbolic states $\bar{s}$ via (de-)fuzzification ($\mu$). Possible discrepancies prompt real-time re-planning.}
\label{fig:integration_timescales}
\end{figure} 

\subsection{Hierarchical Planning-MPC Formulation}\label{sec:hierarchial_formulation}


Consider the original formulation of the planning problem, wherein a set of actions $a(t)\in\mathcal{A}$ with $\left\vert\mathcal{A}\right\vert = A$ must be chosen in order for the state $s_t\in\mathcal{S}$ with $\left\vert\mathcal{S}\right\vert = S$ to proceed in order to satisfy a goal state $x(T)\in\mathcal{G}$ while minimizing some cost. The process is defined as a stochastic planning problem, that is, now the cost functions and state transitions include a random variable. For a contemporary example, see~\cite{flores2020stochastic}. While the challenging field is growing, it is significantly
less mature than classical planning itself. Moreover, we will assume limited probabilistic knowledge a priori as to the distributions of these functions, and instead focus on a 
data driven approach for learning these functions. Specifically, we will use parametric function approximations for defining the expectations of the costs.  

\begin{equation}
\begin{array}{l}
 C(\mathbf{s},\mathbf{a}) = \sum_{t=0}^{n-1} C(s_t, a_t) + C_T(s_n)+\rho(n),\\
 C(s,a) = \mathbb{E}_{\xi}[c(s, a,\xi(s,a))\\
 C_T(s) = \mathbb{E}_{\xi}[c_T(s_n,\xi(s_n,a_n))]
 \end{array}
\end{equation}
where $\rho(n)$ is monotonically increasing with respect to $n$ and is a time penalty of the number of steps until the goal state is reached. The cost, which will depend on some metric of the overall well-being of the human and agreeableness of the robotic service, is treated as stochastic. Note that the randomness itself depends on the state and action taken, modeling the potential variability of, e.g., emotional affect in response to a robot's action across different circumstances. 

Consider that the transitions are governed by the \emph{deterministic} but \emph{unknown} dynamics
\[
s_{t+1}=\psi(s_t(z_t),a_t,z_t,e_t)
\]
where $z_t$ is a continuous state of the robot and $e_t$ is an exogenous set of quantities that affect the transition process. Observe that the complete continuous state $z_t$ completely identifies the discrete symbolic state $s_t$ and $e_t$ is unpredictable and only partially observed. Thus, the actual state transitions are unidentifiable, and the planning computation is performed with a surrogate. The fact that these plans are executed multiple defines in an online sequence set by a scheduler makes the problem tractable. That is, with a one off planning problem, one clearly cannot learn the terminal state cost. However, since each plan is expected to be repeated and a scheduler will define tasks into an indefinite future time period, these function approximations are expected to become increasingly more accurate over time.

This surrogate will be defined by function approximations, introducing parameters $\theta^C$ and $\theta^T$ for the cost functions and $\theta^{\psi}$ for the state transition. 
Specifically, the functional forms for the model the operator uses to compute the plan as:
\begin{equation}\label{eq:planfuns}
\begin{array}{l}
\bar{C}(s_t,a_t) = \Phi^C(\theta^C;s_t,a_t) \\
\bar{C}_T(s_n) = \Phi^T(\theta^T;s_n) \\
\bar{s}_{t+1} = \bar{\psi}(s_t,a_t)=\max\limits_{s\in\mathcal{S}} p_{\psi}^{S}(\theta^{\psi};s_t,a_t) 
\end{array}
\end{equation}
that is, the cost function is some parametrized model with inputs corresponding to the action, previous and subsequent states, and intermediate process, and the state transition is predicted $\bar{s}_{t+1}$ by a probability density also depending on these quantities. While stochastic planning is a possible approach to handling uncertainty, we prefer to use the more comprehensive and reliable toolkit of classical planning algorithms and simply choose perform a softmax of the distribution of the predicted state. 

We will denote the sequence $\{\hat{s}_t,\hat{a}_t\}$ to be the physically realized action and state pair that the procedure executes and transitions through. In addition, at current time $\hat{t}$, the intermediate sequence of reward has been observed as $\{\hat{C}_0,\hat{C}_1,\cdots,\hat{C}_{\hat{t}-1}\}$. We consider that the same planning problem will be solved multiple times by the robot and so historical rewards and states are also available and have informed the models $\bar{C}$, $\bar{C}_T$ and $p_{\psi}^S$.

The lower level control process is defined as follows: in between each transition $t$ to $t+1$, a state-control process running at a faster time scale $\{z_{\tau},v_{\tau}\}$ with $\tau$ ranging from $t,t+h,t+2h,\cdots,t+1$ with $h\ll 1$. We remark that while the planning state and action $(s_t,a_t)$ are members of a discrete set, the lower level control operation $\{z_{\tau},v_{\tau}\}$ is defined by continuous vectors. Naturally, there is a correspondence between the continuous and discrete state, which arises in the overall operation loop when the next state $\hat{s}_{t+1}$ is inferred from the final state in the control problem $\hat{z}_{t+1}$, which we define later. As a simple example, $s_t$ can indicate which out of a finite set of rooms a robot is in, and $z_{\tau}$ the precise coordinates of its location.

We note that nonuniform and adaptive step-sizes, multiple shooting formulations, etc. can be appropriately defined with natural extensions of the formulation.

Below we proceed to detail the functional components of our procedure of integrating the planning, MPC, and learning aspects of the procedure.

\subsection{Selecting the planning action}\label{sec:selecting_action}
We consider real time closed loop planning, in which not only is a plan computed at the initial moment of operation, but for each step in the operating sequence, the plan is recomputed. This is necessary because of the model of the functions used for planning, that is $\Phi^V$ and $p^S_{\psi}$ are updated at each step to reflect discrepancy in the expected and real reward, that is $\bar{C}(s_t,a_t)$ compared to $\hat{C}_t$, as well as the transition model by $\bar{s}_{t+1}$, the predicted state, relative to $\hat{s}_{t+1}$, the actual realized state. For the rest of this section, consider a consider that the current time of real-time online operation is given by $\hat{t}$.  

In this section we present several options with respect to how one computes a plan at each time $\hat{t}$. To this end we assume that at the current time $\hat{t}$ all of the model functions are current, i.e., have been updated with the most recent observations. 

At this stage the operation has updated information that may affect the optimal plan that was computed at $\hat{t}-1$. This necessitates that the planner should then recompute the plan for the remaining time $\{\hat{a}_{\hat{t}}, \bar{s}_{\hat{t}+1},\bar{a}_{\hat{t}+1},
\bar{s}_{\hat{t}+2}, \bar{a}_{\hat{s}+2},\cdots \bar{s}_n\}$. Here the quantities $\bar{s}$ and $\bar{a}$ indicated the predicted future states and the computed actions for these predicted states. The action chosen to be implemented once $\hat{t}\ge t$ is given as $\hat{a}_t$. For the choice of algorithm to perform this computation, there are several options which will be outlined below.



\paragraph{Actions as RL Policies}  
A classic RL approach would correspond to the choice $\hat{a}_t$ being computed in accordance with some action-value function $V(\theta^V,a,s)$. This is meant to represent a cost-to-go, that is,
\begin{equation}\label{eq:actionval}
    V(\theta^V,a_t,s_t) = \max\limits_{a\in\mathcal{A}}\left\{
    \begin{array}{ll} & -\bar{C}(s_t,a)-\sum\limits_{t=\hat{t}+1}^{n-1} \bar{C}(\bar{s}_t,\bar{a}_t)-C_T(\bar{s}_n)-\rho(n) \\
    \text{such that }& \bar{s}_n\in \mathcal{G},\\
    & \bar{s}_{t+1} = \bar{\psi}(\bar{s}_t,\bar{a}_t),\,t>\hat{t} 
    \end{array}\right\}
\end{equation}
Observe that learning the intermediate and terminal cost functions czn be used to define the action-value function, or $\theta^V$ can be used to approximate the model on the right, eschewing the need for $\theta^C$, etc.

The standard techniques here correspond to greedy, $\epsilon$-random greedy, and logit:
\begin{equation}\label{eq:planpolicy}
\begin{array}{lr}
\hat{a}_t=\hat{a}_t^{g}=\arg\max_{a} V(\theta^V_t,a,s_t) & \text{greedy}\\
\hat{a}_t=\hat{a}^{\epsilon}_t=\left\{\begin{array}{lr} \arg\max_a V(\theta^V_t,a,s_t) & \text{w.p.} (1-\epsilon) \\
a' & \text{w.p.}\frac{\epsilon}{A},\,\,\forall a'\in\mathcal{A}\end{array}\right. & \epsilon-\text{greedy}\\
\hat{a}_t=\hat{a}^{\sigma}_t=\hat{a}(\omega^a_t)\sim_a e^{- V(\theta^V_t,a,s_t)}/Z & \text{logit}
\end{array}
\end{equation}

\paragraph{Classical Planning on the Learned Model}

Alternatively, some classical planning algorithm $\mathcal{C}$ can perform a selection of $\hat{a}_t$ based on purely logical representations, treating the learned models as fixed exogenously at the point of computation:
\[
    \hat{a}^c_t = \mathcal{C}(\bar{\psi},\bar{C},\bar{C}_T,\rho,s_t),
\]
where $\mathcal{C}$ may be a learned policy, that is a neural map from the parametrization of the models to the action sequence, or a search procedure, driven by the respective goal satisfaction and cost minimization. Various exact and heuristic search methods, such as the well known A$*$ procedure, can be used to search for a policy. 

This approach is most analogous to the MPC-RL integration described above, as far as the use of classical numerical procedures with firm theoretical guarantees being applied to compute the solution, while neural networks are used to parametrize and model the component objective, or value, functions. 


\subsection{Action-Triggered MPC}\label{sec:action_triggered_mpc}
Consider a real time operation of planning with $\hat{t}$ denoting the current time. A plan $\zeta(s)$ that at time $\hat{t}$ with $\hat{t}>0$ will execute a particular action $\hat{a}_t{\hat{t}}$. This action triggers the solution and execution of an MPC problem. The action defines the problem functions, that is the dynamics, objective, acceptable continuous actions, and any constraints, and the MPC problem is solved in real-time closed loop fashion, delivering a control sequence. This control sequence transitions the continuous state until the terminal state, which is then mapped onto an appropriate discrete state. In addition, any additional reward-relevant information that appeared during the course of the control solution is passed to the planner to assist with updating a more faithful plan at the next discrete time instance.

Formally, the MPC final state is expected to yield, upon its mapping onto the discrete set $\mathcal{S}$, the predicted state that was associated with the choice of action $a_{\hat{t}}$. This is based on the model at time $\hat{t}$, $p^{S}_{\psi}(\theta^{\psi};s_t,a_t)$, this defines the predicted state $\bar{s}_{t+1}$ at the next time step. The details of how $\theta^{\psi}$ is learned and $\bar{s}_{t+1}$ is computed will be the topic of subsequent subsections. 

We define the action of the realized $\hat{a}_t$ to be the formulation and approximate solution of a discretized optimal control problem. The form of the objective $G$ and $W$, dynamic state equation $f$ and constraints $h$ defined in the NLP below depend on the state $\hat{s}_t$ and action $\hat{a}_t$ taken. The initial continuous state is observed $\hat{z}_t$ and it is expected that the final state cost $W$ should deliver a state $z_{t+1}$ that maps on to $\bar{s}_{t+1}$.

The problem is solved in a standard rolling-horizon MPC fashion, that is, sequentially wherein after the first time the solution is computed, the optimal control decision for the first time instant is performed and the state is allowed to evolve, and when the resulting intermediate state has been observed, the process is repeated with the new horizon now shorter by one time step.

We now proceed to describe the components of this operation formally. Real time optimization is performed in the fast time scale, $\{\tau(0),\tau(1),\cdots,\tau(H-1)\}:=\{t,t+h,t+2h,\cdots,t+1-h\}$ with the total number of time steps $H = h^{-1}$. 

We definite parameters for constructing the objective and constraint functions in the MPC problem definition as depending on the action to be 

\[
\left\{p^w(a_t),p^g(a_t),p^f(a_t),p^h(a_t),p^a(a_t)\right\}
\]
The initial physical state $\hat{z}_t$ is observed and the defined optimal control problem is solved below. Note that we can also prescribe a final state set if we want a harder enforcement of landing in the continuous state region corresponding to discrete state space element $\bar{s}_{t+1}$.

Additional details, such as providing warm start solution estimates based on historical planning actions and previous OCP invocations, as well as learning components of the dynamics, are left to a future work on implementation.  

At time $\tau_0=t=\hat{t}$, the first optimal control problem is solved:
\begin{subequations}
\label{eq:MPCinner} 
\begin{align}
\min_{z\in\mathbb{R}^{Hd_z},v\in\mathbb{R}^{Hd_v}}&\quad W(z_{(t+1)};p^w(a_t))  + \sum_{\tau(k)=(\tau(0)=t)}^{\tau(H)=t+1-h}\, G\left(z_{\tau(k)},v_{\tau(k)};p^g(a_t)\right)\,, \label{eq:cost0i}\\
\mathrm{s.t.} &\quad z_{\tau(k+1)} = f\left(z_{\tau(k)},v_{\tau(k)};p^f(a_t)\right),  \label{eq:dyn0i}\\
&\quad h\left(z_{\tau(k)},v_{\tau(k)};p^h(a_t)\right) \leq 0,\quad v_{\tau(k)} \in \tilde{A}(p^a(a_t))\,, \quad  \label{eq:const0i} \\ 
&\quad\quad \forall k\in \{0,\cdots,H-1\}\\ 
& z_{\tau(0)=t} = \hat{z}_t
\end{align}
\end{subequations}
After the optimization problem is solved to compute $\left\{z^*_{t,\tau(:)}(\hat{a}_t,\hat{z}_t),v^*_{t,\tau(:)}(\hat{a}_t,\hat{z}_t)\right\}$, the first control $\hat{v}(t)=v^*_{t,\tau(0)}$ is implemented by the robot controller. 

\begin{figure}[t]
\centering
\resizebox{\textwidth}{!}{%
\begin{tikzpicture}[
  axes/.style={thick,-{Stealth}},
  past/.style={thick, blue, smooth},
  closedloop/.style={thick, blue},
  pred1/.style={thick, dashed, red!40, smooth},
  pred2/.style={thick, dashed, red!70, smooth},
  pred3/.style={thick, dashed, red, smooth},
  block/.style={rectangle, draw, thick, fill=yellow!10, align=center, rounded corners, drop shadow, minimum height=2cm, text width=3cm}
]

  \draw[step=1cm, gray!20, very thin] (-0.5,-0.5) grid (10.5, 4.5);

  \draw[axes] (0,0) -- (11,0) node[right] {Time ($\tau$)};
  \draw[axes] (0,0) -- (0,4.5) node[above] {State ($z$)};
  
  \draw[thick] (2, 0.1) -- (2, -0.1) node[below] {$t$};
  \draw[thick] (4, 0.1) -- (4, -0.1) node[below] {$t\!+\!h$};
  \draw[thick] (6, 0.1) -- (6, -0.1) node[below] {$t\!+\!2h$};
  \node[below] at (8, -0.05) {$\dots$};
  
  \draw[past] (0,1.5) .. controls (1,2.0) .. (2,2.5) node[pos=0.4, above left, font=\footnotesize, text=black] {Past Trajectory};
  
  \draw[pred1] (2,2.5) .. controls (4, 3.5) and (5, 2.6) .. (6, 3.7);
  \node[font=\footnotesize, text=red!50, fill=white, inner sep=1pt] at (6.3, 3.8) {Pred. at $t$};
  
  \draw[closedloop] (2,2.5) .. controls (3, 2.8) .. (4, 2.9);
  \filldraw[red] (2,2.5) circle (3pt);
  \node[font=\footnotesize, fill=white, inner sep=1pt, below right=2pt] at (0.8, 3.1) {Apply $\hat{v}(t)$};
  
  \draw[pred2] (4,2.9) .. controls (6, 3.8) and (7, 3.0) .. (9, 3.6);
  \node[font=\footnotesize, text=red!80, fill=white, inner sep=1pt] at (8.5, 3.8) {Pred. at $t\!+\!h$};
  
  \draw[closedloop] (4,2.9) .. controls (5, 3.0) .. (6, 2.7);
  \filldraw[red] (4,2.9) circle (3pt);
  \node[font=\footnotesize, fill=white, inner sep=1pt, above left=2pt] at (4.6, 3.1) {Apply $\hat{v}(t\!+\!h)$};
  
  \draw[pred3] (6,2.7) .. controls (8, 2.3) and (9, 4.) .. (10.5, 3.1);
  \node[font=\footnotesize, text=red, fill=white, inner sep=1pt] at (10.3, 3.1) {Pred. at $t\!+\!2h$};
  
  \draw[closedloop] (6,2.7) .. controls (7, 2.5) .. (8, 2.8);
  \filldraw[red] (6,2.7) circle (3pt);
  \node[font=\footnotesize, fill=white, inner sep=1pt, below left=2pt] at (6.9, 2.5) {Apply $\hat{v}(t\!+\!2h)$};
  
  \node[font=\footnotesize, rotate=15, text=blue, fill=white, inner sep=1pt] at (2.8, 2.15) {\textbf{Closed-Loop Trajectory}};

  \node[block, right] (rl) at (11.5, 2.25) {
      \textbf{RL Critic} \\ \textit{Gradient Descent} \\ $\theta \leftarrow \theta - \alpha \nabla$ \\ \footnotesize Updates Cost $W, G$
  };
  
  \draw[->, thick, dashed, gray] (8, 2.8) .. controls (9.5, 2.5) and (10.5, 2.25) .. (rl.west);

\end{tikzpicture}%
}
\caption{Action-Triggered RLMPC Loop: The receding horizon principle. At each fast-timescale interval ($t, t+h, t+2h, \dots$), an optimal control problem is solved predicting the state trajectory over horizon $H$ (dashed lines). Only the first control action $\hat{v}(\tau)$ is applied, generating the true closed-loop trajectory (solid blue line). The RL Critic observes this realized performance and updates the cost parameters $\theta$.}
\label{fig:action_triggered_rlmpc}
\end{figure}

Subsequently, the system is evolved until time $t+h$, and the resulting state $\hat{z}(t-1+h)$ is measured. This is expected to be close to the computed estimated state from the optimization solution at $t$, that is $z^*_{t,\tau(1)}$, however, due to real-world state-system mismatch, they are not exactly equal. Still, their proximity enables fast Newton-based methods for computing warm-started solutions. Thus, another optimization problem is solved in real-time. The second problem, solved at and for time $t+h$, is now defined as:
\begin{align*}
\min_{z\in\mathbb{R}^{(H-1)d_z},v\in\mathbb{R}^{(H-1)d_v}}&\quad W(z_{(t+1)};p^w(a_t))  + \sum_{\tau=t+h}^{\tau(k)=(\tau(1)=t+1-h)}\, G\left(z_{\tau(k)},v_{\tau(k)};p^g(a_t)\right)\\
\mathrm{s.t.} &\quad z_{\tau(k+1)} = f\left(z_{\tau(k)},v_{\tau(k)};p^f(a_t)\right),  \\
&\quad h\left(z_{\tau(k)},v_{\tau(k)};p^h(a_t)\right) \leq 0,\quad v_{\tau(k)} \in \tilde{A}(p^a(a_t)),  \\ 
&\quad\quad \forall k\in \{1,\cdots,H-1\}\\ 
& z_t = \hat{z}(\tau(1))
\end{align*}
and again $\hat{v}(t+h)=v^*_{t+h,\tau(1)}$ is implemented and state $\hat{z}(\tau(2))=\hat{z}(t+2h)$ is observed.  This is repeated for $\tau(k)\in\{t-1+h,\cdots, t-2h,t-h\}$ in this rolling horizon to finally obtain the complete realization,
\[
\left\{\left\{\hat{z}(\tau),\hat{v}(\tau)\right\},\,\tau=t,\cdots,t+1-h,\right\}\cup \hat{z}(t+1)
\]
In addition, a real-time reward $\hat{R}^e(t+1)$ is observed externally at the conclusion of the problem, together with the evaluation of the action closed loop cost functional $\hat{R}^o(t+1)$. This quantity $\hat{R}^e(t+1)$ represents mission-relevant information that occurred in the course of the operation. For instance, if, in performing one of the control actions, the robot observed negative emotional affect in the human, this can enter into $\hat{R}^e(t+1)$ in a manner so as to inform the upper layer, and ultimately the task can be planned to be performed more agreeably to the human, that is, without action $\hat{a}_t$.
This action-triggered receding horizon procedure, augmented by the RL critic updates to the cost parameters, is summarized in Figure \ref{fig:action_triggered_rlmpc}.

\subsection{Integrating State representations}\label{sec:integrating_state_reps}

To connect the high-level planning state representation, which consists of logical atoms (Sec.~\ref{sec:planning}), with the low-level numeric state vector representation used in control (Sec.~\ref{sec:MPC}), we introduce a formalization based on the fuzzy membership functions (Sec.~\ref{sec:fuzzy}). This allows the high-level symbolic representation to be grounded in the continuous, low-level control space while preserving interpretability.

With the high-level state defined in terms of fuzzy logic values, we can attempt an integrated planning and control execution cycle. This consists of (i) selecting a high-level action $a_t$ by the planner based on current state $s_t$, (ii) executing the corresponding low-level optimal control process yielding $z_{t+1}$, and (iii) mapping back to a high-level symbolic description $s_{t+1}$ from the observed low-level system state $z_{t+1}$. The process then continues analogously with the next action $a_{t+1}$.

The high-level state $s$ is composed of a set of logical propositions (atoms) $\{L_1, L_2, \dots, L_M\}$ that describe \textit{discrete} properties of the system, such as ``robot A is at place B'' (\texttt{At(A,B)}) or ``object C is grasped by A'' (\texttt{grasped(C,A)}). These logical atoms must be somehow inferred from the underlying low-level state $z$, which consists of continuous physical quantities like position, velocity, or force. To achieve this, each possible logical atom $L_m$ is associated\footnote{possibly associated on the level of predicates from the domain language $\mathcal{L}$ defining $s_t$ (Sec.~\ref{sec:logic}), i.e. all atoms $L_m$ grounded from the same predicate $P(L_m)$ share the same membership function.} with a fuzzy membership function $\mu_{L_m}(z)$ that maps a low-level state to a truth value in $[0,1]$:
\begin{equation}
    \mu_{L_m}(z): \mathbb{R}^{d_z} \to [0,1].
\end{equation}
Here, $\mu_{L_m}(z) = 1$ indicates that the low-level state fully satisfies $L_k$, while $\mu_{L_m}(z) = 0$ indicates complete disagreement, with intermediate values capturing partial satisfaction.

For instance, if $L_1$ corresponds to the atom \texttt{At(RobotA,LocationB)}, its membership function could be defined as
\begin{equation}
\begin{array}{l}
    \mu_{At}(z) = \mathbf{1}\left(\|p_\text{RobotA} - p_\text{LocationB}\|^2\le \bar{d}^0_{At}\right)\\ \qquad +\mathbf{1}\left(\bar{d}^0_{At}\le \|p_\text{RobotA} - p_\text{LocationB}\|^2\le \bar{d}^1_{At}\right)
    \exp\left(-\frac{\|p_\text{RobotA} - p_\text{LocationB}\|^2-\bar{d}^0_{At}}{\sigma^2}\right),
    \end{array}
\end{equation}
where $p_\text{RobotA}$ and $p_\text{LocationB}$ are their respective positions, and $\sigma$ controls the smoothness of the transition. Note the truncation of the distribution: if the robot is within distance $\bar{d}^0_{At}$ of \texttt{LocationB}, we declare that it is \texttt{At} that location with certainty, and if it is a distance at least $\bar{d}^1_{At}$ away, then he is crisply declared to not be at the location, with a fuzzy transition function in the annulus defined by the second indicator in the expression. This translation mechanism is visualized in Figure~\ref{fig:fuzzification_layer}.

\begin{figure}[t]
\centering
\begin{adjustbox}{width=\textwidth}
\begin{tikzpicture}[
  layer/.style={rectangle, draw, thick, minimum width=3.5cm, minimum height=3cm, align=center, fill=white, drop shadow},
  arr/.style={-{Stealth[scale=1.5]}, thick, line width=1.5pt}
]
  \node[layer, fill=gray!5] (continuous) at (0,0) {
      \textbf{Continuous Space} \\ ($z_t \in \mathbb{R}^{2}$) \\ \vspace{0.2cm}
      \tikz{
          \draw[step=0.5cm, gray!50, very thin] (0,0) grid (2,1);
          \fill[red] (1.2,0.6) circle (0.1);
      }
  };
  \node[layer, fill=purple!5] (fuzzy) at (5,0) {
      \textbf{Fuzzification} \\ ($\mu_{L_m}$) \\ \vspace{0.2cm}
      \tikz[xscale=0.5, yscale=0.5]{
          \draw[thick, blue] (-2,0) .. controls (-1,1) .. (0,0);
          \draw[thick, red] (-1,0) .. controls (0,1) .. (1,0);
          \draw[->] (-2.2,0) -- (1.2,0);
      }
  };
  \node[layer, fill=green!5] (discrete) at (10,0) {
      \textbf{Logic Atoms} \\ ($s_t^\mu \in [0,1]^M$) \\ \vspace{0.2cm}
      \footnotesize
      \begin{tabular}{|l|l|}
           \hline
           At(A) & 0.9 \\
           At(B) & 0.1 \\
           \hline
      \end{tabular}
  };
  
  \draw[arr] (continuous) -- node[midway, above, font=\small] {Map $z \to \mu$} (fuzzy);
  \draw[arr] (fuzzy) -- node[midway, above, font=\small] {Agg. $T^s$} (discrete);
\end{tikzpicture}
\end{adjustbox}
\caption{(De-)Fuzzification Layer: The translation from continuous physical coordinates (left), through fuzzy membership functions (middle), to discrete logic atoms (right) for the planner.}
\label{fig:fuzzification_layer}
\end{figure}

The high-level state $s$ can then be represented as an operation on a vector of fuzzy truth values:
\begin{equation}
    s^{\mu} = T^s(\mu_{L_1}(z), \mu_{L_2}(z), \dots, \mu_{L_m}(z);a_t,s_t).
\end{equation}
ensuring a continuous mapping between low-level numeric and high-level symbolic descriptions, allowing logical planning (reasoning) to be performed in a way that adapts to the underlying control dynamics.



\subsection{Updating the Planning State and Model}\label{sec:updating_planning_state}  
After execution of the low-level control, the system reaches some observed state $\hat{z}(t+1)$. The high-level fuzzy state then corresponds to some new fuzzy truth values:
\begin{equation}
    \hat{s}^{\mu}_{t+1} = T^s(\mu_{L_1}(\hat{z}(t+1)), \dots, \mu_{L_M}(\hat{z}(t+1));a_t,s_t).
\end{equation}
where $\hat{s}^{\mu}$ is now a fuzzy realization of the state observed and $T^s:[0,1]^{M}\times \mathcal{A}\times\mathcal{S}\to [0,1]$

The actual state $\hat{s}_{t+1}$ can now be defined as an appropriate defuzzification operation $\mathcal{D}^{\iota_{t+1}}:[0,1]\to \{0,1\}$ where $\iota_{t+1}\in \mathcal{I}_{\mathcal{D}}$ is one of a set of possible operations of size $\left\vert \mathcal{I}_{\mathcal{D}}\right\vert =I_D$.

A principled, albeit computationally extensive, approach to identify the actual high-level state would be to select the state that best matches the observation:
\begin{align*}
    \hat{s}_{t+1} = & \mathcal{D}^{\iota_{t+1}}\left(T^s(\mu_{L_1}(\hat{z}(t+1)), \dots, \mu_{L_M}(\hat{z}(t+1));a_t,s_t)\right)\\ & =  \arg \max_{s \in \mathcal{S}} \bigwedge_{m} \mu_{L_m}(\hat{z}(t+1)).
\end{align*}

Recall that at the initiation of the process, there was an expected state $\bar{s}_{t+1}$. This was defined to be computed as the model estimate:
\[
\bar{s}_{t+1}=\arg\max_{s\in\mathcal{S}} p_{\Psi}^{S}(t)=\arg\max_{s\in\mathcal{S}}p_{\Psi}^{S}(\theta^{\psi};a_t,s_t)  
\]
We can now use the observed state $\hat{s}_{t+1}$ to update the model $\theta^{\psi}$, e.g., by a gradient update with stepsize $\gamma_t$,
\begin{equation}\label{eq:updatepsi}
\begin{array}{l}
\theta^{\psi}_{t+1} = (1-\gamma^{\psi}_t)\theta^{\psi}_t+\gamma^{\psi}_t\frac{\partial \mathbb{P}_{p_{\Psi}^S(t+1)}}{\partial\theta^{\psi}}(\hat{s}_{t+1})-\frac{\gamma^{\psi}_t}{S-1}\frac{\partial \mathbb{P}_{p_{\Psi}^S(t+1)}}{\partial\theta^{\psi}}(\bar{s}_{t+1})
\\\qquad\qquad -\frac{\gamma^{\psi}_t}{(S-1)^2}\sum\limits_{s\in\mathcal{S}\setminus \hat{s}_{t+1}}\frac{\partial \mathbb{P}_{p_{\Psi}^S(t+1)}}{\partial\theta^{\psi}}(s)
\end{array}
\end{equation}


\paragraph{Updating the Upper Level Model}

Consider that a reward $\hat{R}_{t+1}$ is observed upon the realization $(\hat{a}_t,\{\hat{z}_{\tau},\hat{v}_{\tau}\},\hat{s}_{t+1})$. With a robotic aide, this could be defined by either prompting the human for direct feedback, observing the human's expression and estimating emotional affect, as well as externalities on the surrounding environment judged and remarked upon by observers, e.g., in a human-feedback pre-training operation of the system.

\[
\begin{array}{l}
\hat{R}(t+1) = R(\hat{R}^e(t+1),\hat{R}^o(t+1)) \\
\hat{R}^o(t+1)=W(\hat{z}(t+1);p^w(a_t))+\sum\limits_{\tau(k)=(\tau(0)=t)}^{\tau(H)=t+1-h} G(\hat{z}(\tau(k)),\hat{v}(\tau(k));p^g(a_t))
\end{array}
\]

If the computed state $\hat{s}_{t+1}\notin \mathcal{G}$, then we use the intermediate reward to learn the intermediate cost function estimate, say by:
\begin{equation}\label{eq:updateavf}
\theta^C_{t+1}=(1-\gamma^C_t)\Phi^C_{t+1}(\theta^C_t,s_t,a_t) + \frac{\partial \Phi^C(\theta^C_t;a_t,s_t)}{\partial \theta^C}(\hat{R}_{t+1}-\bar{C}_{t+1})
\end{equation}
and if the state is in the goal set $\hat{s}_{t+1}\in\mathcal{G}$, then a similar update is made to $\bar{C}_T$. However, in addition, the loop ceases and the robot is idle until the next task arrives. 

The action-value function can be updated in standard fashion:
\begin{equation}\label{eq:updateavf}
\theta^V_{t+1}=(1-\gamma^V_t)\theta^V_{t+1}V(s_t,a_t) + \frac{\partial \Phi^V(\theta^V_t;a_t,s_t)}{\partial \theta^V}(\hat{R}_{t+1}-\bar{V}_{t+1})
\end{equation}
where now $\bar{V}_{t+1}$ is the action value function estimated by $\bar{V}_{t+1}=\Phi^V(\theta_t^V;\hat{a}_t,\hat{s}_t)$. Observe that the difference ensures the right sign is used to update the parameters relative to the gradient.


\subsection{Updating the MPC Control Parametrization}\label{sec:updating_mpc}
In addition to, or instead of, updating the model of the state transition that operates at the planning layer, that is $p^S_{\psi}(\theta^{\psi};s_{t},a_t)$ from the predicted actual state mismatch ($\bar{s}_{t+1}$ verus $\hat{s}_{t+1}$) and the intermediate cost $\Phi^C(\theta^C;s_a,a_t)$ through observed rewards $\hat{C}_{t}$, we can also use the mismatch to update the action to MPC map $a_t\to (p^w,p^g,p^f,p^h,p^a)$. That is, if the expected and desired final state is not reached upon triggering the MPC associated with action $a_t$, and/or there is unexpected reward or cost, then this information can be used to adjust the specific MPC problem generated by $a_t$.

We can compute the continuous form of the final state mismatch as:
\[
\hat{E}^{\psi}_t := T^{\hat{s}_{t+1}}\bigwedge_{m} \mu_{L_m}(\hat{z}(t+1)) -T^{\bar{s}_{t+1}}\bigwedge_{m} \mu_{L_m}(\hat{z}(t+1)) 
\]
that is, the fuzzy relaxation of $s_{t+1}$ as far as activation of the real state $\hat{s}_{t+1}$ compared to the originally predicted state $\bar{s}_{t+1}$. Through differentiation, we seek to obtain a correspondence between the discrete state transition through the fuzzy relaxation of the final state map given the final continuous state $\hat{z}(t+1)$ to the optimal MPC action that generated $\hat{z}_{t+1}$, through the parameters defining the MPC action as set by $a_t$.

There is also the reward mismatch between the actual realized reward through the low level state operation and the predicted cost associated with this action,
\[
\hat{E}^{C}_t  := \hat{R}_{t+1}-\bar{C}_{t+1}
\]
Now, applying the chain rule, we target the terminal cost in the MPC that we denoted by $W(z_{(t+1)};p^w(a_t))$ for adjustment as far as final state mismatch and the intermediate cost $G(z_{\tau(k)},v_{\tau(k)};p^g(a_t))$ for the reward mismatch, with learning rate $\alpha^M_t$ that may depend on the stage iteration $t$:
\[
\begin{array}{l}
p^w(a_t) \gets p^w(a_t)-\alpha^{M,w}_t \hat{E}^{\psi}_t \frac{d\hat{E}^{\psi}}{dp^w(a_t)} \\
\frac{d\hat{E}^{\psi}}{dp^w(a_t)} = \frac{d\hat{E}^{\psi}}{d\hat{z}_{t+1}}\frac{d\hat{z}_{t+1}}{dp^w(a_t)}\\
p^g(a_t) \gets p^g(a_t)-\alpha^{M,g}_t \hat{E}^{C}_t \sum\limits_{\tau=t+h}^{t+1-h}\left[\frac{dG}{d(z,v)}\frac{d(z_{\tau},v_{\tau})}{dp^g(a_t)}+\frac{\partial G}{\partial p^g(a_t)}\right] 
\end{array}
\]
\paragraph{MPC Sensitivities}\label{sec:mpc_sensitivities}


We now describe how to compute the sensitivity of the fuzzy function $\mu_{L_m}(z^\star_{t+1})$ with respect to the MPC parameters $p_t$, where $z^\star_{t+1}$ is the end state of a simulation of the system in closed-loop with the MPC controller. These sensitivities will be useful to adjust the MPC parameters at the planning layer to ensure that the MPC scheme achieves the goals it is assigned.

We assume that the MPC scheme is run in a possibly complex but differentiable simulation model
\begin{equation}\label{eq:sim}
z^\star_{t+h(i+1)} = f^\star( z^\star_{t+hi}, u_{t+hi}),\quad u_{t+hi} = \pi^\mathrm{MPC}( z^\star_{t+hi},p_t)
\end{equation}
where $\pi^\mathrm{MPC}$ is the MPC policy, and $p_t$ are the MPC parameters selected by the supervisory layer, i.e.
\begin{equation}
    \pi^\mathrm{MPC}( z^\star_{t+hi},p_t) = u_0^\star
\end{equation}
where $u_0^\star$ is the first element of the optimal input sequence $u^\star$ defined by solving the Parametric Nonlinear Program
\begin{subequations}
\label{eq:MPC}
\begin{align}
    \min_{x,u}&\quad T(x_N,p_t) + \sum_{k=0}^{N-1}L(x_k,u_k,p_t) \label{eq:MPC:obj} \\
    \mathrm{s.t.}&\quad x_{k+1} =  f^\star(x_k,u_k,p_t),\qquad x_0= z_{t+hi}^\star \label{eq:MPC:eq}\\
    &\quad h(x_k,u_k,p_t)\geq 0 \label{eq:MPC:ineq} 
\end{align}
\end{subequations}
parametrized with the current simulation state $z^\star_{t+hi}$ and parameters $p_t$. We then establish how to differentiate a closed-loop simulation $\hat z^\star_{t,\ldots,t+h\tau}$ with respect to the MPC parameters $p_t$. Two approaches are useful, depending on the end-goal of the simulation differentiation: the forward and adjoint modes. 
\begin{itemize}
\item \textbf{Forward mode}
This mode is useful if we are interested in the ``overall" sensitivity of the simulation. In that mode, we make a forward pass on the simulation as: 
\begin{align}
   & \frac{\mathrm d z^\star_{t+h(i+1)}}{\mathrm d p_t} = \frac{\partial f^\star }{\partial  z^\star_{t+hi}}\frac{\mathrm dz^\star_{t+hi}}{\mathrm d p_t} + \frac{\partial f}{\partial u_{t+hi} }\frac{\partial \pi^\mathrm{MPC}}{\partial p_t},\qquad  \frac{\mathrm d z^\star_{t}}{\mathrm d p_t} = 0
\end{align}
where all expressions are evaluated at $z^\star_{t+hi}$, $u_{t+hi} = \pi^\mathrm{MPC}(z^\star_{t+hi},p_t)$ and $p_t$ alongside the forward simulations \eqref{eq:sim}.
\item \textbf{Adjoint mode}
If we are only interested in the sensitivity of the scalar quantity $\mu_{L_m}( z^\star_{t+h\tau})$ at the end of the simulation, then the adjoint-mode sensitivities is typically more effective. In that mode, the sensitivity is computed backward in time after the forward simulations \eqref{eq:sim} have been performed. We then define the adjoint variable:
\begin{align}
    \lambda^i = \frac{\mathrm d \mu_{L_m}(z^\star_{t+h\tau})}{\mathrm d z^\star_{t+hi} }^\top,
\end{align}
whose dynamics read as
\begin{align}
    \lambda^{i-1} = \left.\frac{\partial f}{\partial  z^\star_{t+(i-1)h}}^\top\right|_{z^\star_{t+hi},p_t,u_{t+hi}}  \cdot \lambda^{i},\qquad \lambda^\tau =\frac{\partial \mu_{L_m}(z^\star_{t+h\tau})}{\partial z^\star_{t+h\tau}}^\top.
\end{align}
Finally, the sensitivity of $\mu_{L_m}( z^\star_{t+h\tau}) = \mu_{L_m}( z^\star_{t+1})$ to $p_t$ reads as
\begin{align}
\frac{\mathrm d \mu_{L_m}( z^\star_{t+1})}{\mathrm d p_t }^\top = \sum_{i=0}^\tau \left.\left(\frac{\partial f}{\partial u_{t+hi}}\frac{\partial \pi^\mathrm{MPC}}{\partial p_t}\right)^\top\right |_{z^\star_{t+hi},p_t,u_{t+hi}}\cdot\lambda^i,
\end{align}
where $u_{t+hi} = \pi^\mathrm{MPC}(z^\star_{t+hi},p_t)$.
\end{itemize}
Both the forward and adjoint-mode methods require computing the MPC policy sensitivity $\frac{\partial \pi^\mathrm{MPC}}{\partial p_t}$ at every simulation step $i$.
It can be computed as follows.

Consider $w$ collecting the state-input variables $x,u$ of MPC \eqref{eq:MPC}, and\\ $w^\star(z^\star_{t+hi}, p_t)$ the optimal solution to \eqref{eq:MPC} for a give initial condition $z^\star_{t+hi}$ and MPC parameters $p_t$. Consider in addition function $G(w,z^\star_{t+hi}, p_t)$ collecting the equality constraints \eqref{eq:MPC:eq} and function $H(w,z^\star_{t+hi}, p_t)$ collecting the inequality constraints \eqref{eq:MPC:ineq}, and the Lagrange function:
\begin{align}
\mathcal L(\zeta,z^\star_{t+hi}, p_t) = &T(x_N,p_t) + \sum_{k=0}^{N-1}L(x_k,u_k,p_t) + \lambda^\top G(w,z^\star_{t+hi}, p_t) \\&+ \mu^\top H(w,z^\star_{t+hi}, p_t)\nonumber
\end{align}
where $\zeta =\left\{w,\lambda,\mu\right\}$ collects the primal-dual variables associated to \eqref{eq:MPC}, and $\lambda,\mu$ are the dual variables associated to the equality and inequality constraints \eqref{eq:MPC:eq}-\eqref{eq:MPC:ineq}. The primal-dual Interior-Point KKT conditions underlying a solution to \eqref{eq:MPC} read as:
\begin{subequations}
\begin{eqnarray}
\kappa_\nu(\zeta,z^\star_{t+hi}, p_t) = \left[\begin{array}{c}
\nabla_w \mathcal L(\zeta,z^\star_{t+hi}, p_t)\\
G(w,z^\star_{t+hi}, p_t)\\
\mu \circ H(w,z^\star_{t+hi}, p_t) + \nu 
\end{array}\right] = 0
\end{eqnarray}
\end{subequations}
where $\circ$ is the element-wise multiplication, and where $\nu \geq 0$ is the barrier parameter underlying primal-dual Interior Point methods. The NLP is typically solved with $\nu$ small.   
Assuming that the objective \eqref{eq:MPC:obj} is at least twice continuously differentiable, and that the constraints \eqref{eq:MPC:eq}-\eqref{eq:MPC:ineq} are at least continuously differentiable, the Implicit Function Theorem (IFT) guarantees the existence of the solution sensitivities~$\frac{\partial \zeta^\star}{\partial p_t}$ at $p_t$ if the linear independence constraint qualification, second-order sufficient conditions, and strict complementarity are satisfied at the solution~$\zeta^\star$. Under these conditions, the IFT furthermore implies that the solution sensitivity can be computed by solving the linear system:
\begin{align}\label{eq:IFT}
\left.\frac{\partial \kappa_\nu(\zeta,z^\star_{t+hi}, p_t)}{\partial \zeta}\frac{\mathrm d \zeta^\star}{\mathrm d p_t} + \frac{\partial \kappa_\nu(\zeta,z^\star_{t+hi}, p_t)}{\partial p_t}\right|_{\zeta = \zeta^\star} = 0
\end{align}
for $\frac{\mathrm d \zeta^\star}{\mathrm d p_t}$. The sensitivity of the MPC prediction can then be extracted from that solution. Note that the linear system \eqref{eq:IFT} is well-posed if the MPC scheme achieves the LICQ and SOSC conditions.

Recent software developments implement algorithms to both solve the MPC problem and to compute the sensitivity of its solution efficiently \cite{reiter2025synthesis}. The algorithmic parameter $\tau$ can be used to provide a solution manifold to the NLP that is smooth everywhere, such that the sensitivities are well defined for all $p_t$. Note that the approximation introduced by $\nu > 0$ is of the order $\nu$.

\subsection{Fuzzy control relaxation alternative}
Alternatively, instead of explicitly identifying a discrete high-level state $s_{t+1}$ after observing $\hat{z}(t+1)$, one may employ a fuzzy control approach (Sec.~\ref{sec:fuzzy}), where the action $a_{t+1}$ selection is performed directly in the low-level continuous space. In this case, the high-level planner can be replaced by a fuzzy controller $\tilde{\zeta}$ that directly determines the next action based on the computed fuzzy memberships:

\begin{equation}
    \hat{a}_{t+1} =\hat{a}^{\zeta}_{t+1} \tilde{\zeta}(\mu_{L_1}(\hat{z}(t+1)), \dots, \mu_{L_M}(\hat{z}(t+1)).
\end{equation}

This formulation removes the need for explicit high-level symbolic state inference by treating the fuzzy memberships as the effective state representation for decision-making. The mapping $\tilde{\zeta}$ may be implemented akin to a Mamdani fuzzy rule-based system (Sec.~\ref{sec:fuzzy}), for example:

\begin{equation}
\begin{array}{l}
\text{IF } \mu_{L_1} \text{ is high AND } \mu_{L_2} \text{ is low, THEN } a_{t+1} = a_1,\\
\text{IF } \mu_{L_2} \text{ is high, THEN } a_{t+1} = a_2.
\end{array}
\end{equation}

To relax even further, instead of selecting an intermediate high-level action $a_t$, the fuzzy inference system could even directly output the low-level control $v(\tau)$. This would further eliminate any need for the high-level discrete description, as well as the MPC procedure
\begin{equation}
    v(\tau) = \tilde{\zeta}(\mu_{L_1}(\hat{z}(\tau)), \dots, \mu_{L_m}(\hat{z}(\tau))).
\end{equation}
at the expense of interpretability and optimality.

\paragraph{Hybrid Approach}

Ultimately, the policy can be defined as some hybrid approach with respect to the computed policy:
\[
(\hat{a}_t,v_t)\sim \Delta\left(\pi;\hat{a}^g_t,\hat{a}^{\epsilon}_t,\hat{a}^{\sigma}_t,\hat{a}^c_t,\tilde{\zeta}_t\right) 
\]
where $\Delta$ is the unit simplex over the arguments. That is, the Algorithm can compute a classical planning algorithm solution, an RL policy, or a fuzzy control relaxation scheme with some probability across the options.

\subsection{Scheduling-Planning Integration}

We propose that the apex problem is one of Scheduling. Namely, within each daytime period of 16 hours, the robot has a series of tasks that it can perform for the human it is of service to: $\mathcal{T}=\{T_b,T_l,T_d,T_{m_1},T_{m_2},T_w,T_e\}$ corresponding to feeding breakfast, lunch and dinner (i.e., $f(T_b)< s(T_l)$), providing medicine regime 1 and 2, and an event triggered task of assistance to the use the toilet and provision of some entertainment (music, jokes, etc.). We consider that the corresponding events follow Poisson distributions with parameters $\lambda_b$ and $\lambda_e$.

We can consider that the goal is to ensure all tasks are completed within the day, with some regularization towards spreading the tasks evenly throughout the day. We can see that this intuitively yields a solution of, for instance, 30 minutes duration for each task corresponding to arranging to perform $\{T_b,T_{m_1},T_l,T_{m_2},T_d\}$ sequentially with about, e.g., $f(T_b)=s(T_{m_1})-2.5$ hours. However, in addition, we add another objective term $\mathbf{F}$ defining the overall quality of the experience of humans throughout their interaction with the agent as a function of the tasks scheduled. We consider that this function is initially unknown, but is a function of the tasks that are scheduled and when, as well as an objective $F_T$ corresponding to the positive experience of that given task as it happened in real time. Day by day, $\mathbf{F}$ and its components $\{F_T\}$ for $T\in\mathcal{T}$, are learned from interaction with the human, and the Scheduler is optimized for peak human service experience.

If the toilet event is triggered, $T_w$ is performed and any other task ceases and is completed following $f(T_w)$. If the entertainment event is triggered, the task is performed simultaneously with the meals, while it is delayed until after $f(T_{m_1})$ and $f(T_{m_2})$. These rules can all be encoded into logical constraints. Note that if an event is triggered in a manner that interrupts a scheduled task, then the schedule must be recomputed. 

We can consider that there is an event that is triggered by, e.g., erratic behavior, irregular toilet behavior, change in complexion, etc. that notifies an external human health worker. 

We can consider that each task spawns a planning problem. That is, for duration $d(T)$, a planning problem $(\mathcal{P}_T,\mathcal{S}_T,\mathcal{A}_T,\Phi_T)$ (with $\Phi_T=(\psi(T),\phi_G(T))$ being the dynamics) is executed. The goal to be satisfied in the planning problem is defined by completing the task that it corresponds to, which can be defined to be $n\le K(t)$, where $K(t)$ is the maximum number of states that can be traversed in the time duration $d(t)$ defined for task $t$ and $s_n\models \phi_G$.

As an example, for $T_l$, actions could include \texttt{moving to the kitchen}, which must happen first, \texttt{slicing vegetables}, \texttt{boiling stew}, \texttt{sprinkling a teaspoon of dill}, etc. Some of these, in sequence, construct a meal. The previous Section defined the operation of the planning problem that sequences the set of actions to perform the triggered task. These are also associated with a lower level control system. Note that while the physical dynamics of the robot, that is, $f$ above, exhibit significant invariance across Tasks, and thus planning problems, the constraints $h$ as well as the terminal and intermediate costs $W$ and $G$ are expected to be defined in a manner influenced by the valence of the goal(s) associated to the Task.

There is an additional observed reward $\hat{R}(T_l)$ that quantifies how pleasant the experience was for the human. This can be measured by, e.g. a 1-10 rating given by the human after each task and/or reading emotional affect during the operation of the planning task. In addition a live or delayed video recording observed by a health expert can grade the performance in optimality to patient physical and mental health \footnote{We expect that this human-in-the-loop pretraining will be important for adequate alignment to the operation's needs and goals.} This reward is meant to learn an additional Value Function for the choice of Task, but must be learned through interaction with the environment. Recall that for each task/planning problem, there are multiple sequences of actions that yield the goal state, and the goal is to find the one that maximizes the objective $F^{T_l}(s,a;T_l,\mathcal{E}(s_0))$. The incorporation of an additional parametrization by vector $\mathcal{E}$ represents the time of day the task is performed as well as other possibly relevant covariates and and a vector describing the read emotional affect of the human. Then the Task Scheduler can learn the Value of a particular planning problem. Note that it depends on the initial state, and is otherwise determined by the algorithm defining vectors $\{s_{t},a_t\}_{t\in[n]}$. 



\section{Experimental Setup}\label{sec:experiments}
To validate the proposed multilevel architecture, we instantiate the framework in a health service domain. A robot must learn individual patient preferences across two distinct task types - medication delivery and meal preparation - while navigating a realistic physical environment with constrained energy capacity and obstacles to circumnavigate. This setting exercises the full hierarchy: physics-informed MPC at the control level, fuzzy state estimation bridging continuous and discrete representations, multi-objective task planning at the symbolic level, and dual-loop preference learning connecting the levels through differentiable interfaces.
The full implementation is available at \url{https://github.com/alexdifre/hospital-robot-system}.

\subsection{System Architecture}\label{sec:exp_architecture}
\begin{figure}[!htbp]
    \centering
    \includegraphics[width=\linewidth]{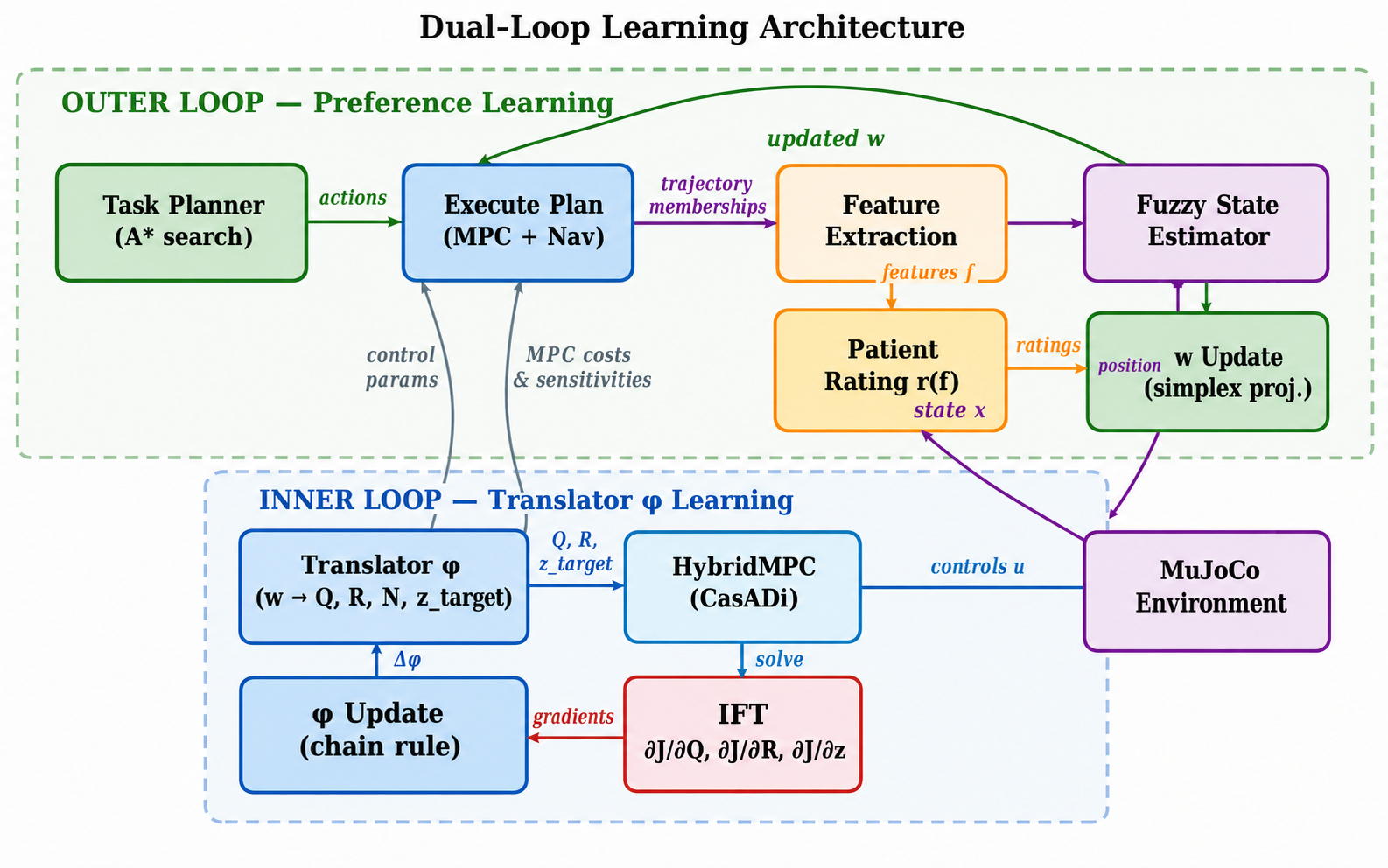}
    \caption{Case study architecture. The \textbf{outer loop} learns patient preference weights $\mathbf{w}$ on the simplex via gradient descent on patient ratings, updating the task planner's multi-objective cost function. The \textbf{inner loop} adapts action-dependent MPC target parameters via IFT sensitivities, enabling end-to-end differentiation through the optimization-based controller.}
    \label{fig:architecture}
\end{figure}
Figure~\ref{fig:architecture} shows the general systems architecture for the symbolic and continuous control loops, coupled through learning.

\paragraph{Preference learning}
The upper (mission-oriented) level seeks to perform medication-delivery and meal-preparation tasks formulated as planning problems that minimize mismatch with the cared-for human's preferences. A task planner searches over the discrete action space using the scalarized multi-objective cost function $c(\mathbf{x}, \mathbf{u}) = \mathbf{w}^\top \mathbf{f}(\mathbf{x}, \mathbf{u})$, where \(\mathbf{w} \in \Delta^4 = \{\mathbf{w} \in \mathbb{R}^5_{\geq 0} : \sum_{i=1}^5 w_i = 1\}\) is a weight vector over the five objectives, and $\mathbf{f}$ extracts normalised cost features along five dimensions: time efficiency, safety, energy consumption, proximity to patient, and approach quality. After each episode, the patient provides dimension-wise ratings $r_i$, which define the observed feedback signals, and the preference weights are updated via projected gradient descent on the simplex to minimise the discrepancy between predicted and observed ratings.
\paragraph{Continuous control loop.}
The lower level corresponding to the faster time scale operates the physical actuator decisions of the robot through MPC.

Specifically,the MPC stage cost is a quadratic tracking objective:
\begin{equation}
    \ell(z_{\tau}, u_{\tau}) = (z_{\tau} - z_{\text{target},\tau}(\hat{\mathbf{w}}))^\top 
    Q(\hat{\mathbf{w}})(z_{\tau} - z_{\text{target},\tau}(\hat{\mathbf{w}})) + 
u_{\tau}^\top R(\hat{\mathbf{w}}) u_{\tau}
    \label{eq:mpc_stage_cost}
\end{equation}
where $Q(\hat{\mathbf{w}}) = \text{diag}(q_1(\hat{\mathbf{w}}), \ldots, 
q_6(\hat{\mathbf{w}}))$ is a diagonal state cost matrix whose entries are 
affine functions of the preference weights, $R(\hat{\mathbf{w}}) = 
\text{diag}(R_1(\hat{\mathbf{w}}), R_2(\hat{\mathbf{w}}), 
R_3(\hat{\mathbf{w}}))$ is a learned diagonal positive-definite control 
cost matrix. The prediction horizon is set to $N = 20$.


\subsection{Health Facility Description}\label{sec:environment}
\begin{figure}[!htbp]
    \centering
    \includegraphics[width=0.85\linewidth]{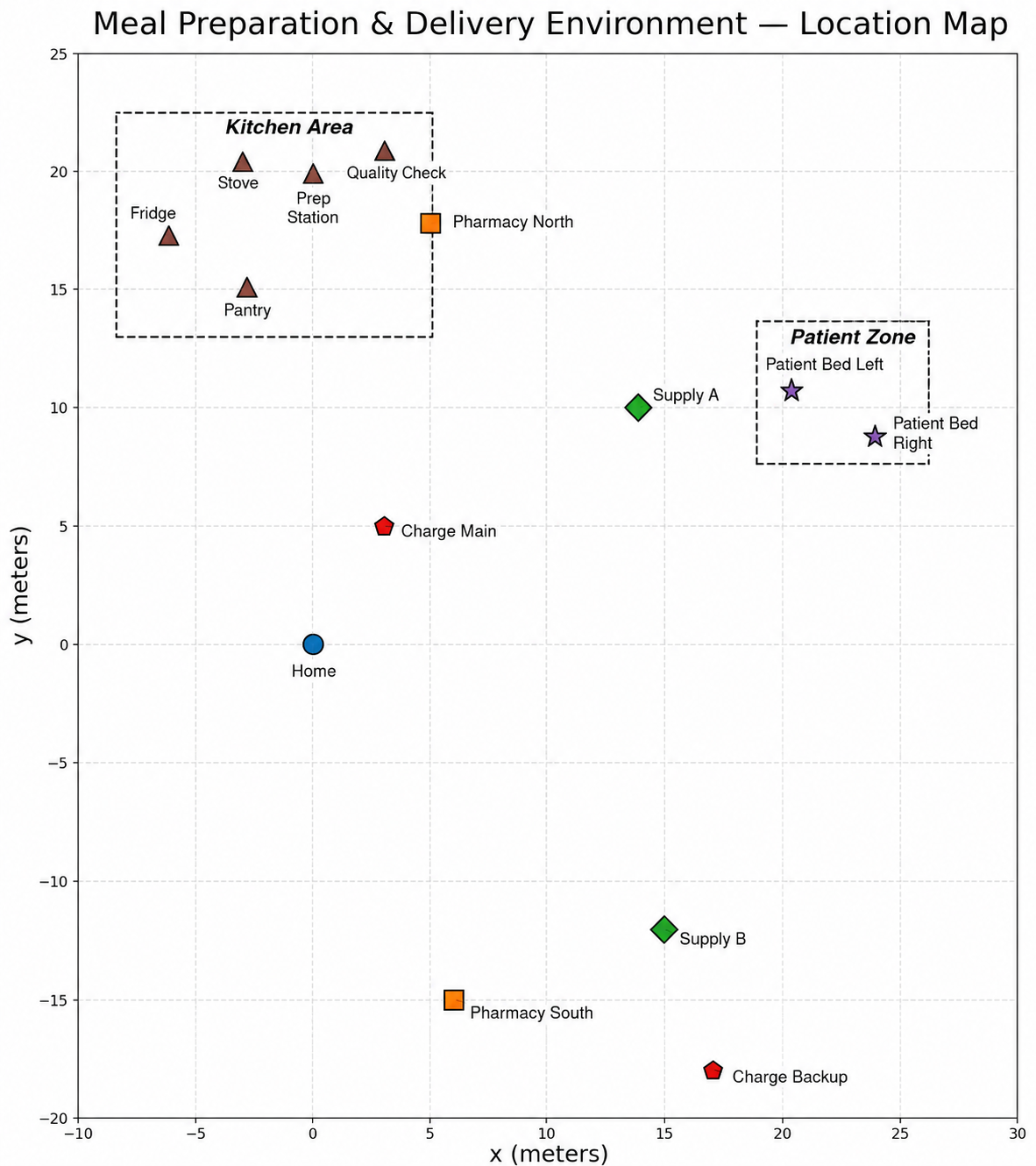}
    \caption{ The layout of the environment showing all navigable locations, colour-coded by functional
type}
    \label{fig:floor_plan}
\end{figure}

The environment (Figure~\ref{fig:floor_plan}) models a care center as a 2D continuous space with named locations spanning six functional categories: pharmacies (north and south), supply rooms (A and B),
a kitchen cluster (pantry, fridge, preparation station, stove, and quality-check station),
patient bedside positions (left and right approach), charging stations (main and backup),
and the home location.

Each transition between two locations is associated with a congestion risk value
\(\rho(i,j) \in (0,1)\), representing the probability of encountering obstacles
or personnel while moving from location \(i\) to location \(j\). The risk is defined with
lower values assigned to safer transitions and higher values assigned to routes
crossing more congested areas, such as the main kitchen workspace or paths near
the stove.

These risk values feed directly into the safety component of the multi-objective cost function, creating a concrete tradeoff between route efficiency and safety.

The robot is modeled as a differential-drive platform with continuous position, orientation, and velocity states. Navigation between symbolic locations is handled by a hybrid planning-control stack: the symbolic planner selects the next task-level target, and the MPC controller tracks a continuous reference toward that target while respecting the robot dynamics and obstacle-avoidance constraints. This realizes the action-triggered MPC scheme: each planning action $a_t$ triggers a rolling-horizon optimal control problem whose target steers the robot toward the continuous-state region corresponding to the predicted next discrete state $\bar{s}_{t+1}$. The differential-drive dynamics serve as the system model $f$, with obstacle avoidance encoded via the inequality constraints $h(z_{\tau}, u_{\tau}) \leq 0$.
\subsection{Task Specifications}\label{sec:tasks}
The simulation incorporates two PDDL planning tasks: medication delivery, encoded by \texttt{domain\_med.pddl} and \texttt{problem\_med.pddl}, and meal preparation, encoded by \texttt{domain\_meal.pddl} and \texttt{problem\_meal.pddl}. Both tasks use the same five-dimensional cost-feature vector $\mathbf{f} \in [0,1]^5$, whose components correspond  to \texttt{nav-time}, \texttt{location-safety}, \texttt{nav-battery}, \texttt{proximity-feature}, and \texttt{approach-feature}. These components are scalarised by the weights \texttt{w0}--\texttt{w4}, where \texttt{w0} is time, \texttt{w1} is safety, \texttt{w2} is battery, \texttt{w3} is proximity, and \texttt{w4} is approach preferences. Formally, $\bar{C}(s_t, a_t) = \mathbf{w}^\top \mathbf{f}(s_t, a_t)$ serves as the learnable stage cost; in the PDDL implementation this quantity is accumulated in \texttt{stage\_cost}, together with \texttt{total\_time} and \texttt{discrete\_steps} in the final metric.
The total cost of the planner is defined as follows:
\begin{align*}
J(\pi)
&= \sum_{t=0}^{n-1} \overline{C}(s_t,a_t)
+ \overline{C}_T(s_n)
+ \psi(n)
\end{align*}

where:
\begin{align*}
\overline{C}(s_t, a_t)
&= \Phi^C(W^C; s_t, a_t) = \mathbf{w}^{\top} \mathbf{f}(s_t, a_t), \\[4pt]
\overline{C}_T(s_n)
&= \Phi^T(W^T; s_n)
= w_0\, T_{\text{total}}
\end{align*}

and:
\begin{itemize}
    \item \(T_{\text{total}}\) is the total execution time.
    \item \(w_0\) is the weight associated with the time cost.
    \item \(\psi(n)\) is a penalty term associated with the number of discrete actions in the plan.
\end{itemize}


\subsubsection{Medication Delivery}
\begin{figure}[!htbp]
    \centering
    \includegraphics[width=\linewidth]{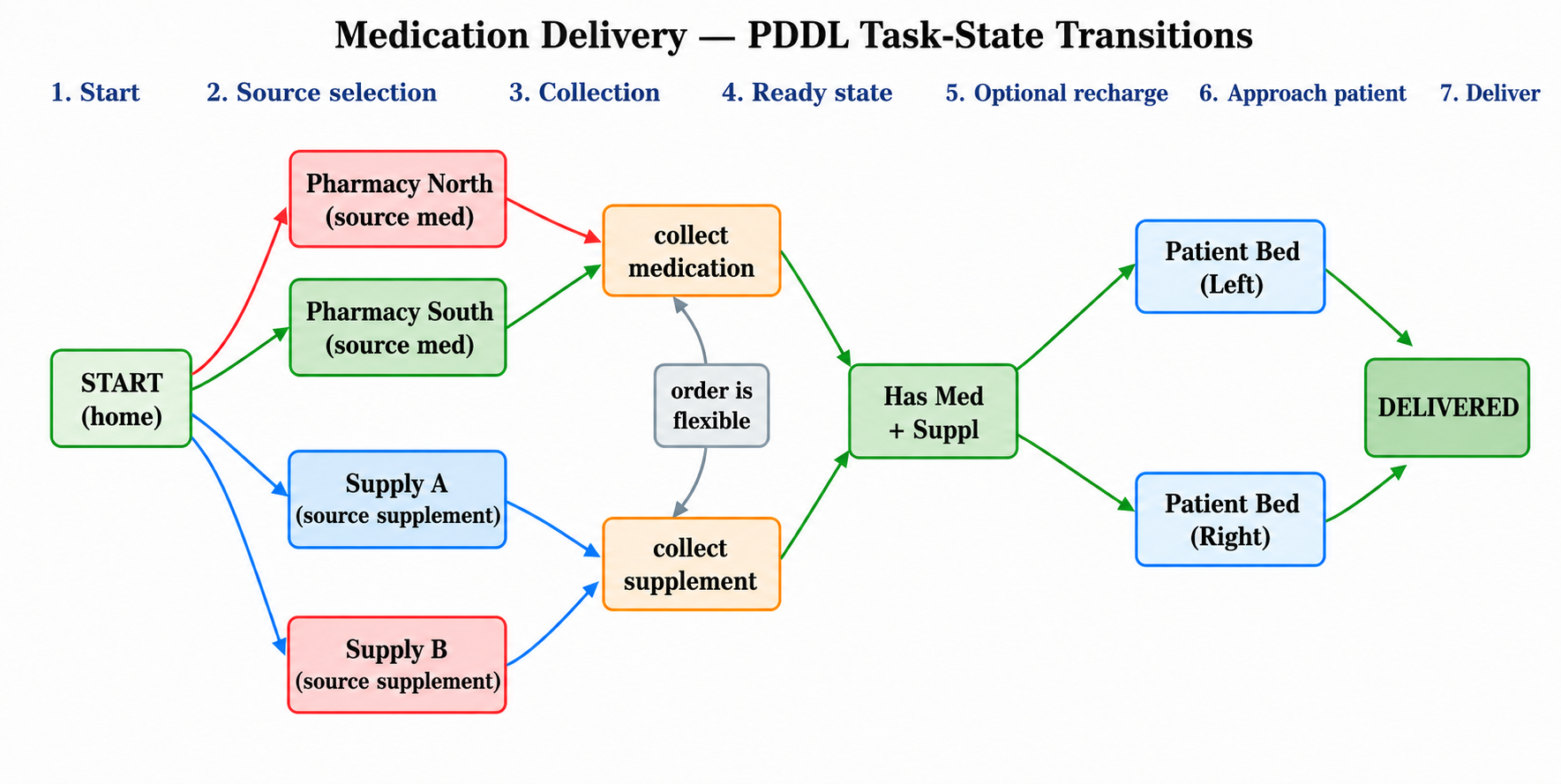}
    \caption{Medication delivery task state transitions. The robot must collect medication from a pharmacy and a supplement from a supply room then deliver both to the patient.}
    \label{fig:med_task}
\end{figure}
The medication delivery task (Figure~\ref{fig:med_task}) is encoded as a PDDL problem specified by the predicates
\pddl{requested-medicine} and \pddl{requested-supplement}. These predicates
define the requested items that the robot must collect before delivery.

The planner can choose among three groups of actions: navigation actions,
collection/checking actions, and delivery actions.

\begin{table}[h]
\centering
\small
\begin{tabularx}{\linewidth}{lX}
\hline
\textbf{Action group} & \textbf{Available actions} \\
\hline
Pharmacy navigation &
\pddl{go\_to\_pharmacy\_north}, \pddl{go\_to\_pharmacy\_south} \\

Supply navigation &
\pddl{go\_to\_supply\_a}, \pddl{go\_to\_supply\_b} \\

Charging navigation &
\pddl{go\_to\_charge\_main}, \pddl{go\_to\_charge\_backup} \\

Bedside approach &
\pddl{approach\_left\_to\_bed}, \pddl{approach\_right\_to\_bed} \\
\hline
\end{tabularx}
\caption{Main navigation choices in the medication-delivery domain.}
\end{table}

After reaching a pharmacy location, the robot can collect the requested
medicine and verify whether it is correct. The intended branch is:

\[
\pddl{collect\_med\_antibiotic}
\;\rightarrow\;
\pddl{check\_medication\_correct}.
\]

The alternative actions \pddl{collect\_med\_analgesic} and
\pddl{collect\_med\_insulin} model wrong-medication branches. If the collected
medicine is incorrect, the planner can execute:

\[
\pddl{check\_medication\_wrong}
\;\rightarrow\;
\pddl{put\_down\_medicine}.
\]

Similarly, after reaching a supply location, the intended supplement branch is:

\[
\pddl{collect\_supp\_vitamin\_d}
\;\rightarrow\;
\pddl{check\_supplement\_correct}.
\]

The alternative actions \pddl{collect\_supp\_electrolyte} and
\pddl{collect\_supp\_omega3} represent wrong-supplement choices. These branches
are handled through:

\[
\pddl{check\_supplement\_wrong}
\;\rightarrow\;
\pddl{put\_down\_supplement}.
\]

Delivery is enabled only when all required conditions hold:

\[
\pddl{has-medication}
\wedge
\pddl{has-supplement}
\wedge
\pddl{correct-medication}
\wedge
\pddl{correct-supplement}
\]

Once these predicates are true, the robot can deliver the items using either:

\[
\pddl{deliver\_on\_bedside\_table\_left}
\quad \text{or} \quad
\pddl{deliver\_on\_bedside\_table\_right}.
\]

The goal in \pddl{problem\_med.pddl} is:

\[
\pddl{(delivered robot1)}.
\]

Plan diversity arises from \emph{route selection} and \emph{branch selection}.
Route selection concerns the choice between alternative pharmacy routes, supply
locations, charging stations, and bedside approaches. Branch selection concerns
whether the robot collects the correct item immediately or enters a correction
branch after detecting a wrong medicine or supplement.

Battery feasibility is encoded directly in the movement-action preconditions.
Each movement action decreases \pddl{battery-soc} by \pddl{nav-battery}, and the
action is allowed only if the remaining charge is sufficient to reach the nearest
charger, represented by \pddl{nearest-charger-battery}. This structure prevents
the planner from selecting movements that would leave the robot unable to
recover. In addition, conditional penalties at $0.15$ and $0.25$ state of charge
make low-battery behavior visible to the weighted cost function.
\subsubsection{Meal Preparation}
\begin{figure}[!htbp]
    \centering
    \includegraphics[width=\linewidth]{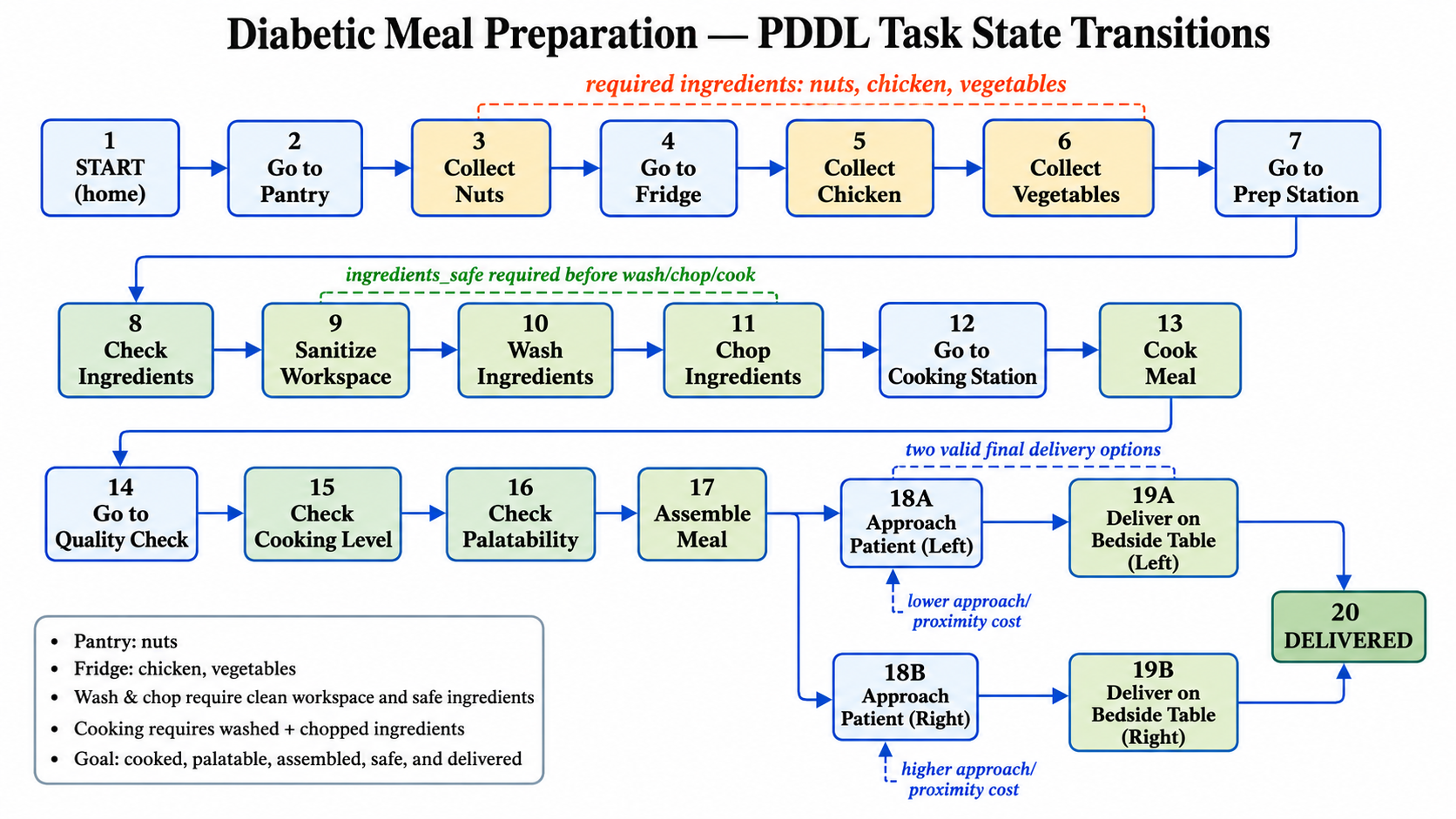}
    \caption{Meal preparation task state transitions. The PDDL problem allows the robot to prepare one of three meals: \texttt{sandwich}, \texttt{salad}, and \texttt{hot\_meal}.}
    \label{fig:meal_task}
\end{figure}

The meal preparation task (Figure~\ref{fig:meal_task}) is encoded as the PDDL problem
\texttt{meal\_preparation\_problem} over the domain
\texttt{meal\_preparation}. The robot \texttt{robot1} starts at
\texttt{home}; the task instance does not fix a single meal, but allows the planner to choose among
\texttt{sandwich}, \texttt{salad}, and \texttt{hot\_meal} through the
actions \texttt{choose\_sandwich}, \texttt{choose\_salad}, and
\texttt{choose\_hot\_meal}, which establish the predicate
\texttt{meal\_to\_prepare}. The required ingredients depend on the
selected meal: \texttt{sandwich} requires \texttt{bread} and
\texttt{vegetables}, \texttt{salad} requires \texttt{nuts} and
\texttt{vegetables}, while \texttt{hot\_meal} requires
\texttt{chicken} and \texttt{vegetables} and is the only meal satisfying \texttt{requires\_cooking}. Ingredient availability is
represented by the predicate \texttt{ingredient\_available}:
\texttt{bread}, \texttt{nuts}, and \texttt{vegetables} are available
at \texttt{pantry}, whereas \texttt{chicken} is available at
\texttt{fridge}. The directed navigation graph is defined through
\texttt{connected} relations linking \texttt{home},
\texttt{pantry}, \texttt{fridge}, \texttt{prep\_station},
\texttt{stove}, \texttt{quality\_check}, and the patient delivery
locations, with recharge actions available through
\texttt{charge\_main} and \texttt{charge\_backup}.

The main symbolic sequence is summarized as follows:

\begin{itemize}
    \item \textbf{Meal selection and ingredient acquisition:}
    \texttt{choose\_sandwich}, \texttt{choose\_salad}, or
    \texttt{choose\_hot\_meal} establishes
    \texttt{meal\_to\_prepare}. Navigation actions such as
    \texttt{go\_to\_pantry} and \texttt{go\_to\_fridge} move the robot
    to the appropriate ingredient locations, while
    \texttt{collect\_ingredient} establishes
    \texttt{collected\_ingredient} and removes the corresponding
    \texttt{missing\_ingredient} facts.

    \item \textbf{Ingredient and workspace safety:}
    \texttt{check\_ingredients} establishes
    \texttt{ingredients\_checked} and
    \texttt{ingredients\_safe}. The action
    \texttt{sanitize\_workspace} establishes
    \texttt{workspace-clean} and
    \texttt{robot-hands-clean}, while removing
    \texttt{cross-contamination-risk}.

    \item \textbf{Meal preparation:}
    \texttt{wash\_ingredients},
    \texttt{chop\_ingredients}, and, for meals satisfying
    \texttt{requires\_cooking}, the action
    \texttt{cook\_meal} are executed, followed by
    \texttt{check\_cooking\_level},
    \texttt{check\_palatability}, and
    \texttt{assemble}. These actions establish the predicates
    \texttt{ingredients\_washed},
    \texttt{ingredients\_chopped},
    \texttt{meal\_cooked} (when applicable),
    \texttt{cooking\_level\_checked},
    \texttt{meal\_palatable}, and
    \texttt{meal\_assembled}.

    \item \textbf{Delivery:}
    \texttt{approach\_to\_left\_side} or
    \texttt{approach\_to\_right\_side} establishes
    \texttt{can\_be\_deliverable}. The meal is then delivered through
    \texttt{deliver\_on\_bedside\_table\_left} or
    \texttt{deliver\_on\_bedside\_table\_right}, establishing the
    predicate \texttt{delivered}.
\end{itemize}

The goal specified in \texttt{problem\_meal.pddl} is the conjunction of
the predicates \texttt{meal\_palatable},
\texttt{meal\_assembled}, and \texttt{delivered}.


\subsection{Patient Profiles}\label{sec:profiles}
\begin{figure}[!htbp]
    \centering
    \includegraphics[width=0.98\linewidth]{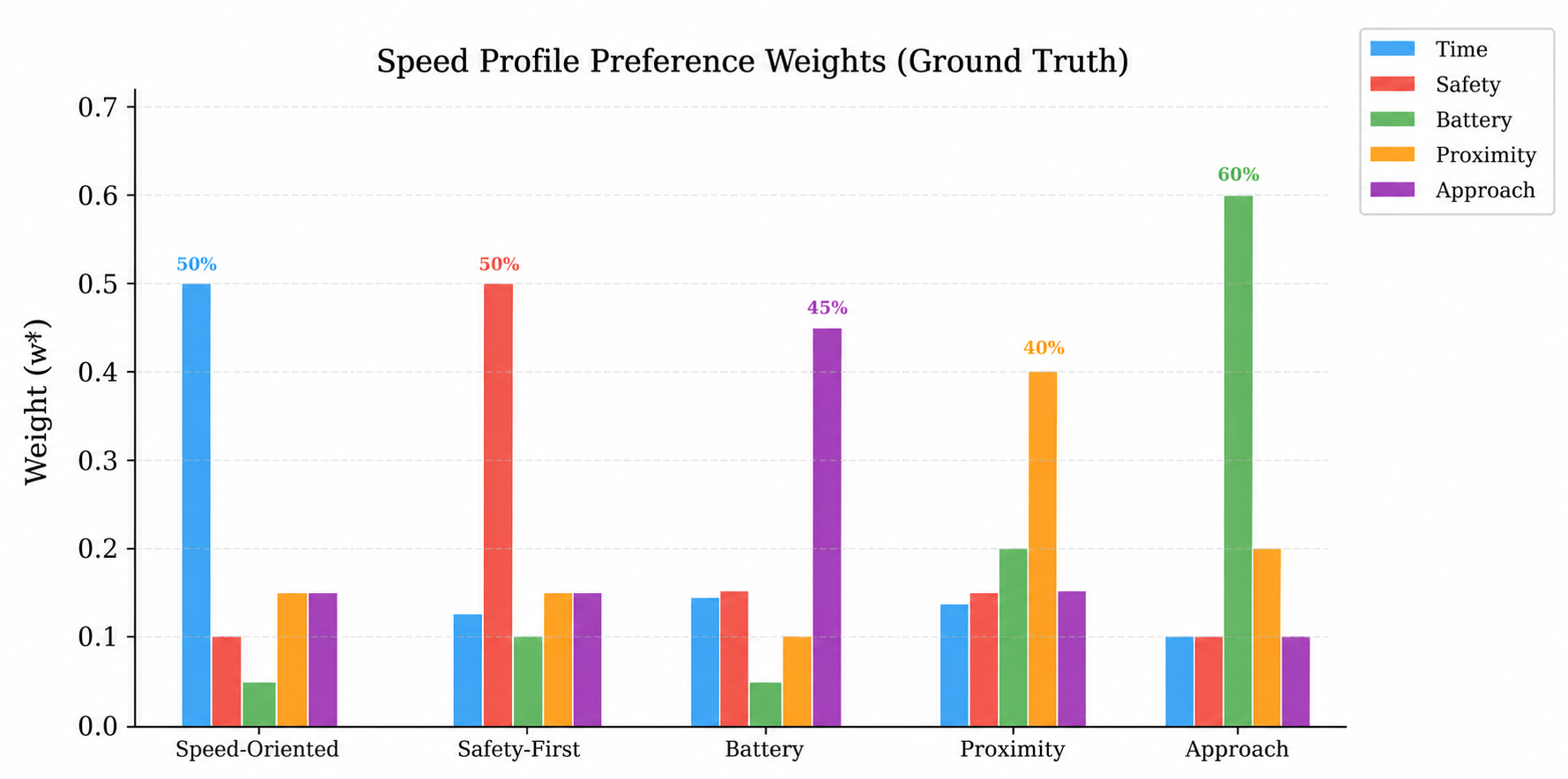}
    \caption{Ground-truth preference weight vectors $\mathbf{w}^*$ for five patient profiles. Each profile has a distinct dominant dimension reflecting the patient's primary concern: time efficiency (speed-oriented), safety from congestion risk (safety-first), battery conservation, proximity to bedside, or meal/delivery presentation quality.}
    \label{fig:profiles}
\end{figure}
We define five patient profiles (Figure~\ref{fig:profiles}), each represented by a ground-truth weight vector $\mathbf{w}^* \in \Delta^4$ on the four-dimensional simplex over five objectives. The profiles are designed to exercise different regions of the preference space; Table~\ref{tab:profiles} lists the corresponding weight vectors.

\begin{table}[ht]
\centering
\color{black} 
\small
\begin{tabular}{lccccc}
\toprule
\textbf{Profile} & \textbf{Time} & \textbf{Safety} & \textbf{Battery} & \textbf{Proximity} & \textbf{Approach} \\
\midrule
Speed-Oriented       & \textbf{0.50} & 0.12 & 0.14 & 0.14 & 0.10 \\
Safety-First         & 0.10 & \textbf{0.50} & 0.15 & 0.15 & 0.10 \\
Energy-Conscious     & 0.15 & 0.15 & \textbf{0.45} & 0.15 & 0.10 \\
Comfort-Focused      & 0.15 & 0.15 & 0.10 & \textbf{0.40} & 0.20 \\
Presentation-Focused & 0.05 & 0.10 & 0.05 & 0.20 & \textbf{0.60} \\
\bottomrule
\end{tabular}
\caption{Patient profile weight vectors. Bold values indicate the dominant dimension. }
\label{tab:profiles}
\end{table}

Online feedback corresponding to human ratings from the patients is simulated as $r_i = 5 - 4 \cdot f_i \cdot w^*_i + \epsilon_i$, where $f_i$ is the normalised cost-feature value and $\epsilon_i \sim \mathcal{N}(0, 0.1)$ is a realization of observation noise. Lower feature values (faster delivery, less risk exposure, less battery usage) yield higher ratings, with penalties weighted by the patient's true preferences. The learning mechanism updates the estimate of the preference scalarization vector from the sequence of these feedback signals through a projected gradient update applied on the simplex.
\subsection{Fuzzy State Estimation}\label{sec:exp_fuzzy}
\begin{figure}[!ht]
    \centering
    \includegraphics[width=0.75\linewidth]{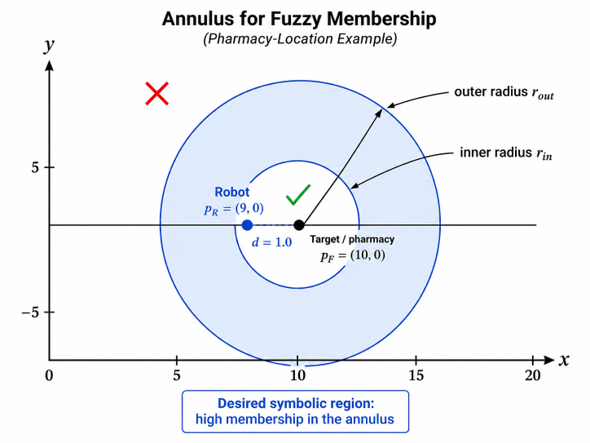}
    \caption{The symbolic predicate associated with a location is represented by a fuzzy radial region centered at the target location, where membership depends on the robot--target distance}
    \label{fig:fuzzy}
\end{figure}
\paragraph{Location membership.}
For each named location $\ell$ with centre position $\mathbf{p}_\ell$, let $d_\ell(\mathbf{p}) = \|\mathbf{p} - \mathbf{p}_\ell\|$. The robot's membership at continuous position $\mathbf{p}$ is computed as a saturated radial distance ramp:
\begin{equation}
    \mu_\ell(\mathbf{p}) =
    \begin{cases}
        1, & d_\ell(\mathbf{p}) \leq 1,\\
        \frac{5 - d_\ell(\mathbf{p})}{4}, & 1 < d_\ell(\mathbf{p}) < 5,\\
        0, & d_\ell(\mathbf{p}) \geq 5.
    \end{cases}
\end{equation}
Thus, a robot closer than one metre to the location centre has full membership, while a robot farther than five metres has zero membership. Task preconditions are evaluated against a threshold $\mu_\ell(\mathbf{p}) \geq 0.7$, determining when the continuous state instantiates the PDDL atom \texttt{(at robot1 }$\ell$\texttt{)} and the robot is considered sufficiently at that location to perform in-place actions.

\subsubsection{Expected Location and Fuzzy Mismatch Detection}

The symbolic planner defines the location that the robot is expected to reach after executing the current action. This location is not obtained from the fuzzy memberships. It comes directly from the planned task. For example, if the current action is a movement toward the pharmacy, the expected location is the pharmacy location. If the action is a delivery to a patient, the expected location is the corresponding patient-bed location.

Let
\[
\ell_{\mathrm{exp}} \in \mathcal{L}
\]
be the expected symbolic location given by the planner, where \(\mathcal{L}\) is the set of all symbolic locations in the environment.

The fuzzy estimator, instead, uses the continuous state of the robot to estimate where the robot actually is. For each symbolic location \(\ell \in \mathcal{L}\), a membership function assigns a value
\[
\mu_{\ell}(x) \in [0,1],
\]
where \(x\) is the current continuous state of the robot.

The fuzzy estimator then selects the symbolic location with the highest membership value, provided that this membership exceeds the same symbolic-instantiation threshold:
\[
\ell_{\mathrm{fuzzy}}
=
\arg\max_{\ell \in \mathcal{L}} \mu_{\ell}(x),
\quad
\max_{\ell \in \mathcal{L}} \mu_{\ell}(x) \geq 0.7.
\]

A mismatch is detected when the location expected by the planner is different from the location predicted by the fuzzy estimator:
\[
\mathrm{Mismatch}
=
\left(
\ell_{\mathrm{fuzzy}} \neq \ell_{\mathrm{exp}}
\right).
\]

In other words, the planner expects the robot to reach \(\ell_{\mathrm{exp}}\), while the fuzzy estimator predicts that the robot is actually in \(\ell_{\mathrm{fuzzy}}\). If the two locations are the same, the execution is considered consistent with the symbolic plan. 

When a mismatch is detected, the mismatch is logged together with the final geometric error and the fuzzy membership values. In the current implementation, the symbolic task state is then aligned with the planner's expected location before the next replanning step, while the logged geometric mismatch is used for continuous target adaptation.

\subsection{Stopping Criteria for MPC Execution}

The execution process has two stopping levels. The first one concerns a single MPC optimization, while the second one concerns the closed-loop execution of a complete task leg.

At the lower level, each MPC call solves one finite-horizon optimal control problem using Acados. In the current implementation, Acados is configured with a full SQP nonlinear solver and tolerance
\[
\varepsilon_{\mathrm{tol}} = 10^{-4}.
\]
A single MPC solve is considered successful when Acados returns status 0. In that case, the first control input of the optimized trajectory is applied to the system. If Acados fails, the controller reports a solver failure.

At the higher level, each task leg is executed in closed loop over a direct waypoint reference generated for the current symbolic action. At every iteration, the current robot state is observed, the MPC problem is solved for the active waypoint reference, the first control input is applied, and the environment state is updated. The execution advances to the next waypoint when the robot is within the waypoint tolerance, or when the per-waypoint iteration allocation is exhausted. The leg stops when all waypoint targets have been processed, when the global closed-loop step budget is exhausted, or when an Acados solver failure occurs.

\newpage

\section{Simulation Results}\label{sec:results}
This section presents the experimental results for the multilevel architecture described in Section~\ref{sec:experiments}.

\subsection{Convergence Dynamics}\label{subsubsec:plan_diversity}

\begin{figure}[H]
    \centering
    \includegraphics[width=0.85\textwidth]{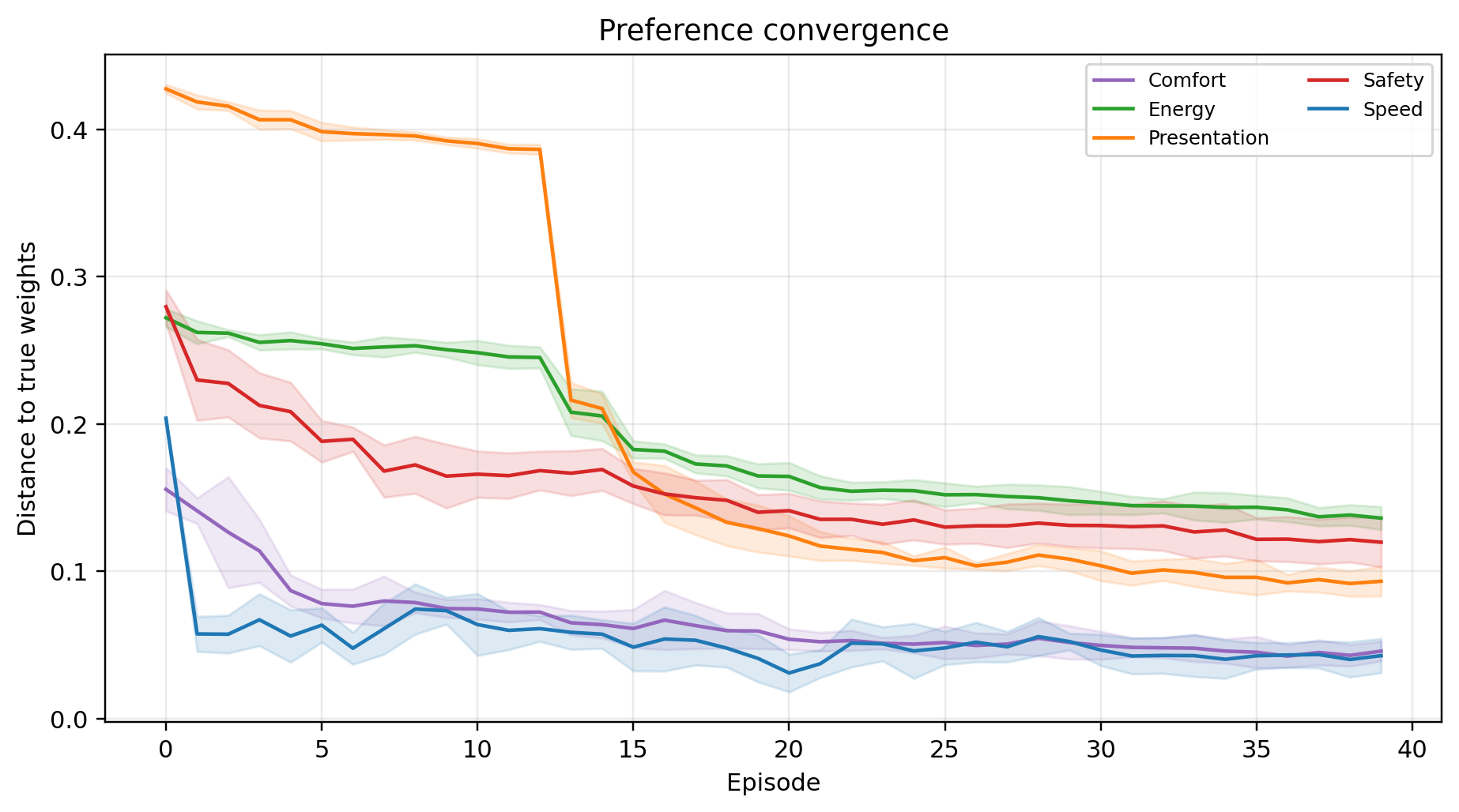}
    \caption{Distance to true preference vector $d = \|\hat{w} - w^{*}\|_{2}$ over 40 episodes for all five patient profiles. Solid lines show the median across 5 seeds; shaded bands denote the interquartile range.}
    \label{fig:weight-distance}
\end{figure}

All experiments use 40 episodes per run with alternating medication delivery and meal preparation tasks, 5 random seeds per condition, and the hospital environment described in Section~\ref{sec:environment}. Convergence is assessed by the $L_2$ distance $d = \|\hat{\mathbf{w}} - \mathbf{w}^*\|_2$ between the learned and true preference vectors, with a convergence threshold of $d \leq 0.15$ for all profiles. Hyperparameters are held constant across all profiles and seeds and are set to: preference learning rate $\eta = 0.12$, initial exploration noise $\sigma_0 = 0.15$ with decay factor $0.2$, and translator learning rate $0.002$.

The system was evaluated across all five patient profiles with 5 random seeds each, producing 25 independent runs of 40 episodes. Of these, \textbf{25 out of 25 runs achieved convergence}, with all profiles reaching $d \leq 0.15$ across all seeds. A convergence episode of $0$ indicates that the uniformly initialised vector already satisfies the distance threshold for that profile.


\begin{table}[htbp]
    \centering
    \small
    \begin{tabular}{lccccccc}
        \toprule
        \textbf{Profile} & \textbf{$d_{\mathrm{thresh}}$} &
        \textbf{Conv. Ep.} & \textbf{Best $d$} & \textbf{Final $d$} & \textbf{Conv.} & \textbf{Plans} \\
        \midrule
        Speed    & 0.15 & 1 (1--1)    & 0.023 $\pm$ 0.008 & 0.043
                 & 5/5 & 3 \\
        Safety   & 0.15 & 9 (9--17)   & 0.118 $\pm$ 0.018 & 0.120
                 & 5/5 & 3 \\
        Present. & 0.15 & 16 (16--17) & 0.084 $\pm$ 0.006 & 0.093
                 & 5/5 & 2 \\
        Comfort  & 0.15 & 0 (0--1)    & 0.037 $\pm$ 0.006 & 0.046
                 & 5/5 & 1 \\
        Energy   & 0.15 & 29 (21--30) & 0.136 $\pm$ 0.008 & 0.136
                 & 5/5 & 3 \\
        \midrule
        \textbf{Overall} & & & & & \textbf{100\%} &  \\
        \bottomrule
    \end{tabular}
    \caption{Full-system convergence results (5 profiles $\times$ 5 seeds,
    40 episodes each). Conv.\ Ep.\ = median episodes to first reach
    $d_{\mathrm{thresh}}$ (IQR in parentheses). Task Rate = fraction of
    episodes with successful delivery. Plans =  structural plan
    variants observed across episodes.}
    \label{tab:results}
\end{table}

\newpage

\subsection{Weight Learning Accuracy}

\begin{figure}[H]
    \centering
    \includegraphics[width=0.7\textwidth]{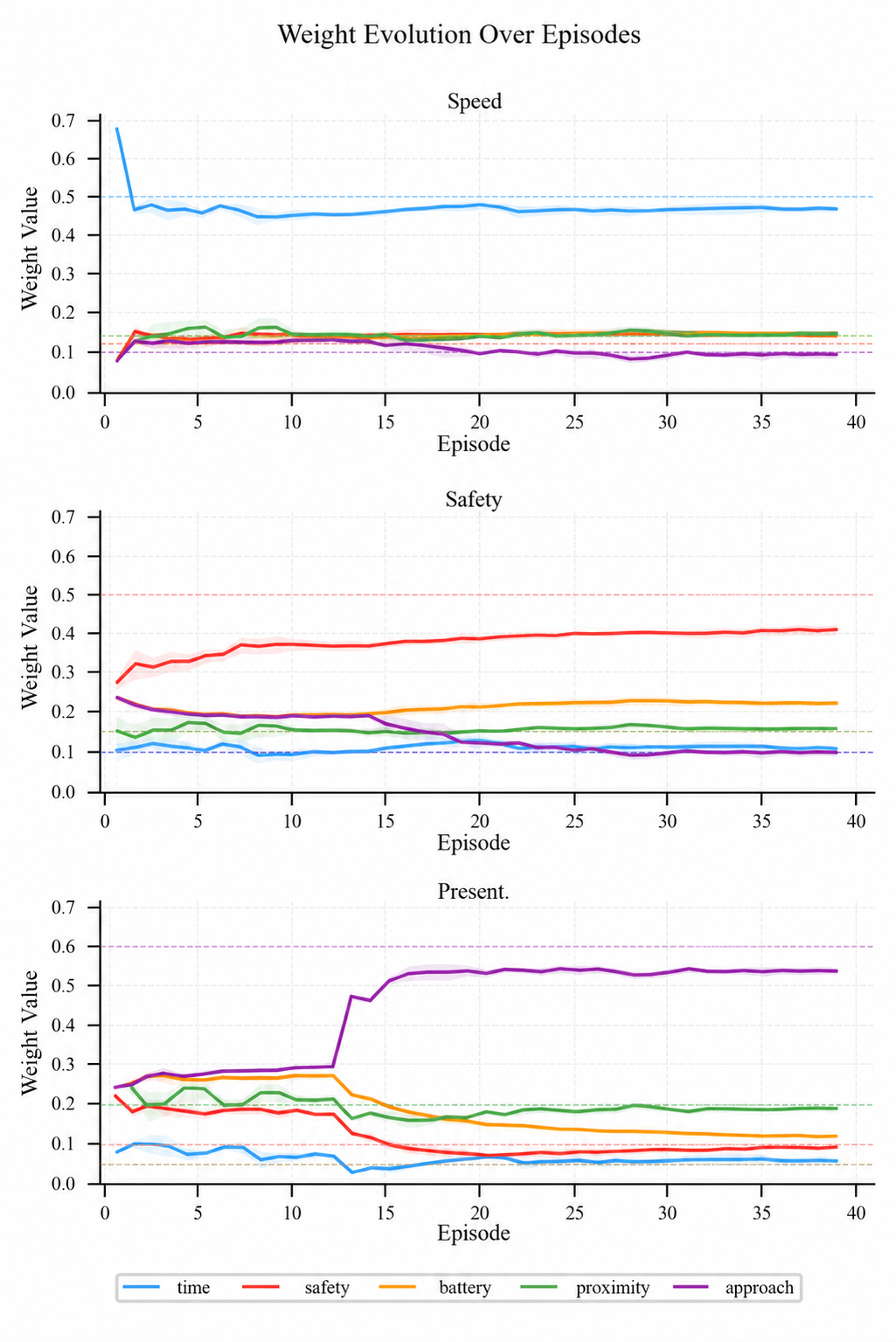}
    \caption{Weight evolution over episodes for the speed-oriented and presentation-focused profiles. Dashed lines indicate the corresponding true weights.}
    \label{fig:weight-trajectories}
\end{figure}

Figure~\ref{fig:weight-trajectories} shows the weight trajectories for a fast-converging speed-oriented profile and the presentation-focused profile. For the speed-oriented setting, the time weight moves rapidly from the uniform initialisation of $0.20$ toward the target value of $0.50$, reaching an average value of about $0.47$ by the final episode. The profile converges consistently across all seeds, with a median convergence episode of $1$ and a final mean distance of $0.043$.

For the presentation-focused profile, the approach weight increases more gradually from about $0.24$ in the first episode to about $0.54$ by episode $40$, remaining slightly below the true target value of $0.60$. The current JSON results show that this profile does converge under the standard setting: all $5$ seeds converge, with convergence episodes between $16$ and $19$ and a median of $16$. This indicates that the approach signal is harder to learn than the speed-oriented preference, but it is still recovered reliably enough to satisfy the convergence threshold.

\begin{figure}[H]
    \centering
    \includegraphics[width=0.95\textwidth]{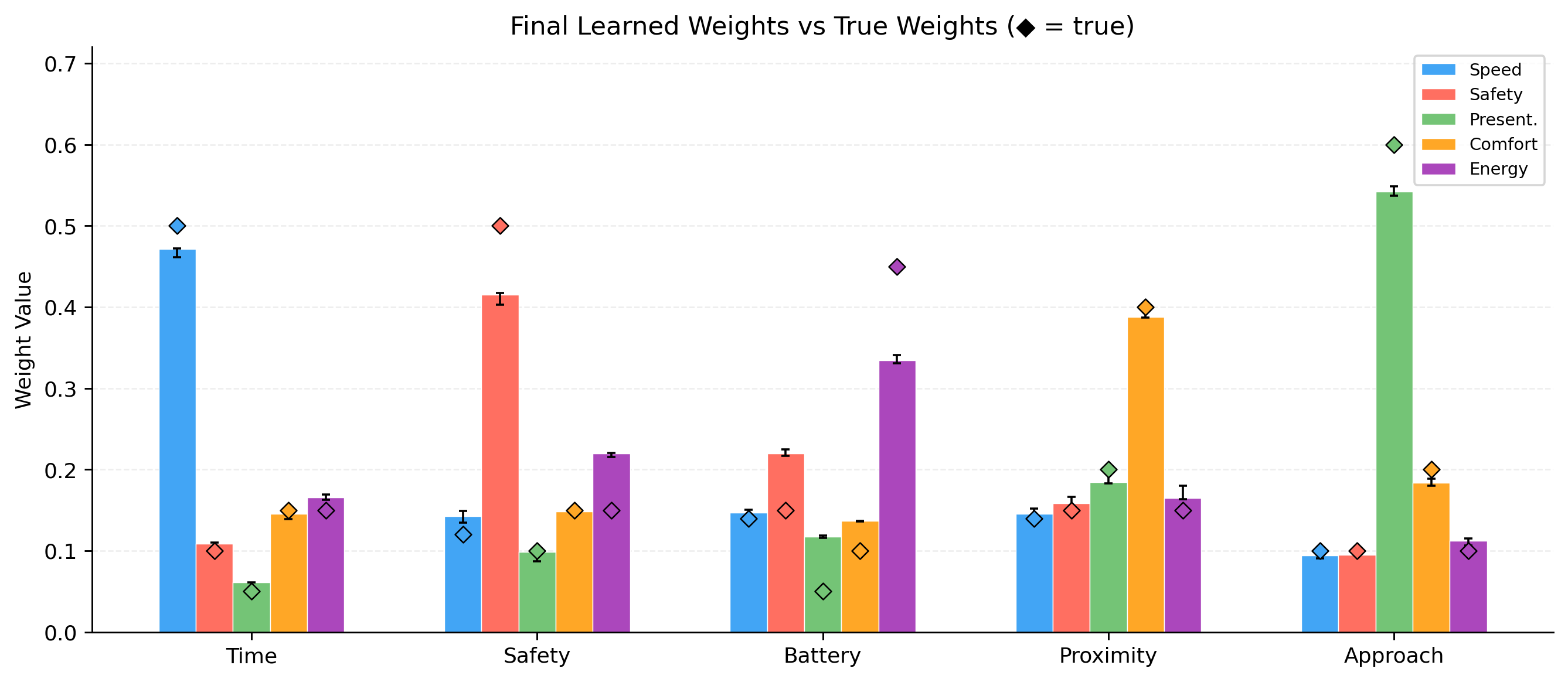}
    \caption{Final learned weights compared with the corresponding true weights. Diamonds indicate the true values.}
    \label{fig:final-weights}
\end{figure}

As shown in Figure~\ref{fig:final-weights}, the dominant preference dimension is correctly identified in all $25$ runs, covering all $5$ profiles across all $5$ seeds. This remains the operationally critical result: a misidentified dominant preference would produce systematically inappropriate plan selections, such as prioritising speed for a safety-first patient. Even for the presentation-focused profile, whose approach weight remains slightly underestimated, the approach dimension is correctly identified as dominant in all $5$ seeds.

\subsection{Plan Adaptation and Diversity}\label{subsubsec:plan_diversity}

The task planner's sensitivity to learned preferences is evaluated through the distribution of plan variants generated across episodes. A plan variant is defined by the structural decisions made by the high-level planner. For medication delivery, these include pharmacy selection, supply room selection, visit ordering, and approach side at the patient bedside. For meal preparation, variants are defined by meal type and approach side.

Across the evaluated runs, the planner produces multiple structural variants rather than a single repeated execution pattern. As shown in Table~\ref{tab:results}, all profiles achieve a 100\% task success rate and full convergence across the five runs. However, the number of distinct plan variants differs across profiles: the speed-oriented and safety-first profiles produce 4 variants, the presentation-focused and comfort-focused profiles produce 3 variants, and the energy-conscious profile produces 5 variants. This shows that the learned preference vector does not only modify the scalar cost optimized by the planner, but also influences the symbolic decisions selected during task execution.

\paragraph{Medication delivery.}
Medication delivery episodes complete successfully across all profiles and seeds, while also exhibiting profile-specific differences in route structure. Speed-oriented runs tend to select shorter routes, favouring the pharmacy and supply-room combination that minimises total travel time, even when this increases exposure to higher-risk corridors.

Safety-first runs instead favour lower-risk alternatives, accepting longer routes in exchange for reduced congestion or hazard exposure. This behaviour is consistent with the construction of the environment, where shorter route options are generally associated with higher risk exposure. Therefore, the safety-first profile is expected to deviate from the shortest available route when the shorter path crosses areas with higher risk or congestion.

In navigation-heavy medication delivery, energy-conscious runs show a pattern that is close to the speed-oriented profile. This is expected because travelled path length affects both execution time and energy consumption in a similar way: shorter and more direct routes tend to be both faster and less energy-consuming. As a result, the route choices that minimise travel time often also reduce battery usage, producing partial overlap between speed-oriented and energy-conscious planning behaviour.

Approach-side selection also contributes to profile-specific planning behaviour. In the experimental setup, left-side delivery is considered riskier than right-side delivery. Therefore, safety-first profiles tend to select right-side delivery when possible, because this reduces the risk component of the plan, even if it conflicts with the patient's preferred delivery side.

\paragraph{Meal preparation.}
Meal preparation episodes also complete successfully across all profiles and seeds. Speed-oriented profiles tend to favour faster meal options, such as sandwich preparation, because these minimise task duration. 

For the energy-conscious profile, meal selection makes the behaviour more distinguishable because it is strongly constrained by the energy cost of the available preparation options. In particular, the hot-meal is substantially more energy-consuming than simpler options such as sandwich or salad preparation. As a result, it is discarded by the energy-conscious profile, since its higher preparation cost is not compensated by a sufficient gain in the battery-related objective.

The plan-diversity results therefore complement the convergence results.

\subsection{Learning Continuous Target Parameters from Symbolic--Control Mismatches}

\begin{figure}[H]
\centering
\includegraphics[width=0.95\linewidth]{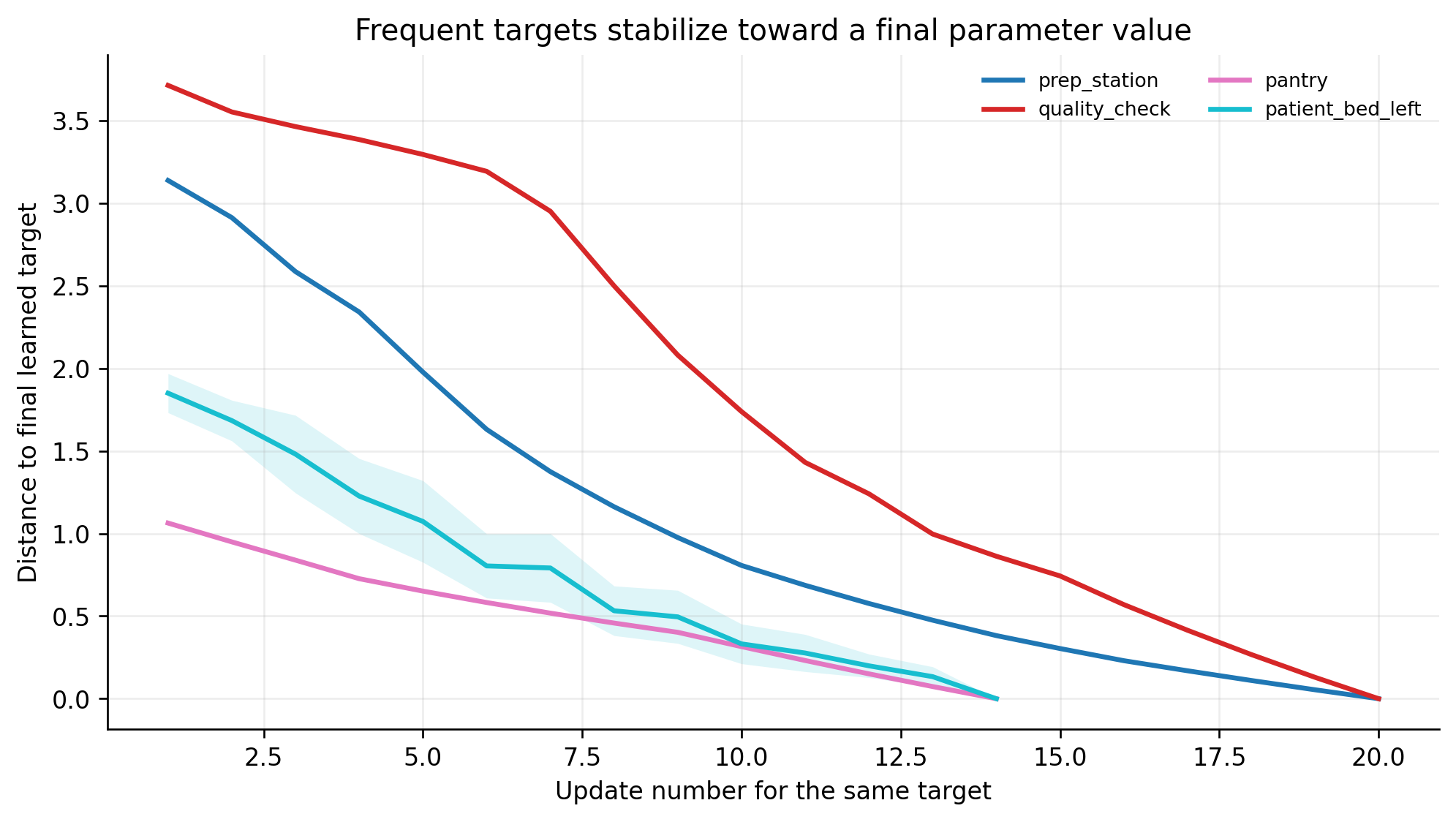}
\caption{
Distance between the current learned target parameter and the final learned value
as a function of the number of updates for the same target. }
\label{fig:target-convergence}
\end{figure}

This experiment evaluates the update mechanism connecting the symbolic planning layer
with the continuous MPC layer. In the simulated environment, we artificially perturb the
final continuous position associated with selected symbolic targets. As a consequence, the
nominal target used by the planner no longer coincides with the effective target that should
be used by the low-level controller. This creates a controlled mismatch between the symbolic
goal and the continuous state reached after MPC execution.

Each symbolic action activates a corresponding
MPC problem through an action-dependent parametrization. Therefore, when the observed
continuous outcome differs from the expected one, the mismatch can be used to update the
parameters of the MPC problem associated with that symbolic action. In this experiment,
the learned MPC-side parameters are the final target coordinates used by the controller. The purpose
is to verify whether repeated executions allow the system to move away from the incorrect
nominal target and stabilize toward a corrected continuous reference.

\begin{figure}[H]
\centering
\includegraphics[width=0.8\linewidth]{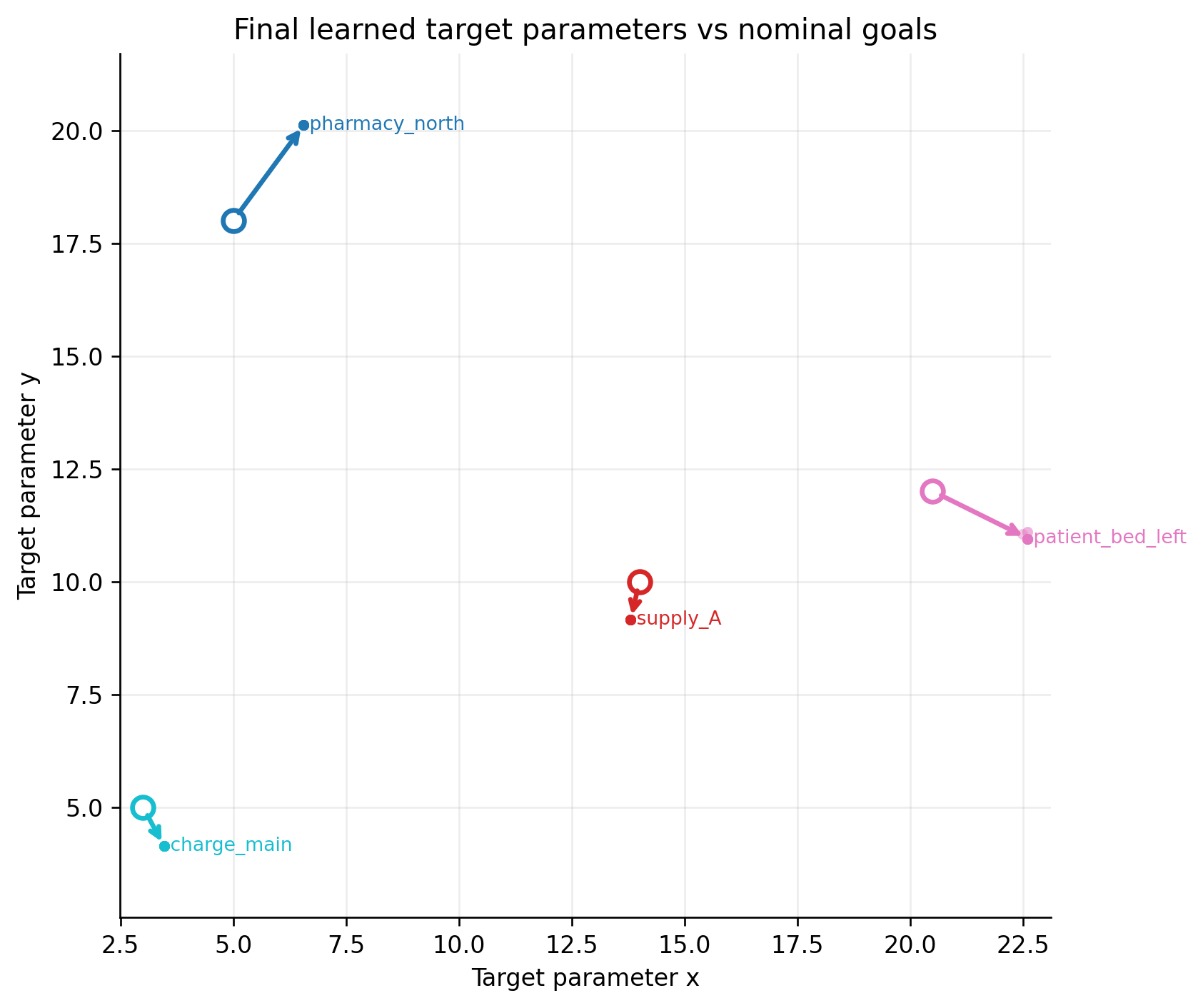}
\caption{
Comparison between nominal target positions and final learned target parameters.
Circles denote nominal goals, while arrows indicate the displacement toward the
learned continuous target. The learned parameters deviate from the nominal
positions in response to the artificial perturbation introduced in the environment.
}
\label{fig:target-map}
\end{figure}

Figure~\ref{fig:target-convergence} shows the distance between the current learned target
parameter and its final learned value as a function of the number of updates for the same
target. The expected behavior is observed: for frequently visited targets, this distance
decreases progressively and eventually reaches zero. This indicates that the system does not
continue to track the original nominal target indefinitely. Instead, it accumulates evidence
from repeated executions and converges toward a new operational target parameter.

The convergence speed differs across targets. Targets such as \texttt{pantry} and
\texttt{patient-bed-left} approach their final learned values after approximately
13--14 updates, whereas \texttt{prep-station} and especially \texttt{quality-check}
require more updates, stabilizing closer to 20 updates. This behavior is consistent with a
local data-driven update rule: the more often a target is visited, the more information is
available about the terminal mismatch, and the more accurately the corresponding MPC
target parameter can be corrected.

The experiment also supports the role of the fuzzification interface between the symbolic
and continuous layers. The continuous state $\hat{z}(t+1)$ is mapped into a symbolic or
fuzzy representation through membership functions $\mu_L(z)$. This makes it possible to
compare the state expected by the planner with the state estimated from the controller outcome.
Without this fuzzy representation, the terminal error would remain a purely geometric
continuous error. Through fuzzification, instead, it becomes an interpretable
symbolic--continuous mismatch signal that can be used to update the MPC parametrization.

\FloatBarrier

\section{Discussion}\label{sec:discussion}

\subsection{Real World Challenges}\label{subsec:hallenges}

An embodied AI system is defined as a system that purposefully exchanges both \emph{energy} (e.g., forces and torques) and \emph{information} (e.g., commands and sensor input) with its physical environment \cite{koditschek2021robotics, roy2021machine}.
In this section, we base our discussion on these two exchanges. We identify  challenges faced by purely RL-based solutions and analyze how the scheme in this work addresses the challenges.

The energy exchange requires the system to accurately model the dynamics of both itself and the physical environment, as the forces and torques involved must be predicted and controlled to achieve meaningful physical interaction. However, these dynamics often cannot be fully characterized until the actual deployment \cite{aljalbout2025reality}. Manufacturing inconsistencies, for example, may introduce mass imbalances in rotating components, producing unmodelled gyroscopic effects in a robot's motion, while unpredictable environmental conditions, such as varying weather, may alter the friction coefficient experienced by an autonomous vehicle~\cite{manas2025knowledge}.
Purely RL-based solutions often require pre-training in simulation before deployment in the real environment; consequently, unmodeled dynamics in simulation inevitably degrade performance and even compromise safety during actual deployment. This issue is commonly referred to as the sim-to-real gap~\cite{zhao2020sim, da2025survey}.
The proposed scheme in this paper addresses this issue by enabling the system to learn directly in its actual operating environment: the incorporation of a physical dynamics model ensures immediate functional deployment \rev{(evidenced by the 100\% task success rate across all episodes of all runs
reported in Sec.~\ref{subsubsec:plan_diversity} attained from the first episode of every run, before any preference or parameter adaptation)}, while RL-based tuning of the MPC controller and a continuously updated symbolic planner progressively improve performance by adapting to the specific characteristics of the physical environment \rev{(Figures~\ref{fig:target-convergence} and \ref{fig:target-map} show that the MPC parameters move from the nominal target to a corrected continuous reference over updates in response to environmental perturbation)}. 

The energy exchange also implies that the agent has a significant influence on the physical environment, making the control problem inherently safety-critical. In high-stakes domains such as autonomous vehicles and aerial drones, constraint violations can lead to collisions with severe consequences for bystanders and surrounding infrastructure~\cite{li2021differentiable, sormoli2024survey}. In human-robot interaction, safety beyond the physical domain and into the psychological domain is widely overlooked: poor alignment between robot behavior and user expectations can cause discomfort and erode the sense of trust~\cite{martinetti2021redefining}. 
The trial-and-error learning paradigm of purely RL-based solutions inevitably leads to violations of safety constraints during training, before the policy converges to safe and reliable behavior~\cite{lin2024almost}. Furthermore, purely RL-based solutions learn completely opaque policies, for which safety analysis and formal guarantees are not feasible.
To address this limitation, the proposed scheme embeds physical dynamics and constraints within the MPC formulation to guarantee physical safety, while addressing psychological alignment by generating user-interpretable task sequences via fuzzy logic mapping and incorporating user feedback to align the underlying cost functions in the optimization problems \rev{(This can be seen in Figure~\ref{fig:final-weights} where the true preferences are recovered by the learned weights)}.

Meanwhile, the information exchange in embodied AI systems requires efficient adaptation to novel tasks that arise during deployment~\cite{yuan2023hierarchical, jaquier2025transfer}. Consider a robotic care assistant performing scheduled tasks such as delivering meals and administering medication. When a new task arises, such as fetching supplies from a new storage room, the system should not need to re-learn how to operate in its environment from scratch. 
However, purely RL-based solutions typically require end-to-end retraining when faced with new tasks, as the learned policy is specific to the dynamics and reward structure of the original task~\cite{andreas2017modular, zhu2023transfer}.
The scheme in this paper addresses this by exploiting the fact that both the physical dynamics model in the MPC and the discrete action vocabulary in the symbolic planner are largely invariant across tasks; only the cost functions and constraints need to be redefined for the new task. The adaptation to a novel task is therefore reduced to specifying a new planning problem within the existing framework, rather than learning from scratch. \rev{For example, the medication delivery and meal preparation tasks in Sec.~\ref{sec:experiments} are structurally distinct planning problems, yet they share the dynamics $f$, the environment's location map, and the fuzzy membership functions grounding the location predicates. Since the same five-dimensional cost-feature vector scores both tasks, a single preference vector is learned across episodes alternating between them.}


The aforementioned embodied AI challenges, namely, imperfect modeling, safety-critical control, and rapidly emerging tasks, are fundamental to many real-world applications. To show the practical relevance of this paper, we now discuss several representative domains.

\vspace{1em}

\subsection{Other Applications}\label{subsec:implementation}

\paragraph{Autonomous Driving.}
Autonomous driving vehicles navigate urban traffic environments with minimal or no human intervention, which requires the ability to perceive multi-modal sensory data, control vehicle dynamics, and respond to unforeseen situations~\cite{zhao2025survey}. For purely RL-based methods, inevitable modeling errors in simulated training and the variance of real-world environment compromise the safety of nearby vehicles and pedestrians. This calls for in-environment learning while respecting safety constraints, which this paper addresses by combining model-based control and adaptive updates.

\paragraph{Drones for Remote Sensing and Agriculture.}
The development of unmanned aerial vehicle has enabled drones to perform various remote sensing and agricultural tasks, such as surveying and mapping of land, crop health monitoring, and precision spraying of pesticides~\cite{bilal2025consumer, kartal2025comprehensive}. For these complicated tasks with potentially sparse reward signals, the end-to-end RL policies are difficult to train and lack the interpretability. The proposed scheme decomposes these tasks into tractable and transferable subtasks, and uses interpretable symbolic planning to enable formal safety analysis.

\paragraph{Automated Manufacturing.}
Automated manufacturing systems perform industrial tasks such as assembly, welding, and material handling, sometimes in collaboration with human operators~\cite{mattera2025optimal, ngwu2025reinforcement}. When significantly different production tasks arise, policies tailored to the previous task cannot generalize. In contrast, the proposed framework enables efficient task reconfiguration by leveraging shared dynamics models across tasks \rev{(as in the two-task instantiation of Sec.~\ref{sec:experiments})}. Furthermore, preference learning via updating the cost functions of the planner and the MPC improves alignment at the interface between automated manufacturing and human operators.

\paragraph{Domestic Robots.}
\rev{Domestic robots are designed to perform household tasks, such as daily chores and healthcare support~\cite{soni2024advancing}, closely related to the health-service case study of Sec.~\ref{sec:experiments}}. A capable agent ensures safety, fosters user trust, and adapts efficiently to changes in the home environment~\cite{newman2023towards}. These requirements are difficult to satisfy with purely RL-based solutions, which lack safety guarantees and produce opaque, task-specific policies. In contrast, the scheme proposed in this work formally enforces safety constraints, generates user-interpretable action sequences via fuzzy logic mapping, and adapts to new tasks by reusing learned system dynamics and user preferences.

\subsection{Limitations and Future Work}\label{subsec:limitations}

\textcolor{blue}{Realizing the proposed scheme in practice raises several open challenges. RLMPC
learns from real-world interaction, which is expensive to obtain in many
applications and particularly so in the safety-critical ones that motivate the
approach. Practical deployment therefore hinges on learning efficiently from
small offline or off-policy batches, including data logged under earlier
parameterizations, under teleoperation, or during unrelated tasks, which in turn
demands offline and off-policy algorithms developed specifically for RLMPC. At
present, RLMPC schemes largely inherit general-purpose model-free algorithms such
as PPO and SAC, and consequently neither exploit the structure of the underlying
optimal control problem nor account for the ways in which differentiating through
an MPC scheme differs from backpropagating through a neural network. Recent work has begun to address this \cite{anand2023painless}, but
further algorithmic development is needed before such methods scale to complex
real-world tasks. Beyond the learning algorithm itself, the scheme commits to
point estimates of its learned parameters and is therefore equally confident in
regimes it has explored extensively and in those it has barely visited;
quantifying uncertainty in the resulting decisions, and incorporating risk
measures directly into the underlying optimization, remain to be addressed. Finally, deployment would require a certifiable
framework that verifies after each learning step that the adapted parameters
preserve the constraint satisfaction guarantees of the MPC scheme, so that
improvements in empirical performance cannot silently erode the safety properties
that motivate the use of MPC in the first place.}

\section{Conclusion}\label{sec:conclusion}

This paper addressed a fundamental tension in deploying learning-enabled robots in safety-critical environments: black-box reinforcement learning can discover effective policies but offers no guarantees on safety, interpretability, or predictability --- properties that are non-negotiable in domains such as hospital care. We presented a multilevel control architecture that achieves the adaptive personalisation of RL while maintaining the deterministic, interpretable structure required for autonomous systems operating alongside vulnerable populations.

The bilevel structure couples an outer loop that performs classical planning - an interpretable and reliable algorithmic toolkit for choosing discrete combinatorial action sequences that achieve a goal. Each action defines a parametrization of an optimal control problem. An MPC controller defines continuous actuator operation of an autonomous system that appropriately incorporates known physics engineering. Learning is performed between the layers through real-time reward signals and discrete-continuous fuzzy map discrepancies. These propagate, through MPC sensitivities, as gradient information, closing the entire system loop. 

We validated the system in a hospital service robot domain with two structurally distinct task types - medication delivery and meal preparation - across five patient profiles spanning different regions of a  5-dimensional preference space. Consistently successful task completion and adaptation to patient preferences was observed in the simulation validation.

The increasing interest of deploying AI in real world scenarios, from autonomous vehicles and drones to robotics in manufacturing and care, presents the necessity of enforcing reliability through first principles domain knowledge. Increasingly, there is a complex layered taxonomy of tasks and behaviors, and interpretability is essential for fostering trust and robust performance. The current default of scaling Deep RL for operating these autonomous systems presents significant risks to safety as well as exhibiting limitations in agreeable deployment in human-sensitive settings due to its opaque nature. In addition, the requirement of increasingly deeper models and training computation presents unsustainable energy and economic costs, as far as the desired scale of deployment of such systems. Incorporation of crisp interpretable planning together with physics-driven control is, in the authors' view, essential towards comprehensive reliability and satisfying performance of autonomous AI systems.

This paper presents a first attempt at simultaneously integrating classical planning, which encodes crisp goals in a solid logical framework, through to MPC, the most state of the art the technology for fast real time control, while incorporating classic RL tools for learning the rewards from the environment. We defined the entire end-to-end Scheduling to Planning to Control, including both the operation and the co-learning across layers. We hope the work opens up a productive research program towards the management and control of autonomous agents. 

The potential for future work extending the presented framework is vast - mathematical theory proving guarantees of approximate optimality, stability, sample complexity and other criteria, validation in physical settings, synergistic bilevel algorithmic development, expansion across use cases. We hope that this work triggers a powerful and effective research and development program towards reliable and effective autonomous systems broadly.

\section*{Data Availability}
All data used in this manuscript is synthetic and publicly available at \\\texttt{https://github.com/alexdifre/hospital-robot-system}.

\bibliographystyle{unsrt}
\bibliography{refs.bib}
\end{document}